\documentclass[11pt]{article}

\usepackage{graphicx}
\usepackage{subcaption}
\usepackage{longtable}
\usepackage{wrapfig}
\usepackage{rotating}
\usepackage[normalem]{ulem}
\usepackage{amsmath}
\usepackage{amssymb}
\usepackage{mathtools}
\usepackage{capt-of}
\usepackage{hyperref}
\usepackage{fullpage}
\usepackage{xcolor}
\usepackage{booktabs} % For toprule, midrule and bottomrule
\usepackage{bm}
\usepackage{tocloft}

\DeclareMathOperator{\real}{Re}
\DeclareMathOperator{\imag}{Im}
\DeclareMathOperator{\sign}{sign}
\DeclareMathOperator{\midint}{mid}
\newcommand{\conj}[1]{\overline{#1}}
\newcommand{\normv}{\mathrm{v}}
\newcommand{\Space}{\mathcal{X}}
\newcommand{\inter}[1]{\bm{#1}}

\newtheorem{definition}{Definition}[section]
\usepackage{amsthm}
\newtheorem{remark}[definition]{Remark}
\newtheorem{theorem}[definition]{Theorem}
\newtheorem{conjecture}[definition]{Conjecture}
\newtheorem{proposition}[definition]{Proposition}
\newtheorem{lemma}[definition]{Lemma}

\newcommand{\codenumber}[1]{#1} % Use this when publishing
\usepackage[backend=biber,url=false,eprint=false,maxbibnames=99]{biblatex}
\hypersetup{
  bookmarksopen=true,
  bookmarksdepth=subsubsection
}
\graphicspath{{figures/}}

\title{Self-Similar Singular Solutions to the Nonlinear Schrödinger and the Complex Ginzburg-Landau Equations}
\author{Joel Dahne, Jordi-Llu\'is Figueras}
\date{}

\begin{document}

\maketitle

\begin{abstract}
  We prove the existence of several branches of radial self-similar
  singular solutions for the mass supercritical complex
  Ginzburg--Landau equation, including the nonlinear Schrödinger
  equation as a special case. We are also able to control the number
  of monotone intervals of the profiles associated with these
  self-similar solutions. In particular, we prove the existence of
  monotone radial self-similar singular solutions for the
  three-dimensional cubic nonlinear Schrödinger equation. The paper
  combines sharp analytic bounds on the self-similar profile at
  infinity with computer-assisted bounds around zero and their
  matching at an intermediate point.
\end{abstract}

\tableofcontents

\section{Introduction}
\label{sec:introduction}
We address the problem of self-similar blowup for the complex
Ginzburg--Landau (CGL) equation,
\begin{equation}
  \label{eq:CGL}
  i \frac{\partial u}{\partial t} + (1 - i\epsilon)\Delta u + (1 + i\delta)|u|^{2\sigma}u = 0
  \quad \text{ in }\mathbb{R}^{d} \times (0, T).
\end{equation}
Here the parameters \(\epsilon\) and \(\delta\) are nonnegative real
numbers and \(\sigma\) is positive. Our interest lies in the existence
of branches of such self-similar solutions when varying the parameter
\(\epsilon\) (fixing \(\delta = 0\)). These branches include the point
\(\epsilon = \delta = 0\), corresponding to the focusing nonlinear
Schrödinger (NLS) equation. As a result, we also prove self-similar
blowup for the NLS equation.

\subsection{The NLS case}
\label{sec:nls-case}
The study of self-similar blowup for the NLS equation has a long
history; let us therefore start by describing our result in the
context of this special case. The NLS equation is one of the
fundamental equations in the study of nonlinear waves and the
questions of well-posedness and finite-time blowup are classical
problems. The existence of blowup solutions and sufficient conditions
for their occurrence were established by Glassey in
1977~\cite{Glassey77}. Since then, many researchers have sought to
characterize the behavior of these solutions.

The nature and existence of the blowup depend on the value of the
criticality constant \(s_{c} := \frac{d}{2} - \frac{1}{\sigma}\). In
the energy-subcritical case, \(s_{c} < 1\), the Cauchy problem is
locally well-posed in the energy space \(H^{1}\). In the
mass-subcritical regime, \(s_{c} < 0\), there is global existence for
\(H^{1}\) initial data, and hence no blowup occurs,
see~\cite{Cazenave88, Cazenave2003}. For both the mass-critical
(\(s_{c} = 0\)) and mass-supercritical but energy-subcritical
(\(0 < s_{c} < 1\)) cases, a large number of numerical studies in the
1980s and 1990s indicated the possibility of radial self-similar
blowup. An excellent reference is the comprehensive monograph by Sulem
and Sulem~\cite{Sulem1999}.

For the critical case, the numerics indicated that these solutions
were only asymptotically self-similar. The questions of existence,
rate of convergence and stability for these solutions were finally
settled in a remarkable series of papers by Merle and
Raphaël~\cite{Merle2005a, Merle2004a, Merle2006, Merle2004b} (see
also~\cite{Perelman2001}), which provide a detailed description of the
phenomena.

Our main interest in this paper is the supercritical case. For
slightly supercritical parameters, \(0 < s_{c} \ll 1\), the existence
of self-similar blowup has been established,
see~\cite{Rottschafer2002, Rottschafer2003, Merle2010, Merle2014,
  Bahri2021}. For larger values of \(s_{c}\), which in particular
include the three-dimensional cubic NLS equation
(\(s_{c} = \frac{1}{2}\)), the question has remained unsolved for a
long time and the conjectured behavior is primarily based on numerical
simulations. However, while this manuscript was in preparation,
Donninger and Schörkhuber~\cite{Donninger2024} gave a remarkable proof
of self-similar blowup for the three-dimensional cubic NLS equation,
which we discuss further below.

Before stating our results, let us describe the existing numerics for
the supercritical case in more detail. For the radial
self-similar solutions that we are interested in,
Zakharov~\cite{Zakharov84} proposes the profile
\begin{equation}
  \label{eq:blowup-profile}
  u(x, t) = \dfrac{1}{(2\kappa(T - t))^{\frac{1}{2}\left(\frac{1}{\sigma} + i \frac{\omega}{\kappa}\right)}}
  Q\left(\dfrac{|x|}{(2\kappa(T - t))^{\frac{1}{2}}}\right),
\end{equation}
where \(\kappa, \omega \in \mathbb{R}\) and
\(Q \colon [0, \infty) \to \mathbb{C}\) is a profile with asymptotic behavior
\begin{equation*}
  Q(\xi) \sim \xi^{-\frac{1}{\sigma} - i \frac{\omega}{\kappa}} \text{ as } \xi \to \infty.
\end{equation*}
Inserting the ansatz~\eqref{eq:blowup-profile} into
Equation~\eqref{eq:CGL} gives the singular ODE
\begin{align}
  \label{eq:Q}
  0 &= (1 - i\epsilon)\left(Q'' + \frac{d - 1}{\xi}Q'\right) + i\kappa\xi Q'
  + i \frac{\kappa}{\sigma}Q - \omega Q + (1 + i\delta)|Q|^{2\sigma}Q,\\
  Q'(0) &= 0,\nonumber\\
  Q(\xi) &\sim \xi^{-\frac{1}{\sigma} - i \frac{\omega}{\kappa}} \text{ as } \xi \to \infty.\nonumber
\end{align}
For arbitrary values of \(\kappa\) and \(\omega\), the equation does
not have a solution; the unknowns are hence
\((Q(\xi), \kappa, \omega)\). For \(\lambda > 0\), there is a scaling
symmetry taking a solution \((Q(\xi), \kappa, \omega)\) to another
solution
\((\lambda^{1/\sigma}Q(\lambda\xi), \lambda^{2}\kappa,
\lambda^{2}\omega)\). This allows us to fix \(\omega = 1\) and look
for solutions \((Q(\xi), \kappa)\). Numerical simulations by Budd,
Chen and Russell~\cite{Budd99} indicate the existence of a countable
set of \(\kappa\) values for which solutions exist. These come in the
form of ``multi-bump'' solutions, which (in some parameter regimes) seem
to be characterized by the number of monotone intervals (the monotone
index) of \(|Q|\). Of these, only the first solution, the monotone
one, appears to be (nonlinearly) stable. Sulem and
Sulem~\cite{Sulem1999} refer to solutions of~\eqref{eq:Q} for which
\(|Q|\) is monotonically decreasing as ``admissible solutions''. Based
on these numerical results, the following conjecture can be
formulated:

\begin{conjecture}[NLS]
  The mass-supercritical, but energy-subcritical, NLS equation has a
  countable number of nontrivial, radial, self-similar singular
  solutions associated with solutions of Equation~\eqref{eq:Q}. The
  first of these solutions has a monotone profile \(|Q|\) and is
  stable, while the remaining ones are unstable.
\end{conjecture}

Our first contribution in this paper is to provide a partial answer to
this conjecture: we prove the existence of a finite number of such
solutions for two different sets of parameters. Specifically, we adopt
the choices considered in the work by Plecháč and
Šverák~\cite{Plech2001}, which we discuss further below, corresponding
to the following two cases:
\begin{description}
\item[Case I] \(d = 1\) and \(\sigma = 2.3\), corresponding to
  \(s_{c} \approx 0.065\).
\item[Case II] \(d = 3\) and \(\sigma = 1\), corresponding to
  \(s_{c} = \frac{1}{2}\). This is the three-dimensional cubic case.
\end{description}
More precisely, we have the following two results:

\begin{theorem}[\textbf{Case I}: NLS equation]\label{thm:NLS-intro}
  There exist 8 nontrivial radial self-similar singular solutions
  \(u\) of the form~\eqref{eq:blowup-profile} with
  \(Q \in C^{\infty}([0, \infty)) \cap L^3([0, \infty))\)
  satisfying Equation~\eqref{eq:CGL} for \(d = 1\) and
  \(\sigma = 2.3\) (and \(\epsilon = \delta = 0\)). For the first of these
  solutions, \(|Q|\) is monotone, while the other solutions are
  labeled by the number of monotone intervals of \(|Q|\).
  \end{theorem}

\begin{theorem}[\textbf{Case II}: NLS equation]
  There exist 2 nontrivial radial self-similar singular solutions
  \(u\) of the form~\eqref{eq:blowup-profile} with
  \(Q \in C^{\infty}([0, \infty)) \cap L^2([0, \infty))\)
  satisfying Equation~\eqref{eq:CGL} for \(d = 3\) and
  \(\sigma = 1\) (and \(\epsilon = \delta = 0\)). For the first solution,
  \(|Q|\) is monotone, while for the other solution, \(|Q|\) has two
  intervals of monotonicity.
\end{theorem}

These two theorems are expanded upon in Section~\ref{sec:nls}.

\begin{remark}
  In the two theorems above, the solutions are characterized by their
  monotone index, similar to the results of Budd, Chen and
  Russell~\cite{Budd99}. In general, however, this does not appear to
  be the case. For example, the third solution for the
  three-dimensional cubic NLS equation appears to be monotone, see
  e.g.~\cite[Figure 3.9]{Plech2001}.
\end{remark}

\subsection{The CGL case}
\label{sec:cgl-case}

Next, we investigate these solutions as one moves into the regime of
the CGL equation. The CGL equation has a long history in physics,
serving as a model equation for a vast variety of phenomena, from
nonlinear waves and second-order phase transitions, to
superconductivity and superfluidity. We refer the reader to the review
paper~\cite{Aranson2002} for many more details. Our interest in this
problem comes from the work by Plecháč and Šverák~\cite{Plech2001}.
One of their motivations for studying the CGL equation is the analogy
between the three-dimensional cubic CGL equation and the
Navier--Stokes equation. The two equations share the same scaling
properties and energy identity. The existence and partial regularity
theory for the two equations are also similar~\cite{Yan1999}. In fact,
a number of works in the 1990s studied the CGL equation as a model for
the Navier--Stokes equation, see
e.g.~\cite{Bartuccelli1990,Doering1994}.

Using a combination of asymptotic analysis and numerical simulations,
Plecháč and Šverák provided strong evidence that the self-similar
solutions of the NLS equation persist for small values of
\(\epsilon\), forming branches of solutions as \(\epsilon\) is varied
(see Figure~\ref{fig:branches-numerical}). Some of these solutions
were further studied numerically by Budd, Rottschäfer and
Williams~\cite{Budd2005}. In the case of asymptotically small
\(\kappa\), self-similar solutions of this type were proved to exist
by Rottschäfer~\cite{Rottschafer2008,Rottshafer2013}. The CGL equation
is also known to support other, non-self-similar, types of
blowup~\cite{Masmoudi2008, Cazenave2014, Duong2022, Duong2023,
  Duong2023FlatBlowup}. For global existence results,
see~\cite{Shimotsuma2016, Correia2019, Correia2021}.

Our second contribution in this paper is to prove the existence of
some of the branches of self-similar singular solutions studied by
Plecháč and Šverák. For the one-dimensional case, a summary of the
results is given in the following theorem.

\begin{theorem}[\textbf{Case I}: CGL equation]
  The solutions in Theorem~\ref{thm:NLS-intro} extend into branches of
  solutions for the CGL equation. Along these branches, the number of
  monotone intervals of the profile \(|Q|\) remains constant.
\end{theorem}

We refer to Section~\ref{sec:branches} and
Theorem~\ref{thm:existence-cgl-case-1} for a more precise statement.
In \textbf{Case II}, the problem is computationally more challenging;
see Theorem~\ref{thm:existence-cgl-case-2} for the precise results we
prove in this case. We also give some extended results for
\textbf{Case II} in Section~\ref{sec:pointwise-results-cgl}, where we
prove pointwise existence for a larger parameter range, but do not
prove that they necessarily connect into branches.

\begin{remark}
  Previous numerical studies of these branches show that when plotted
  in the \(\epsilon\)-\(\kappa\) space, they have a single turning
  point in \(\epsilon\), but otherwise appear monotone in these
  parameters. Our results indicate that this does not, in fact,
  capture the full behavior for all branches. In particular, for the
  three-dimensional cubic case, our numerical results indicate that
  the third branch has two turning points in \(\kappa\), and the fifth
  branch has multiple turning points in both \(\epsilon\) and
  \(\kappa\) (see Figure~\ref{fig:branches-d3-numerical}). This
  behavior is not seen in~\cite{Plech2001} because the continuation of
  the branches stops right before reaching the first turning point in
  \(\kappa\) (a likely explanation being that they did not expect this
  behavior, and their numerical scheme was therefore not adapted to
  handle it). In~\cite{Budd2005}, only the first two branches are
  studied, which do not show this behavior. Our rigorous results
  presented here only cover the first turning point in \(\kappa\) for
  the third branch. Further numerical exploration of these phenomena
  would be desirable. For example, it is at this point not clear if
  the fifth branch eventually turns back to start approaching the
  origin.
\end{remark}

\begin{remark}
  Both Plecháč and Šverák~\cite{Plech2001} as well as Budd,
  Rottschäfer and Williams~\cite{Budd2005} study the stability of the
  solutions along the branches. For the three-dimensional cubic case,
  their results indicate that the top half of the first branch and the
  bottom half of the second branch are stable, and the other parts of
  the branches are unstable. The parallels between the CGL equation
  and the Navier--Stokes equation make the existence of unstable
  self-similar blowup particularly interesting, as this is one of the
  possible scenarios regarding the question of well-posedness for
  Navier--Stokes. This is, for example, mentioned as a possibility
  in~\cite{Wang2025}, where they numerically isolate unstable
  self-similar solutions for the Euler equation. However, there may be
  constraints on the extent to which sufficiently unstable solutions
  can exist~\cite{Constantin2026}.
\end{remark}

\subsection{Proof strategy}
\label{sec:proof-strategy}

Our strategy, discussed in more detail in
Section~\ref{sec:shooting-method}, is based on a shooting method
inspired by the work of Plecháč and Šverák. The main idea is to
construct two solutions, \(Q_{0}\) and \(Q_{\infty}\), to
Equation~\eqref{eq:Q}, with \(Q_{0}\) satisfying the boundary
condition at zero and \(Q_{\infty}\) the condition at infinity. We
then need to show that with the right choice of parameters, these two
solutions match at an intermediate point, forming a global solution
satisfying both boundary conditions. The approach shares similarities
with that used by Bahri, Martel and Raphaël~\cite{Bahri2021} and
Rottschäfer~\cite{Rottshafer2013}, an important difference being that
we cannot rely on smallness of any of the parameters to prove the
required matching condition. Instead, we make use of computer-assisted
methods to establish the necessary conditions. This requires obtaining
effective, not only asymptotic, bounds for \(Q_{0}\) and
\(Q_{\infty}\).

One of the main challenges is that the equation has an irregular
singular point at \(\infty\). To handle this, we take inspiration from
the asymptotic analysis done by Plecháč and Šverák, computing
asymptotic expansions for the solution \(Q_{\infty}\). However,
obtaining effective bounds from these expansions requires bounds on
all the remainder terms involved. Acquiring these effective bounds is
one of the most challenging aspects of the proof. This difficulty is
amplified by the fact that we need bounds that are uniform near
\(\epsilon = 0\). Away from zero, one can use simpler bounds, but they
deteriorate as \(\epsilon\) approaches zero. Treating only the case
\(\epsilon = 0\) would also simplify the work, since several
parameters become exactly zero.

For \(Q_{0}\), the approach is split into two parts. Near the origin,
the singularity at \(0\) is handled using a Taylor expansion with
explicit bounds for the remainder. Away from the origin, the solution
is then continued to sufficiently large \(\xi\) values using a
rigorous numerical ODE solver. At this stage, we are dealing with an
ODE on a finite interval bounded away from \(0\), meaning that there
are no singularities in the path. For this setting, there is a large
amount of work available on the rigorous computation of solutions. We
make use of a solver from the CAPD library~\cite{Kapela2021} based on
taking successive steps using Taylor expansions with rigorous error
bounds. Another common approach is the use of spectral methods, see
e.g.~\cite{Lessard2014, Hungria2015, vandenBerg2021, Breden2026}.

The primary need for computer assistance arises within the
intermediate regime that bridges the local Taylor expansion at the
origin with the asymptotic series at infinity. This part is inherently
nonperturbative, and obtaining effective bounds without the help of
the computer seems difficult. To connect the intermediate regime with
the two endpoints, we must also explicitly evaluate the effective
bounds coming from the respective expansions. This evaluation is,
however, fairly standard, requiring only the evaluation of a couple of
special functions and a (large) number of explicit constants.

While most of the mathematical difficulty lies in the derivation of
the effective bounds for the asymptotic expansion of \(Q_{\infty}\),
the technical difficulties for a complete implementation of the
computer-assisted proof are also substantial. In particular, tracking
the branches for the CGL equation requires a significant amount of
work. For example, the highly non-uniform computational difficulty
along the branches requires the use of adaptive bisection approaches
that scale across hundreds of cores.

In a recent independent work, Donninger and
Schörkhuber~\cite{Donninger2024} proved the existence of a
self-similar singular solution to the three-dimensional cubic NLS
equation, corresponding to the first profile of our second theorem.
Their approach similarly relies on a combination of quantitative
analysis and computer assistance, but operates under a different
paradigm. They compactify the domain and formulate a fixed-point
problem centered around an explicit polynomial approximation. The
fixed-point setup is complicated by the fact that the associated
linear operator itself depends on the unknown \(\kappa\). To handle
this, they introduce a specialized functional that vanishes at the
correct parameter value. The existence of the fixed point is then
proved by computing several bounds related to the linear operator and
the functional. The explicit polynomial provides the global control
required to connect the two endpoints, serving a role analogous to the
rigorous ODE solver used for the intermediate regime in our approach.

For an introduction to the type of computer assistance and rigorous
interval arithmetic that is used in this paper, see, for example, the
book~\cite{Tucker2011} for a general introduction and~\cite{Plumbook,
  GomezSerranoSEMA} for more PDE-oriented discussions. These
techniques have yielded significant results; see, for example,
classical results such as the existence of universality in the
quadratic map (period-doubling cascade)~\cite{Lanford}, the
double-bubble conjecture~\cite{Hass_1995}, the existence of strange
attractors in the Lorenz system~\cite{Tucker_1999} or the solution to
the Kepler conjecture~\cite{Hales2005}. Moreover, computer-assisted
proofs have been successful in proving the existence of homoclinic and
heteroclinic orbits in different systems~\cite{Zgliczynski98,
  Wilczak2020, Henot2023, Mireles2024}. Other results for PDEs appear
in~\cite{Fefferman86, Fefferman1996, Arioli2021}, and in Hamiltonian
systems in~\cite{Rana1991, Celletti1995, Figueras2016, Danesi2023,
  Capinski2023}.

\subsection{Organization of the paper}
In Section~\ref{sec:shooting-method}, we present the general approach
we use for proving the existence of a solution \(Q\) to the boundary
value problem~\eqref{eq:Q}. As discussed above, it is based on
matching two solutions, \(Q_{0}\) and \(Q_{\infty}\), satisfying the
boundary conditions at \(0\) and \(\infty\), respectively. Then, in
Sections~\ref{sec:nls} and~\ref{sec:branches}, we present in more
detail the results obtained for both the NLS and CGL equations.
Sections~\ref{sec:proof-nls} and~\ref{sec:proof-cgl} contain a
detailed explanation of how the results are obtained for these two
equations. Section~\ref{sec:solution-infinity} discusses the
asymptotic expansion used for controlling \(Q_{\infty}\), and in
particular how to obtain effective bounds on the error. In
Section~\ref{sec:solution-zero}, we deal with \(Q_{0}\), both the
expansion at zero and some of the details of the rigorous ODE solver.
Section~\ref{sec:pointwise-results-cgl} presents additional results on
the existence of solutions for the CGL equation for \textbf{Case II}.
Section~\ref{sec:implementation-details} discusses the implementation
of the computer-assisted proofs. Finally, in
Appendix~\ref{sec:function-bounds-proofs}, we give some auxiliary
lemmas used for the proofs in Section~\ref{sec:solution-infinity}.

The full code on which the computer-assisted parts of the proofs are
based and notebooks presenting the results and figures are
available in the repository~\cite{CGL.jl}.

\begin{remark}
  In this paper, tight intervals will be represented with
  subsuperscript notation (e.g.\
  \([1.147721, 1.147734] = 1.1477_{21}^{34}\)), while wide intervals
  will be represented with both-endpoint notation (e.g.\
  \([-2.4, 10.2]\)).
\end{remark}

\section{A rigorous shooting method}
\label{sec:shooting-method}
The shooting method is a classical numerical method for solving
boundary value problems by reducing them to a combination of an
initial value problem and a zero finding problem: one searches for
initial conditions on one of the boundaries for which the solution
also satisfies the boundary condition at the other boundary. We give a
rigorous version of this approach that, under the right conditions,
allows us to prove the existence of a solution to the boundary value
problem. Our version mirrors the approach used by Plecháč and
Šverák~\cite{Plech2001} when numerically computing the branches for
the CGL equation, the main difference being the use of rigorous error
bounds and a slight variation in how the solution at infinity is
handled.

Our starting point is to construct two partial solutions to
Equation~\eqref{eq:Q}, the first one satisfying the boundary condition
at zero and the second one satisfying the boundary condition at
infinity. Both of these functions depend on several parameters and
the goal is to find parameters such that these solutions agree at an
intermediate point, forming a global solution satisfying all required
boundary conditions.

We denote the solution starting from infinity by \(Q_{\infty}(\xi)\). It satisfies
\begin{equation}
  \label{eq:problem-infinity}
  \begin{split}
    (1 - i\epsilon)\left(Q_{\infty}'' + \frac{d - 1}{\xi}Q_{\infty}'\right) + i\kappa\xi Q_{\infty}'
    + i \frac{\kappa}{\sigma}Q_{\infty} - \omega Q_{\infty} + (1 + i\delta)|Q_{\infty}|^{2\sigma}Q_{\infty} &= 0,\\
    Q_{\infty}(\xi) &\sim \xi^{-\frac{1}{\sigma} - i\frac{\omega}{\kappa}} \quad\text{as } \xi \to \infty.
  \end{split}
\end{equation}
The boundary condition at infinity is, however, not enough to
determine the solution uniquely. As discussed by Plecháč and
Šverák~\cite{Plech2001}, one expects a one dimensional (complex)
manifold of solutions. We will (locally) parameterize this manifold by
\(\gamma \in \mathbb{C}\) (the details of this parameterization are
discussed in Section~\ref{sec:solution-infinity}). For our purposes,
the dependence on the parameters \(\kappa\) and \(\epsilon\) is also
important; to emphasize this dependence, we use the notation
\(Q_{\infty}(\xi) = Q_{\infty}(\gamma, \kappa, \epsilon; \xi)\).

We denote the solution starting from zero by \(Q_{0}(\xi)\). It is the
unique solution to the following initial value problem
\begin{equation}
  \label{eq:problem-zero}
  \begin{split}
    (1 - i\epsilon)\left(Q_{0}'' + \frac{d - 1}{\xi}Q_{0}'\right) + i\kappa\xi Q_{0}'
    + i \frac{\kappa}{\sigma}Q_{0} - \omega Q_{0} + (1 + i\delta)|Q_{0}|^{2\sigma}Q_{0} &= 0,\\
    Q_{0}(0) &= \mu,\\
    Q_{0}'(0) &= 0.
  \end{split}
\end{equation}
Again, we emphasize the dependence on the parameters \(\mu\),
\(\kappa\) and \(\epsilon\) with the notation
\(Q_{0}(\xi) = Q_{0}(\mu, \kappa, \epsilon; \xi)\). Solutions could be
defined on \([0, \infty)\), however, generically they have large
oscillations as \(\xi \to \infty\), and hence do not have finite
energy.

Since our ODE is of second order, the required matching condition is
that the values and derivatives of \(Q_{0}\) and
\(Q_{\infty}\) should match at some point \(0 < \xi_{1} < \infty\). If
indeed we have parameters \(\mu\), \(\gamma\), \(\kappa\) and
\(\epsilon\) such that
\begin{equation}
  \label{eq:Q-match}
  Q_{0}(\mu, \kappa, \epsilon; \xi_{1}) = Q_{\infty}(\gamma, \kappa, \epsilon; \xi_{1})
  \text{ and }
  Q_{0}'(\mu, \kappa, \epsilon; \xi_{1}) = Q_{\infty}'(\gamma, \kappa, \epsilon; \xi_{1}),
\end{equation}
then the two functions can be glued together to form a global function
satisfying both boundary conditions. With this condition, it is
natural to define the following map
\begin{equation}
  \label{eq:G}
  G(\mu, \gamma, \kappa, \epsilon) = \left(
    Q_{0}(\mu, \kappa, \epsilon; \xi_{1}) - Q_{\infty}(\gamma, \kappa, \epsilon; \xi_{1}),
    Q_{0}'(\mu, \kappa, \epsilon; \xi_{1}) - Q_{\infty}'(\gamma, \kappa, \epsilon; \xi_{1})
  \right) \colon \mathbb{R} \times \mathbb{C} \times \mathbb{R} \times \mathbb{R} \to \mathbb{C}^{2}.
\end{equation}
By identifying \(\mathbb{C}\) with \(\mathbb{R}^{2}\), we can
interpret \(G\) as a map from \(\mathbb{R}^{5}\) to
\(\mathbb{R}^{4}\). A zero of the function \(G\) then corresponds to a
matching as in~\eqref{eq:Q-match}, and ultimately implies the
existence of a self-similar solution to the CGL equation.

In the NLS case, which fixes \(\epsilon = 0\), \(G\) is a function
from \(\mathbb{R}^{4}\) to \(\mathbb{R}^{4}\), and the zero set is
expected to consist of isolated points. The existence of such zeros is
asserted in Theorem~\ref{thm:G-nls} in Section~\ref{sec:nls}, and
discussed more in depth in Section~\ref{sec:proof-nls}. In the general
case, which leaves \(\epsilon\) free, the zero set is expected to be
given by curves in \(\mathbb{R}^{5}\). These curves form the branches
of self-similar singular solutions that are the topic of
Sections~\ref{sec:branches} and~\ref{sec:proof-cgl}.

\section{Self-similar singular solutions for the NLS equation}
\label{sec:nls}
In this section we give a more detailed account of the results on the
existence of self-similar singular solutions to the nonlinear
Schrödinger equation, corresponding to the case
\(\epsilon = \delta = 0\) in Equation~\eqref{eq:CGL}. Following
Plecháč and Šverák~\cite{Plech2001}, we split our results into
two cases: \textbf{Case I} corresponding to \(d = 1\) and
\(\sigma = 2.3\) and \textbf{Case II} corresponding to \(d = 3\) and
\(\sigma = 1\).

As discussed in Section~\ref{sec:shooting-method} above, the proof of
the existence of self-similar singular solutions is based on a
shooting method and ultimately reduces to proving the existence of
zeros of the function \(G\) in Equation~\eqref{eq:G}. The results for
the NLS equation are based on the following theorem regarding zeros of
\(G\), the proof of which is given in Section~\ref{sec:proof-nls}.
\begin{theorem}
  \label{thm:G-nls}
  Let \(\epsilon = \delta = 0\) and consider the function
  \(G\)~\eqref{eq:G}. In \textbf{Case I} (\(d = 1\), \(\sigma = 2.3\)),
  the function \(G\) has (at least) 8 locally unique zeros, in
  \textbf{Case II} (\(d = 3\), \(\sigma = 1\)), it has (at least) 2
  locally unique zeros. Enclosures for these zeros are given in
  Table~\ref{tab:results-nls-d-1} for \textbf{Case I} and
  Table~\ref{tab:results-nls-d-3} for \textbf{Case II}.
\end{theorem}

\begin{table}[ht]
  \centering
  \begin{tabular}{ccccc}
    \toprule
    \(j\) & \(\mu_{j}^{I}\) & \(\gamma_{j}^{I}\) & \(\kappa_{j}^{I}\) & \(\xi_{1}\)\\
    \midrule
    \(1\) & \(1.2320375_{03}^{47}\) & \(0.758150_{50}^{83} - 0.43437_{09}^{13} i\) & \(0.853108_{71}^{96}\) & \(10\)\\
    \(2\) & \(0.78307_{66}^{76}\) & \(1.7916_{36}^{49} - 1.1049_{10}^{24} i\) & \(0.49322_{25}^{32}\) & \(15\)\\
    \(3\) & \(1.12384_{17}^{38}\) & \(4.331_{25}^{32} - 1.5341_{15}^{72} i\) & \(0.34675_{30}^{43}\) & \(20\)\\
    \(4\) & \(0.8838_{74}^{81}\) & \(9.73_{77}^{82} + 0.046_{31}^{61} i\) & \(0.2667_{59}^{61}\) & \(25\)\\
    \(5\) & \(1.07955_{89}^{95}\) & \(18.1568_{25}^{71} + 9.263_{16}^{23} i\) & \(0.216402_{14}^{29}\) & \(40\)\\
    \(6\) & \(0.92718_{34}^{54}\) & \(21.94_{08}^{13} + 36.282_{30}^{84} i\) & \(0.18183_{77}^{81}\) & \(45\)\\
    \(7\) & \(1.054410_{14}^{77}\) & \(-8.221_{20}^{510} + 87.50_{44}^{51} i\) & \(0.156677_{88}^{99}\) & \(60\)\\
    \(8\) & \(0.95114_{57}^{62}\) & \(-130.05_{45}^{55} + 126.992_{05}^{54} i\) & \(0.137564_{37}^{43}\) & \(75\)\\
    \bottomrule
  \end{tabular}
  \caption{\textbf{Case I} (\(d = 1\), \(\sigma = 2.3\),
    \(\epsilon = 0\)): Enclosures of \(\mu\), \(\gamma\) and
    \(\kappa\) corresponding to zeros of \(G\)~\eqref{eq:G} and
    self-similar singular solutions of NLS. The last column gives the
    value of \(\xi_{1}\) used in the computations (compare
    with~\cite[Table 3.2]{Plech2001}). The $j$ index denotes the number of monotone intervals of $|Q|$.}
  \label{tab:results-nls-d-1}
\end{table}

\begin{table}[ht]
  \centering
  \begin{tabular}{cccccccc}
    \toprule
    \(j\) & \(\mu_{j}^{II}\) & \(\gamma_{j}^{II}\) & \(\kappa_{j}^{II}\) & \(\xi_{1}\)\\
    \midrule
    \(1\) & \(1.88565_{69}^{71}\) & \(1.713600_{84}^{98} - 1.491793_{77}^{91} i\) & \(0.9173561_{05}^{81}\) & \(60\)\\
    \(2\) & \(0.83995_{89}^{96}\) & \(13.8526_{16}^{32} + 6.034_{49}^{53} i\) & \(0.3212_{39}^{41}\) & \(140\)\\
    \bottomrule
  \end{tabular}
  \caption{\textbf{Case II} (\(d = 3\), \(\sigma = 1\),
    \(\epsilon = 0\)): Enclosures of \(\mu\), \(\gamma\) and
    \(\kappa\) corresponding to zeros of \(G\)~\eqref{eq:G} and
    self-similar singular solutions of NLS. The last column gives the
    value of \(\xi_{1}\) used in the computations (compare
    with~\cite[Table 3.5]{Plech2001}). The $j$ index denotes the number of monotone intervals of $|Q|$.}
  \label{tab:results-nls-d-3}
\end{table}

\begin{remark}
  The local uniqueness of the zeros of \(G\) is within the chosen
  parametrization of the manifold at infinity for \(Q_{\infty}\). It
  does not give local uniqueness for the self-similar parameter
  \(\kappa\).
\end{remark}

As a consequence, we get the following result about the existence of
self-similar singular solutions to the nonlinear Schrödinger
equation.
\begin{theorem}
  \label{thm:existence-nls}
  The nonlinear Schrödinger equation,
  \begin{align*}
    i \frac{\partial u}{\partial t} + \Delta u + |u|^{2\sigma}u &= 0 \text{ in } \mathbb{R}^{d} \times (0, T),\\
    u(x, 0) &= u_{0}(x),
  \end{align*}
  supports radial self-similar singular solutions of the form
  \begin{equation*}
    u(x, t) = \dfrac{1}{(2\kappa(T - t))^{\frac{1}{2}\left(\frac{1}{\sigma} + i \frac{\omega}{\kappa}\right)}}
    Q\left(\dfrac{|x|}{(2\kappa(T - t))^{\frac{1}{2}}}\right),
  \end{equation*}
  for \(Q \colon [0, \infty) \to \mathbb{C}\) satisfying \(Q(0) = \mu > 0\)
  and \(Q'(0) = 0\). In \textbf{Case I} (\(d = 1\), \(\sigma = 2.3\)),
  such self-similar singular solutions exist for (at least) 8
  different pairs, \(\{(\mu_{j}^{I}, \kappa_{j}^{I})\}_{j = 1}^{8}\),
  and in \textbf{Case II} (\(d = 3\), \(\sigma = 1\)), they exist for
  (at least) 2 different such pairs,
  \(\{(\mu_{j}^{II}, \kappa_{j}^{II})\}_{j = 1}^{2}\). Enclosures of
  \(\{(\mu_{j}^{I}, \kappa_{j}^{I})\}_{j = 1}^{8}\) and
  \(\{(\mu_{j}^{II}, \kappa_{j}^{II})\}_{j = 1}^{2}\) are given in
  Tables~\ref{tab:results-nls-d-1} and~\ref{tab:results-nls-d-3}
  respectively.
\end{theorem}

Moreover, the quantitative control from Theorem~\ref{thm:G-nls} allows
us to establish the following qualitative result for the profiles
\(|Q|\).
\begin{theorem}
  If the self-similar singular solutions of
  Theorem~\ref{thm:existence-nls} in \textbf{Case I} and \textbf{Case
    II} are (separately) ordered by decreasing value of \(\kappa\)
  (such as in Tables~\ref{tab:results-nls-d-1}
  and~\ref{tab:results-nls-d-3}), then for solution \(j\), the profile
  \(|Q|\) has exactly \(j\) intervals of monotonicity on the interval
  \((0, \infty)\). In particular, for \(j = 1\), \(|Q|\) is
  monotonically decreasing.
\end{theorem}

\section{Branches of self-similar singular solutions for the CGL equation}
\label{sec:branches}
In this section we consider the general case, \(\epsilon \geq 0\)
(still fixing \(\delta = 0\)). We give results for the existence of
branches of self-similar singular solutions, originating from the
corresponding solutions to the NLS equation (\(\epsilon = 0\))
discussed in Section~\ref{sec:nls} above. The starting point for our
results is the numerical studies done by Plecháč and
Šverák~\cite{Plech2001}.

\subsection{Numerical approximations}
\label{sec:numerical-approximations-cgl}
\begin{figure}
  \centering
  \begin{subfigure}[t]{0.4\textwidth}
    \includegraphics[width=\textwidth]{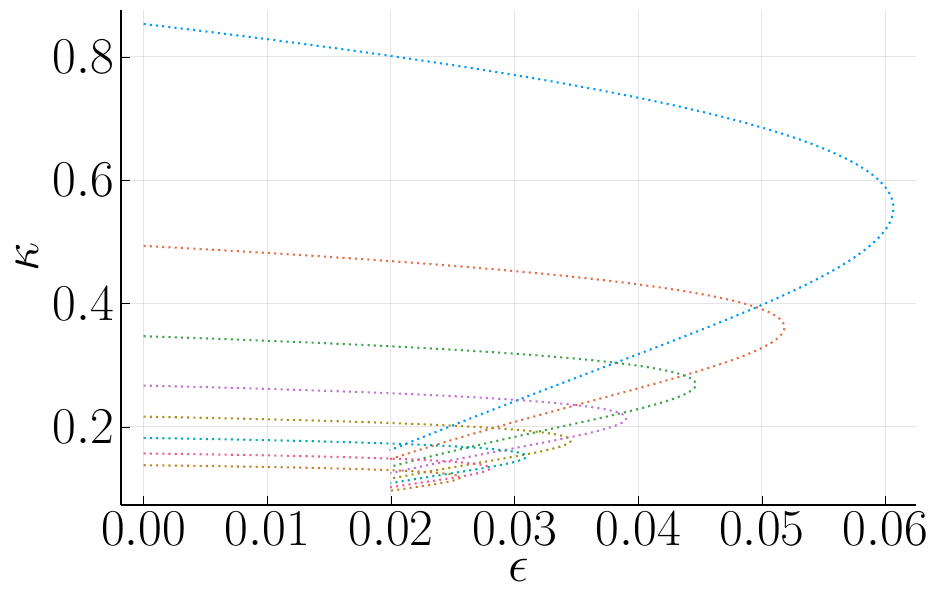}
    \caption{\textbf{Case I} (\(d = 1\), \(\sigma = 2.3\))}
    \label{fig:branches-d1-numerical}
  \end{subfigure}
  \hspace{0.05\textwidth}
  \begin{subfigure}[t]{0.4\textwidth}
    \includegraphics[width=\textwidth]{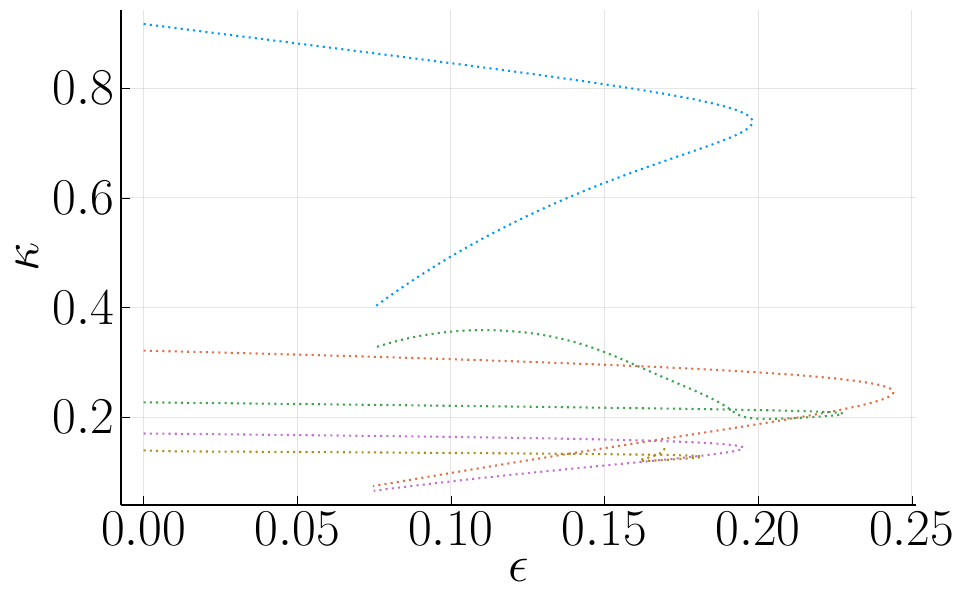}
    \caption{\textbf{Case II} (\(d = 3\), \(\sigma = 1\))}
    \label{fig:branches-d3-numerical}
  \end{subfigure}
  \caption{Numerical approximations of branches, compare
    with~\cite[Figures 3.1 and 3.6]{Plech2001}.}
  \label{fig:branches-numerical}
\end{figure}

Plecháč and Šverák studied the branches using a combination of
numerical and analytical tools. Implementing a similar procedure
allows us to compute numerical approximations of the branches. These
numerical approximations are given in
Figure~\ref{fig:branches-numerical} and have been computed with the
help of the Julia package
\texttt{BifurcationKit.jl}~\cite{BifurcationKit} (see
Section~\ref{sec:implementation-details} for more details about the
implementation).

For \textbf{Case I}, our figure agrees with the one obtained by Plecháč
and Šverák. However, for \textbf{Case II}, our figure indicates a
qualitative behavior that is missing from that of Plecháč and Šverák.
The third and fifth branches appear to not be monotonically decreasing
in \(\kappa\), instead they make a turn upwards shortly after the turn
in \(\epsilon\). The behavior of the third branch indicates that they
eventually turn back and, like the other branches, ultimately tend
towards the origin. For the fifth branch we are not able to run the
continuation for as far as the other branches and, while it seems
likely that it would ultimately tend towards the origin, the numerical
results do not make this clear. It should be noted that our
observations do not contradict those made by Plecháč and Šverák, their
continuation of the third and fifth branches stops right before they
turn in \(\kappa\).

\subsection{Rigorous results}
\label{sec:rigorous-results-cgl}
\begin{figure}
  \centering
  \begin{subfigure}[t]{0.4\textwidth}
    \includegraphics[width=\textwidth]{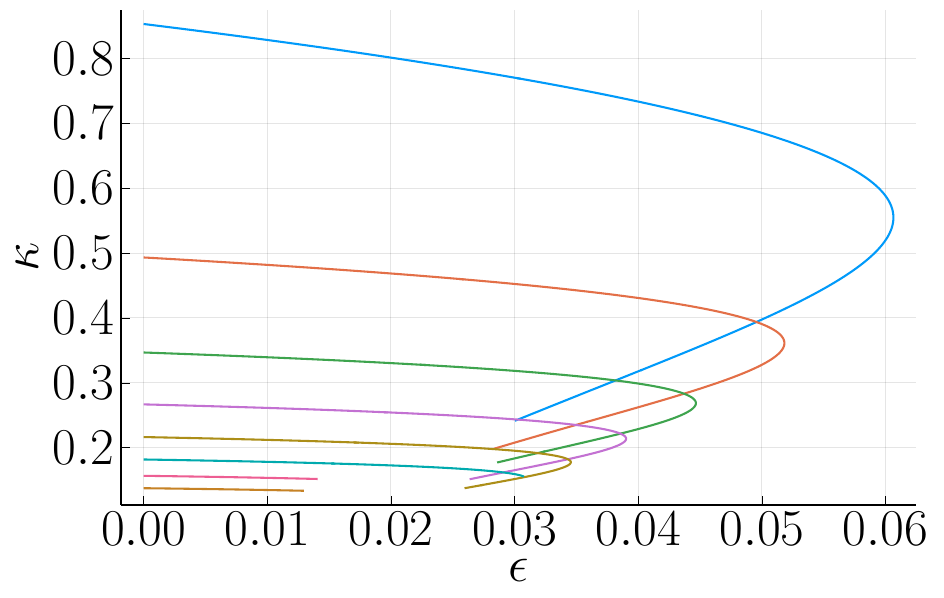}
    \caption{\textbf{Case I} (\(d = 1\), \(\sigma = 2.3\))}
    \label{fig:branches-d1}
  \end{subfigure}
  \hspace{0.05\textwidth}
  \begin{subfigure}[t]{0.4\textwidth}
    \includegraphics[width=\textwidth]{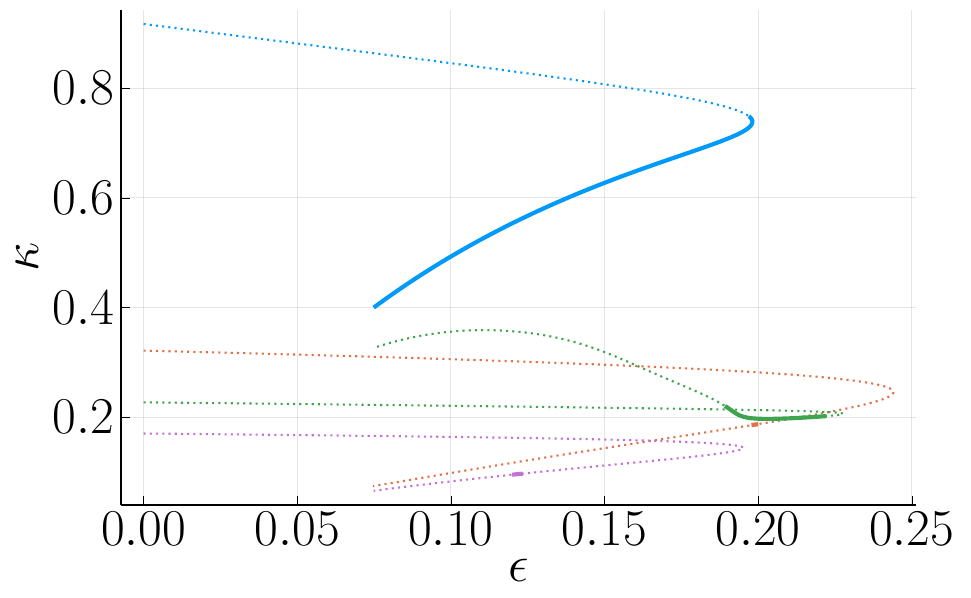}
    \caption{\textbf{Case II} (\(d = 3\), \(\sigma = 1\))}
    \label{fig:branches-d3}
  \end{subfigure}
  \caption{Rigorously verified branches of self-similar singular
    solutions for the CGL equation. For \textbf{Case II}, the filled-in
    parts of the curves are the verified parts, whereas the dotted
    parts are the numerical approximations from
    Figure~\ref{fig:branches-d3-numerical}. For \textbf{Case I}, all
    shown parts are verified.}
  \label{fig:branches}
\end{figure}

Let us now focus on our rigorous results for these branches. Our
results regarding the existence of the branches for \textbf{Case I} are
summarized in the following theorem, for which the proof is discussed
in Section~\ref{sec:proof-cgl}.
\begin{theorem}
  \label{thm:existence-cgl-case-1}
  In \textbf{Case I} (\(d = 1\), \(\sigma = 2.3\)), there exist (at
  least) 8 continuous curves,
  \(\{(\epsilon_{j}^{I}(s), \mu_{j}^{I}(s), \kappa_{j}^{I}(s))\}_{j =
    1}^{8}\), such that for each \(s \in [0, 1]\) and each \(j\), the
  CGL equation with \(\epsilon = \epsilon_{j}^{I}(s)\) has a
  self-similar singular solution of the form
  \begin{equation*}
    u(x, t) = \dfrac{1}{(2\kappa_{j}^{I}(s)(T - t))^{\frac{1}{2}\left(\frac{1}{\sigma} + i \frac{\omega}{\kappa_{j}^{I}(s)}\right)}}
    Q\left(\dfrac{|x|}{(2\kappa_{j}^{I}(s)(T - t))^{\frac{1}{2}}}\right),
  \end{equation*}
  for \(Q \colon [0, \infty) \to \mathbb{C}\) satisfying
  \(Q(0) = \mu_{j}^{I}(s)\) and \(Q'(0) = 0\). The curves satisfy
  \((\epsilon_{j}^{I}(0), \mu_{j}^{I}(0), \kappa_{j}^{I}(0)) = (0,
  \mu_{j}^{I}, \kappa_{j}^{I})\), with \(\kappa_{j}^{I}\) and
  \(\mu_{j}^{I}\) as in Theorem~\ref{thm:existence-nls}, as well as
  \(\epsilon_{j}^{I\prime}(0) > 0\). Moreover, along the whole curve
  \(j\), the profile \(|Q|\) has \(j\) intervals of monotonicity on
  the interval \((0, \infty)\). The \((\epsilon, \kappa)\)-projections
  of these curves are depicted in Figure~\ref{fig:branches-d1}.
\end{theorem}

\begin{remark}
  The theorem as stated here is very light on quantitative properties
  of the parameter curves
  \((\epsilon_{j}^{I}(s), \mu_{j}^{I}(s), \kappa_{j}^{I}(s))\). The
  procedure for the proof does however give us very strong
  quantitative control of the curves at a \(C^{0}\) level: it produces
  a thin tube in which the curve is guaranteed to lie.
  Figure~\ref{fig:branches-d1} is based on this tube. Precise
  numerical data for these curves are provided in the repository
  accompanying the computer-assisted proof~\cite{CGL.jl}.
\end{remark}

\begin{remark}
  The reason that some of the curves in Figure~\ref{fig:branches-d1}
  stop earlier than the numerical approximations of the curves given
  in Figure~\ref{fig:branches-d1-numerical} is that rigorous
  verification is much more computationally demanding, and to avoid
  excessive computational time we stop the computations earlier. More
  details about the computational cost are given in
  Section~\ref{sec:other-cases-cgl}.
\end{remark}

For \textbf{Case II} the problem is significantly more computationally
demanding, in particular near the beginning of the branches. For that
reason our rigorous results only cover \(1 \leq j \leq 4\), and only
parts of those branches. Some of the difficulties encountered are
discussed in Section~\ref{sec:pointwise-results-cgl}, where we give
partial, pointwise, results about existence of self-similar singular
solutions along the branches. For this same reason we also skip the
verification of the number of intervals of monotonicity for the
profile \(|Q|\) in this case. The proof of the theorem is discussed in
Section~\ref{sec:proof-cgl}.

\begin{theorem}
  \label{thm:existence-cgl-case-2}
  In \textbf{Case II} (\(d = 3\), \(\sigma = 1\)), there exist (at
  least) 4 continuous curves,
  \(\{(\epsilon_{j}^{II}(s), \mu_{j}^{II}(s), \kappa_{j}^{II}(s))\}_{j
    = 1}^{4}\), such that for each \(s \in [0, 1]\) and each \(j\), the
  CGL equation with \(\epsilon = \epsilon_{j}^{II}(s)\) has a
  self-similar singular solution of the form
  \begin{equation*}
    u(x, t) = \dfrac{1}{(2\kappa_{j}^{II}(s)(T - t))^{\frac{1}{2}\left(\frac{1}{\sigma} + i \frac{\omega}{\kappa_{j}^{II}(s)}\right)}}
    Q\left(\dfrac{|x|}{(2\kappa_{j}^{II}(s)(T - t))^{\frac{1}{2}}}\right),
  \end{equation*}
  for \(Q \colon [0, \infty) \to \mathbb{C}\) satisfying
  \(Q(0) = \mu_{j}^{II}(s)\) and \(Q'(0) = 0\). The
  \((\epsilon, \kappa)\)-projections of these curves are depicted,
  together with the numerical approximations for larger parts of the
  curves, in Figure~\ref{fig:branches-d3}.
\end{theorem}

\section{Proofs for the NLS equation}
\label{sec:proof-nls}
As stated in Section~\ref{sec:nls}, we have several results for the
existence of self-similar solutions for the NLS equation. All of them
follow the same strategy but with different initial data and
constants, so we detail an example below to illustrate our procedure.

For the proof of existence, the example goes through all parts of the
proof, except for one crucial step: How to compute rigorous interval
enclosures of the function \(G\)~\eqref{eq:G} and its Jacobian. How
these interval enclosures can be computed is the topic of
Section~\ref{sec:solution-infinity}, which handles the computation of
\(Q_{\infty}\), and Section~\ref{sec:solution-zero}, which handles the
computation of \(Q_{0}\). Together, these two sections constitute the
main part of the paper.

\subsection{A detailed example}
\label{sec:detailed-example-nls}
We go through the proof of existence of the self-similar singular
solution, as well as the verification of the monotonicity, for the
first (\(j = 1\)) solution in \textbf{Case II} (\(d = 3\),
\(\sigma = 1\)).

Recall that the proof of existence is based on proving the existence
of a zero of the map \(G(\mu, \gamma, \kappa)\), defined in
Equation~\eqref{eq:G} (here we fix \(\epsilon = 0\) since we are
considering the NLS case). By splitting \(\gamma\) into real and
imaginary parts, we can treat \(G\) as a map from \(\mathbb{R}^{4}\) to
\(\mathbb{R}^{4}\). To prove the existence of a root, we make use of the
\emph{Krawczyk interval Newton method}, see e.g.~\cite{Moore1977}. To
apply the interval Newton method we need to find a box
\(X \subset \mathbb{R}^{4}\), such that
\begin{equation}
  \label{eq:newton-condition}
  \midint(X) - YG(\midint(X)) + (I - YJ_{G}(X))(X - \midint(X)) \subseteq \operatorname{int}(X),
\end{equation}
where \(\midint(X)\) denotes the midpoint of \(X\),
\(\operatorname{int}(X)\) denotes the interior of \(X\), \(J_{G}(X)\)
denotes a box enclosing the Jacobian of \(G\) and \(Y\) is a
preconditioner that in practice is taken to be an approximate inverse
of the midpoint of \(J_{G}(X)\). If we find \(X\) satisfying this
condition, then the function \(G\) has a unique zero in \(X\), and
this zero is contained in the set given by the left-hand side of the
above expression.

The first step in finding the set \(X\) is to find a good numerical
(noninterval) approximation of the zero. This is done by starting with
a rough approximation, which is refined using a few standard Newton
iterations. The resulting approximation, rounded to 16 digits, is
\begin{equation}\label{eq:approximation}
  (\mu_{0}, \gamma_{0}, \kappa_{0}) =
  (1.885656963035410, 1.713600956210395 -1.491793873532361i, 0.9173561176773725).
\end{equation}
Using interval arithmetic to enclose
\(G(\mu_{0}, \gamma_{0}, \kappa_{0})\), we obtain the interval box
\begin{multline*}
  G(\mu_0, \gamma_0, \kappa_0) \in
  \big([-6.0 \cdot 10^{-12}, 9.2 \cdot 10^{-12}], [-1.8 \cdot 10^{-12}, 1.4 \cdot 10^{-11}],\\
  [-4.0 \cdot 10^{-11}, 7.1 \cdot 10^{-10}], [-4.7 \cdot 10^{-10}, 2.9 \cdot 10^{-10}]\big).
\end{multline*}

\begin{remark}
That these enclosures are very small, and straddle zero, is an
indication that we have a good approximation, and that what is
limiting is the error bounds in the evaluation of \(G\).
\end{remark}

Similarly, we can compute an enclosure of the Jacobian at the
approximation~\eqref{eq:approximation}, obtaining
\begin{equation*}
  J_{G}(\mu_0, \gamma_0, \kappa_0) \in
  \begin{pmatrix}
    -0.002175519_{7906}^{8269} & 0.010390_{7685}^{96410} & 0.001160_{3508}^{5473} & 0.01974_{4274}^{8832} \\
    -0.004673487_{6795}^{7156} & -0.001160_{3469}^{5434} & 0.010390_{7686}^{9651} & -0.0112_{6820}^{7276} \\
    0.003856604_{1674}^{5468} & -0.0001_{8879}^{99510} & 0.0001_{6388}^{7469} & -0.017_{2705}^{5209} \\
    0.00943802_{2710}^{30810} & -0.0001_{6388}^{7469} & -0.0001_{8901}^{9982} & -0.005_{3686}^{61810} \\
  \end{pmatrix}.
\end{equation*}

With the approximation~\eqref{eq:approximation} in our hands, the next
step is to find a set \(X \subset \mathbb{R}^{4}\) around this
approximation, such that the condition of the interval Newton
method~\eqref{eq:newton-condition} can be verified. The set \(X\) is
found using a heuristic method that is described in more detail in
Section~\ref{sec:implementation-details}. In this section we take
\(X\) as given, and only focus on the a posteriori verification of the
interval Newton condition. More precisely, we
take\footnote{Technically, the computations use a slightly different
  value for \(X\), which, rounded outwards, becomes this value.}
\begin{equation*}
  X = 1.88565_{6923}^{7051} \times
  1.713600_{8426}^{97510} \times
  -1.491793_{7777}^{9044} \times
  0.9173561_{0507}^{8097}.
\end{equation*}
For this \(X\) we can compute the enclosures
\begin{equation*}
  G(\operatorname{mid}(X)) \in \left([\pm 7.6 \cdot 10^{-12}], [\pm 7.5 \cdot 10^{-12}], [\pm 3.8 \cdot 10^{-10}], [\pm 3.8 \cdot 10^{-10}]\right),
\end{equation*}
and
\begin{equation*}
  J_{G}(X) \subseteq
  \begin{pmatrix}
    -0.0021_{6875}^{8228} & 0.010390_{7541}^{9798} & 0.001160_{3462}^{5499} & 0.0197_{0036}^{9269} \\
    -0.0046_{6696}^{79910} & -0.001160_{3423}^{54510} & 0.010390_{7542}^{9799} & -0.011_{2258}^{3152} \\
    0.0038_{2835}^{8571} & -0.0001_{8879}^{99510} & 0.0001_{6388}^{74610} & -0.017_{0841}^{7075} \\
    0.0094_{0974}^{6594} & -0.0001_{6387}^{7469} & -0.0001_{8900}^{9982} & -0.005_{1846}^{8002} \\
  \end{pmatrix}.
\end{equation*}
We take \(Y\) to be an approximation of the inverse of the midpoint of
\(J_{G}(X)\), giving us, rounded to 8 digits,
\begin{equation*}
  Y =
  \begin{pmatrix}
    1.5944878 & 2.7710729 & -39.0238 & 123.64427 \\
    96.53897 & -13.610825 & 130.32019 & -37.745564 \\
    10.827354 & 97.82849 & -77.125946 & 82.458626 \\
    -0.6188105 & 1.7183731 & -68.34296 & 28.638468 \\
  \end{pmatrix}.
\end{equation*}
With this, we obtain
\begin{multline*}
  \operatorname{mid}(X) - YG(\operatorname{mid}(X)) + (I - YJ_{G}(X))(X - \operatorname{mid}(X))\\\subseteq
  \left(1.88565_{6924}^{7051}, 1.713600_{8432}^{9754}, -1.491793_{7783}^{9038}, 0.9173561_{0552}^{8052}\right).
\end{multline*}

It can now be verified that this is contained in the interior of
\(X\), so the interval Newton condition~\eqref{eq:newton-condition} is
fulfilled. The interval Newton method then implies that the function
\(G\) has a unique root in the set \(X\) and this root satisfies that
\begin{equation*}
  \mu \in 1.88565_{6924}^{7051},\quad
  \gamma \in 1.713600_{8432}^{9754} - 1.491793_{7783}^{9038} i,\quad
  \kappa \in 0.9173561_{0552}^{8052}.
\end{equation*}

\begin{figure}
  \centering
  \begin{subfigure}[t]{0.45\textwidth}
    \includegraphics[width=\textwidth]{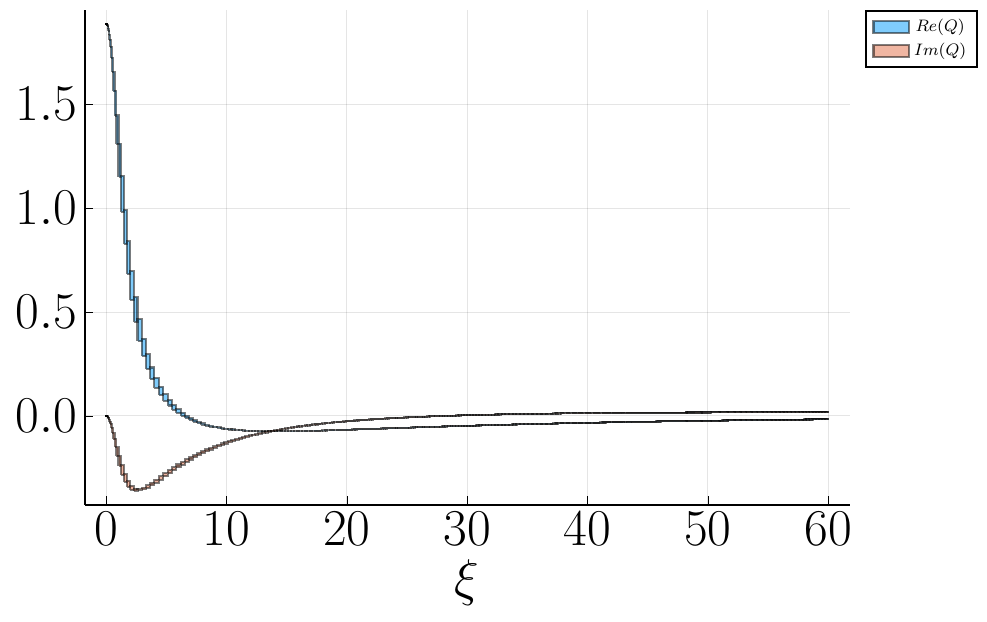}
  \end{subfigure}
  \hspace{0.05\textwidth}
  \begin{subfigure}[t]{0.45\textwidth}
    \includegraphics[width=\textwidth]{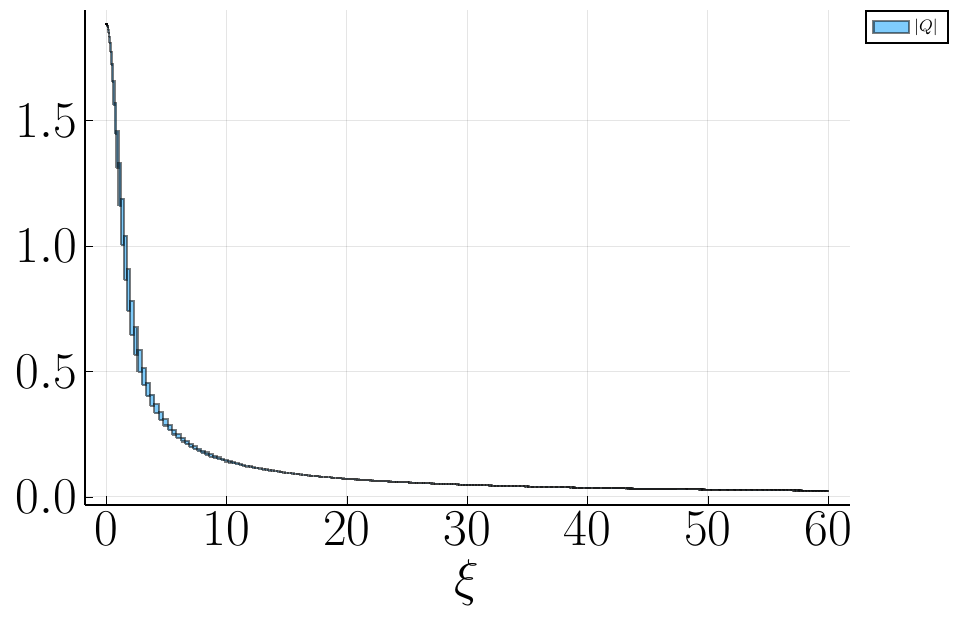}
  \end{subfigure}
  \caption{Rigorous enclosures of the profile \(Q\) corresponding to
    the first (\(j = 1\)) self-similar singular solution of the NLS
    equation in \textbf{Case II} (\(d = 3\) and \(\sigma = 1\)). The
    plots consist of thin boxes, and the profile is proved to be
    contained inside these boxes.}
  \label{fig:nls-example}
\end{figure}

With the rigorous enclosures of \(\mu\), \(\gamma\) and \(\kappa\) we
can, using the same rigorous numerical integration as when enclosing
\(G\), see Section~\ref{sec:solution-zero}, compute enclosures of the
profile \(Q\) on the whole interval \([0, \xi_{1}]\). The enclosures
are depicted in Figure~\ref{fig:nls-example}. The plot of \(|Q|\)
indicates that it is decreasing, at least on the interval
\([0, \xi_{1}]\); our next step is to prove that this is indeed the
case.

To prove that \(|Q|\) is decreasing, it is enough to verify that
\(\frac{d}{d\xi}|Q|\) is nonzero on the interval \((0, \infty)\). In
practice it is easier to work with \(\frac{d}{d\xi}|Q|^{2}\), and what
we will prove is that \(\frac{d}{d\xi}|Q|^{2}\) is nonzero on
\((0, \infty)\). We split the interval \((0, \infty)\) into three
parts, \((0, \xi_{0})\), \([\xi_{0}, \xi_{1}]\) and
\((\xi_{1}, \infty)\).

The verification on the interval \((\xi_{1}, \infty)\) is based on the
bounds in Section~\ref{sec:verify-monotonicity}, with the final
condition given in Lemma~\ref{lemma:monotonicity-infinity}.
Lemma~\ref{lemma:abs2-Q-derivative} tells us that for
\(\xi \geq \xi_{1}\), we have
\begin{equation*}
  \frac{d}{d\xi}|Q|^{2}
  = 2|p_{Q}(\xi)|^{2}\real\left(P(\xi)\conj{P'(\xi)}\right)
  + 2R_{\mathrm{mon}}(\xi)\xi^{(2\sigma + 1)\normv - \frac{2}{\sigma} - 3}.
\end{equation*}
Here the function \(P\) is the confluent hypergeometric function,
which determines the leading order behavior of \(Q\) at infinity, see
Section~\ref{sec:solution-infinity} for more details. To prove
monotonicity, it suffices to verify that
\begin{equation*}
  \left|\real\left(P(\xi)\conj{P'(\xi)}\right)\right| > \frac{|R_{\mathrm{mon}}(\xi)|}{|p_{Q}(\xi)|^{2}}\xi_{1}^{(2\sigma + 1)\normv - \frac{2}{\sigma} - 3}.
\end{equation*}
Taking into account that \(P(\xi)\conj{P'(\xi)}\) decays like
\(\xi^{-\frac{2}{\sigma} - 1}\), we write this as
\begin{equation*}
  \left|\real\left(P(\xi)\conj{P'(\xi)}\xi^{\frac{2}{\sigma} + 1}\right)\right| > \frac{R_{\mathrm{mon}}(\xi)}{|p_{Q}(\xi)|^{2}}\xi_{1}^{(2\sigma + 1)\normv - 2}.
\end{equation*}
In this case the left-hand side has a non-zero limit at infinity and
using Lemma~\ref{lemma:real-P-P-dxi-expansion}, we can compute the
enclosure
\begin{equation*}
  \left|\operatorname{real}\left(P(\xi)\overline{P'(\xi)}\xi^{\frac{2}{\sigma} + 1}\right)\right|
  \in 0.393_{28}^{68}.
\end{equation*}
For the right-hand side we can use Lemmas~\ref{lemma:Q-expansion}
and~\ref{lemma:abs2-Q-derivative} to compute the upper bound
\begin{equation*}
  \frac{|R_{\mathrm{mon}}(\xi)|}{|p_{Q}(\xi)|^{2}}\xi_{1}^{(2\sigma + 1)\mathrm{v} - 2} \leq 0.00025_{08}^{14}.
\end{equation*}
Here the inequality with the interval means that there exists a number
in the interval such that this inequality holds. One can then verify
that the required inequality is satisfied.

\begin{figure}
  \centering
  \begin{subfigure}[t]{0.45\textwidth}
    \includegraphics[width=\textwidth]{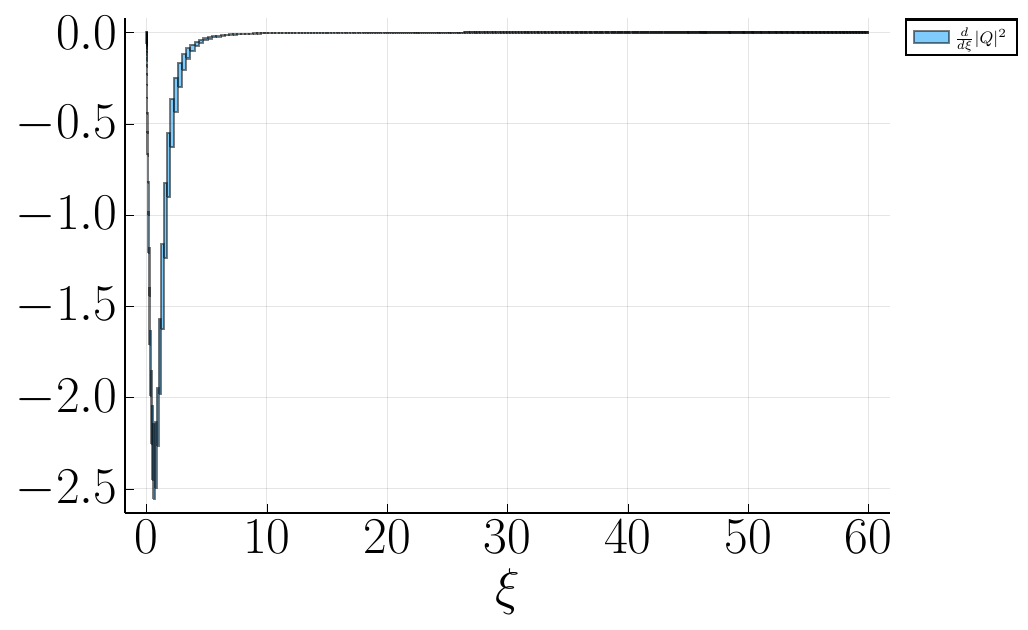}
  \end{subfigure}
  \hspace{0.05\textwidth}
  \begin{subfigure}[t]{0.45\textwidth}
    \includegraphics[width=\textwidth]{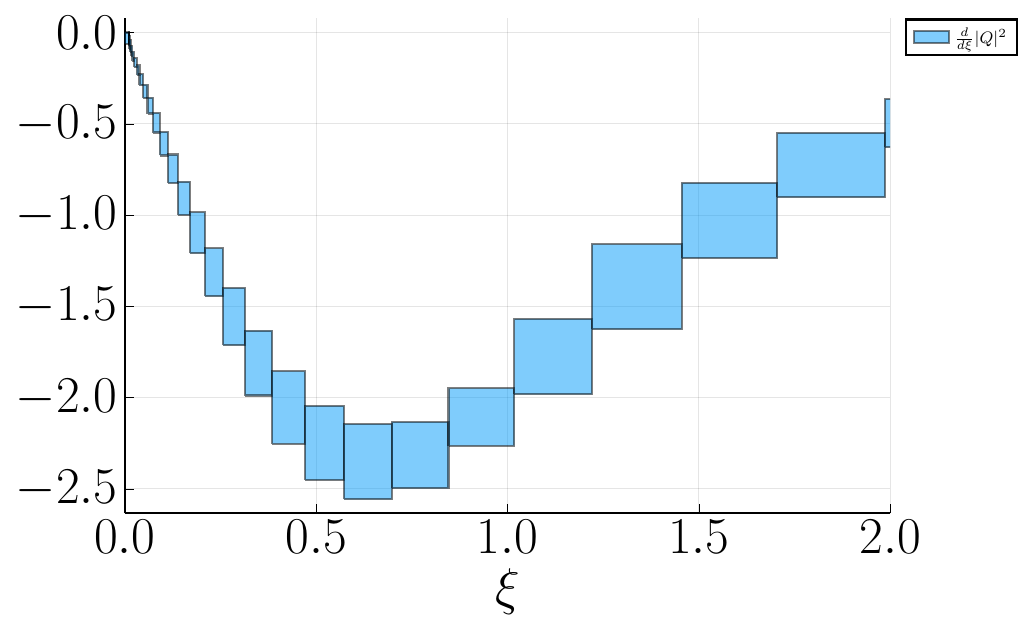}
  \end{subfigure}
  \caption{Rigorous enclosure of \(\frac{d}{d\xi}|Q|^{2}\), on the
    interval \([0, \xi_{1}]\) to the left, and near zero to the
    right.}
  \label{fig:nls-example-derivative}
\end{figure}

What remains is checking the monotonicity on the intervals
\((0, \xi_{0})\) and \([\xi_{0}, \xi_{1}]\). As a first step, we use
the rigorous numerical integrator, see
Section~\ref{sec:solution-zero}, to enclose \(\frac{d}{d\xi}|Q|^{2}\)
on the interval \([0, \xi_{1}]\), the resulting enclosure is depicted
in Figure~\ref{fig:nls-example-derivative}. Except for the two
enclosing boxes closest to zero, all the enclosing boxes can be
verified to not contain zero. Taking \(\xi_{0}\) to be given by the
upper bound of \(\xi\) for the second box, giving
(rounded\footnote{The exact value is
  \(\xi_0 =
  0.010091552734488686178526695158552684006281197071075439453125\)} to
10 digits) \(\xi_{0} = 0.01009155273\), is then enough for us to be
able to verify the monotonicity on \([\xi_{0}, \xi_{1}]\).

The final step is to verify the monotonicity on \((0, \xi_{0})\). Due
to the boundary conditions at \(\xi = 0\), we have that
\(\frac{d}{d\xi}|Q|^{2}\) is zero at \(\xi = 0\). To verify that
\(\frac{d}{d\xi}|Q|^{2}\) is nonzero on \((0, \xi_{0})\), it is
therefore enough to verify that \(\frac{d^{2}}{d\xi^{2}}|Q|^{2}\) is
nonzero on \((0, \xi_{0})\). For this we use the Taylor expansion at
\(\xi = 0\) to compute an enclosure on \([0, \xi_{0}]\) giving us
that, over this interval, \(\frac{d^{2}}{d\xi^{2}}|Q|^{2}\) is
contained in \(-6.05_{57}^{95}\). Since this enclosure is nonzero, it
guarantees the monotonicity of \(|Q|\) on \((0, \xi_{0})\).

Combining the above monotonicity results on \((0, \xi_{0})\),
\([\xi_{0}, \xi_{1}]\) and \((\xi_{1}, \infty)\), gives us that
\(|Q|\) is monotone on the entire interval \((0, \infty)\).

\begin{remark}
  For the non-monotone solutions, \(j \geq 2\), slightly more work is
  required to correctly count the number of zeros of \(|Q|'\) on the
  interval \([\xi_{0}, \xi_{1}]\). We isolate the zeros by observing
  where the enclosures of \(\frac{d}{d\xi}|Q|^{2}\) change sign, the
  local uniqueness of the zeros can then be verified by checking that
  \(\frac{d^{2}}{d\xi^{2}}|Q|^{2}\) has a constant (nonzero) sign
  between the sign changes of \(\frac{d}{d\xi}|Q|^{2}\). That
  \(\frac{d^{2}}{d\xi^{2}}|Q|^{2}\) is nonzero also ensures that the
  zero of \(\frac{d}{d\xi}|Q|^{2}\) indeed corresponds to a zero of
  \(|Q|'\) (and not \(|Q|\)).
\end{remark}

\section{Proofs for the CGL equation}
\label{sec:proof-cgl}
Similarly to the NLS equation, we go through a detailed example that
illustrates the procedure for proving the existence of the branches of
self-similar singular solutions that were presented in
Section~\ref{sec:branches}. Afterwards, we give some more detailed
information about the computed branches and discuss some of the
difficulties in \textbf{Case II}.

\subsection{A detailed example}
\label{sec:detailed-example-cgl}
We go through a detailed example for the first (\(j = 1\)) branch in
\textbf{Case I} (\(d = 1\), \(\sigma = 2.3\)).

\begin{figure}
  \centering
  \begin{subfigure}[t]{0.4\textwidth}
    \includegraphics[width=\textwidth]{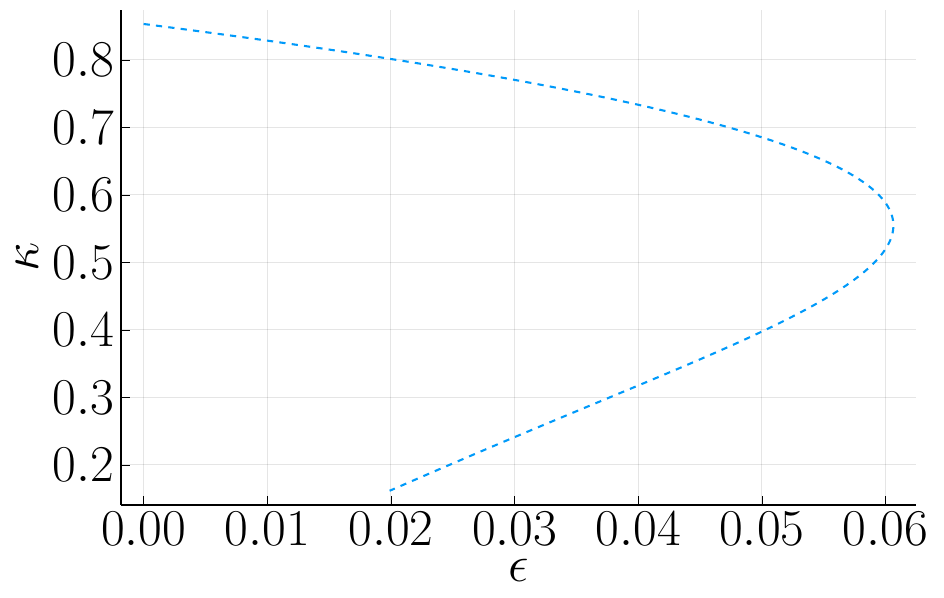}
  \end{subfigure}
  \hspace{0.05\textwidth}
  \begin{subfigure}[t]{0.4\textwidth}
    \includegraphics[width=\textwidth]{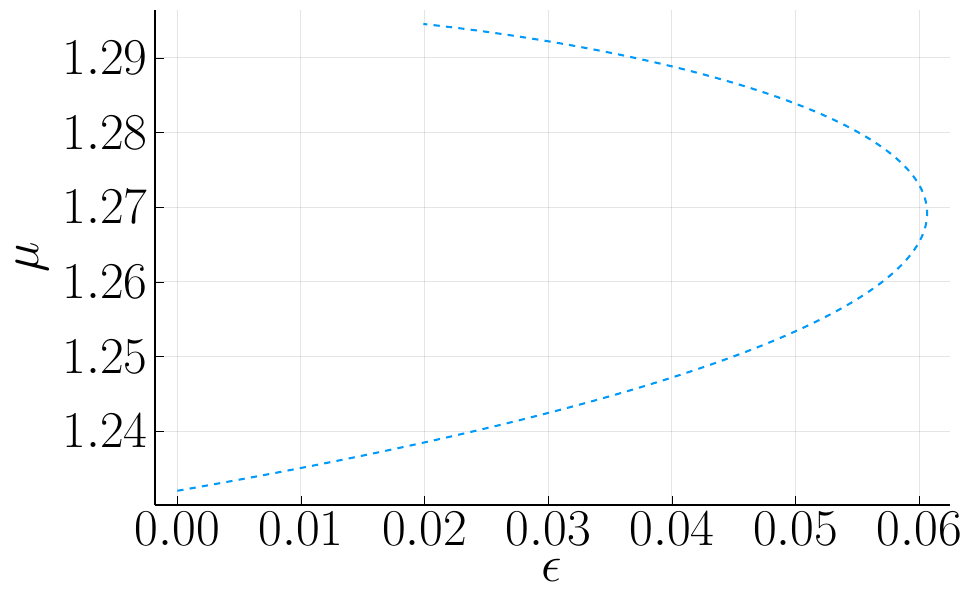}
  \end{subfigure}
  \caption{Approximation of the first branch of self-similar singular
    solutions for the CGL equation in \textbf{Case I}. To the left the
    \((\epsilon, \kappa)\)-projection of the curve, and to the right
    the \((\epsilon, \mu)\)-projection.}
  \label{fig:CGL-example-approximation}
\end{figure}

The first step is to compute a numerical approximation of the
parameter curve. The heavy lifting is done by the Julia package
\texttt{BifurcationKit.jl}~\cite{BifurcationKit}, which handles
automatic continuation of solutions. The result is a curve in
\((\epsilon, \mu, \kappa)\), depicted in
Figure~\ref{fig:CGL-example-approximation}. Note that the parameter
\(\gamma\) (used to parameterize the manifold at infinity) does not
directly appear in the numerical approximation, it is only used as an
internal variable at this stage.

\begin{figure}
  \centering
  \begin{subfigure}[t]{0.45\textwidth}
    \includegraphics[width=\textwidth]{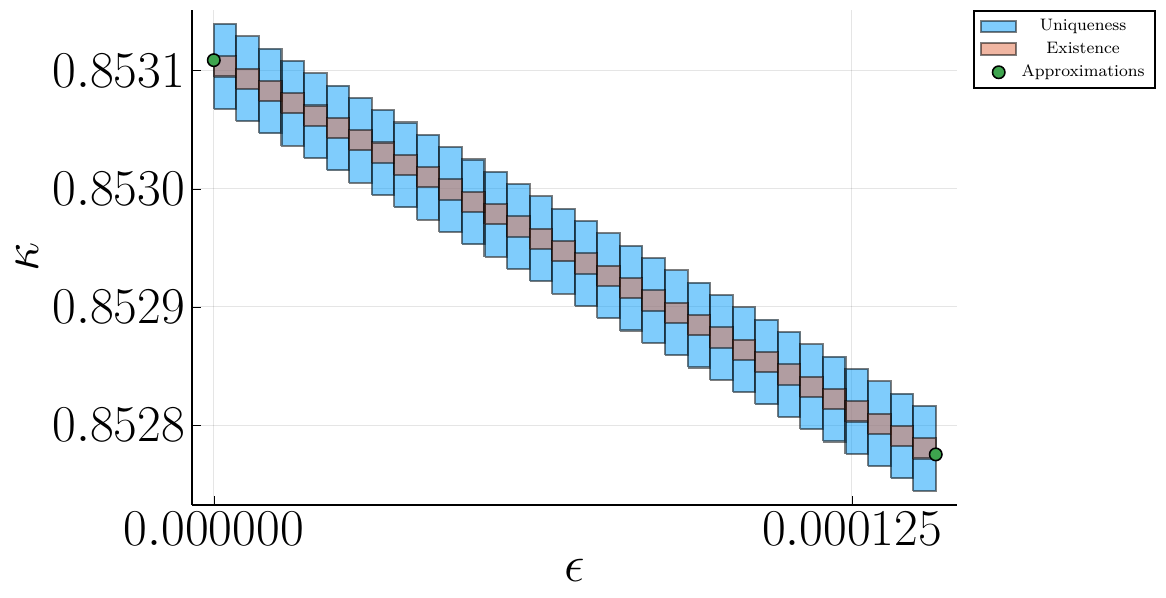}
    \caption{Boxes enclosing the beginning of the first branch in
      \textbf{Case I}.}
    \label{fig:CGL-example-beginning}
  \end{subfigure}
  \hspace{0.05\textwidth}
  \begin{subfigure}[t]{0.45\textwidth}
    \includegraphics[width=\textwidth]{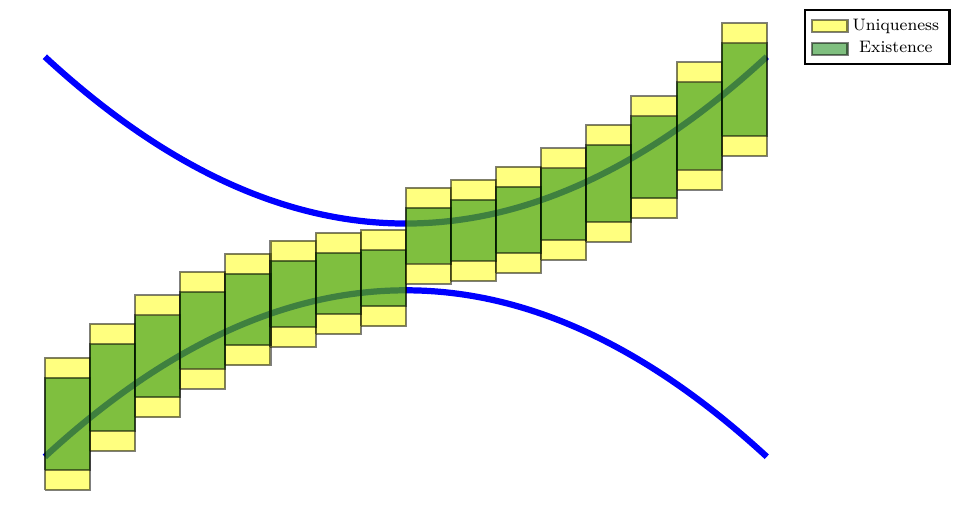}
    \caption{Sketch of example where separate curve segments fail to
      form one continuous curve.}
    \label{fig:CGL-example-uniqueness-issue}
  \end{subfigure}
  \caption{}
\end{figure}

Let us start by considering the problem of verifying the existence at
the beginning of the curve, before we discuss the verification of the
full curve. The first two points on the branch in our numerical
approximation are given by
\begin{equation*}
  \epsilon_1 = 0.0,\quad
  \mu_1 = 1.232037523210031,\quad
  \kappa_1 = 0.8531088807225596
\end{equation*}
and
\begin{equation*}
  \epsilon_2 = 0.0001414213562373095,\quad
  \mu_2 = 1.2320784725879474,\quad
  \kappa_2 = 0.8527751258772226.
\end{equation*}
We will prove the existence of a curve of solutions, starting at
\(\epsilon_1\), and ending at \(\epsilon_2\). This proof will come
together as a collection of boxes covering the interval, where each
box is proved to contain a solution for each \(\epsilon\) in the box,
see Figure~\ref{fig:CGL-example-beginning}. More precisely, the goal
is to prove that there exists a continuous curve,
\((\epsilon(s), \mu(s), \kappa(s))\), such that for each point on this
curve, the CGL equation has a self-similar singular solution. The
\((\epsilon, \kappa)\)-projection of this curve will be contained
inside the existence-boxes of Figure~\ref{fig:CGL-example-beginning}.
For the beginning of the branch, which is what we are considering now,
we can take the curve to be parameterized by \(\epsilon\). We are then
looking for a curve of the form
\((\epsilon, \mu(\epsilon), \kappa(\epsilon))\), with
\(\epsilon \in [\epsilon_1, \epsilon_2]\). In practice, the curve we
compute will also contain the value of \(\gamma\), so the curve is
given by
\((\epsilon, \mu(\epsilon), \gamma(\epsilon), \kappa(\epsilon))\), but
in the statement of Theorem~\ref{thm:existence-cgl-case-1} the
\(\gamma\) is dropped.

To prove the existence, the interval \([\epsilon_{1}, \epsilon_{2}]\)
is split into \(N\) consecutive subintervals
\begin{equation*}
  [\epsilon_1, \epsilon_2] = \bigcup_{i = 1}^{N} \inter{\epsilon}_{i}.
\end{equation*}
For each subinterval \(\inter{\epsilon}_{i}\) the proof of existence
and uniqueness follows exactly the same approach as in the example for
the NLS equation in Section~\ref{sec:detailed-example-nls}. The only
difference being that for the NLS equation we fixed \(\epsilon = 0\),
whereas in this case we do the computations with \(\epsilon\)
represented by the interval \(\inter{\epsilon}_i\). The properties of
interval arithmetic then ensures that any results we get out of
this procedure, in this case enclosures of existence and uniqueness,
are valid for any \(\epsilon \in \inter{\epsilon}_i\). For the
computations to succeed the subintervals \(\inter{\epsilon}_{i}\)
cannot be too wide, otherwise we cannot compute tight enough
enclosures to verify the interval Newton
condition~\eqref{eq:newton-condition}. In the case we are considering
here, it is enough to split \([\epsilon_{1}, \epsilon_{2}]\) uniformly
into \(N = 32\) subintervals.

The procedure produces, for each \(\inter{\epsilon}_{i}\), two boxes,
\begin{equation*}
  \inter{\mu}_{i}^{(e)} \times \inter{\gamma}_{i}^{(e)} \times \inter{\kappa}_{i}^{(e)}
  \subsetneq
  \inter{\mu}_{i}^{(u)} \times \inter{\gamma}_{i}^{(u)} \times \inter{\kappa}_{i}^{(u)},
\end{equation*}
where for all \(\epsilon \in \inter{\epsilon}_{i}\) the function \(G\)
is proved to have a zero in the inner box, and this zero is unique in
the outer box. The boxes
\(\inter{\epsilon}_{i} \times \inter{\kappa}_{i}^{(e)}\) and
\(\inter{\epsilon}_{i} \times \inter{\kappa}_{i}^{(u)}\) are the ones
depicted in Figure~\ref{fig:CGL-example-beginning}. Moreover, the
Jacobian of \(G\) is guaranteed to be invertible for all
\(\epsilon \in \inter{\epsilon}_{i}\), and parameters in the outer
box. Since \(G\) is continuously differentiable with respect to all of
its arguments, it follows by the implicit function theorem that there
exists a continuous curve of solutions
\((\epsilon, \mu_{i}(\epsilon), \gamma_{i}(\epsilon),
\kappa_{i}(\epsilon))\) defined for
\(\epsilon \in \inter{\epsilon}_{i}\).

This gives us \(N\) segments of curves. To prove that we have a
continuous curve on the whole interval \([\epsilon_1, \epsilon_2]\) we
need to verify that these \(N\) segments in fact join together to form
one long curve. What one needs to avoid is a situation as in
Figure~\ref{fig:CGL-example-uniqueness-issue}, where we have two
separate curves and our enclosures go from enclosing the bottom curve,
to enclosing the upper curve. To ensure that this situation does not
happen, it is enough to verify that
\begin{equation*}
  \inter{\mu}_{i}^{(e)} \times \inter{\gamma}_{i}^{(e)} \times \inter{\kappa}_{i}^{(e)}
  \subseteq
  \inter{\mu}_{i - 1}^{(u)} \times \inter{\gamma}_{i - 1}^{(u)} \times \inter{\kappa}_{i - 1}^{(u)}
\end{equation*}
holds for all \(i = 2, \dots, N\), i.e., the enclosure of existence
for the \(i\)th interval is included in the enclosure of uniqueness
for the preceding interval. Note that in
Figure~\ref{fig:CGL-example-uniqueness-issue}, this condition is not
satisfied for the middle parts. In our case, this condition is however
easily verified, and the \(N\) segments of curves can be joined
together into one curve.

In some cases it will however not be possible to directly verify this
condition, this happens further down on some of the branches. In that
case we divide the intervals \(\inter{\epsilon}_i\) into smaller
subintervals, this gives tighter enclosures for the existence, which
makes the condition easier to verify.

\begin{remark}
  It is possible to use higher order methods to enclose the branches,
  see e.g.\ the works by Sander and Wanner~\cite{Sander2016} as well
  as Lessard, Sander and Wanner~\cite{Wanner2017}. While in general
  more efficient, these methods are based on bounds of the
  second-order derivatives, which in our case would require
  significant amount of work to get.
\end{remark}

The last remaining part is to count the number of critical points of
the profile \(|Q|\) along the curve. This is done for each box
independently, and the procedure is the same as for the example with
the NLS equation in Section~\ref{sec:detailed-example-nls}. For the
beginning of the branch, the process works fine for all subintervals,
but in some cases further down the branch the error bounds are too large
and the process sometimes fails. In these cases we split each
subinterval into smaller pieces, allowing us to get tighter
enclosures, and increasing the chance of success. This extra splitting
comes with a computational cost, but does allow us to count the number
of critical points along the whole branch.

\subsubsection{Verification of the full branch}
The verification of the full branch follows the same procedure as
above, but there are, however, a few more things that need to be taken
into account.

\begin{figure}
  \centering
  \includegraphics[width=0.45\textwidth]{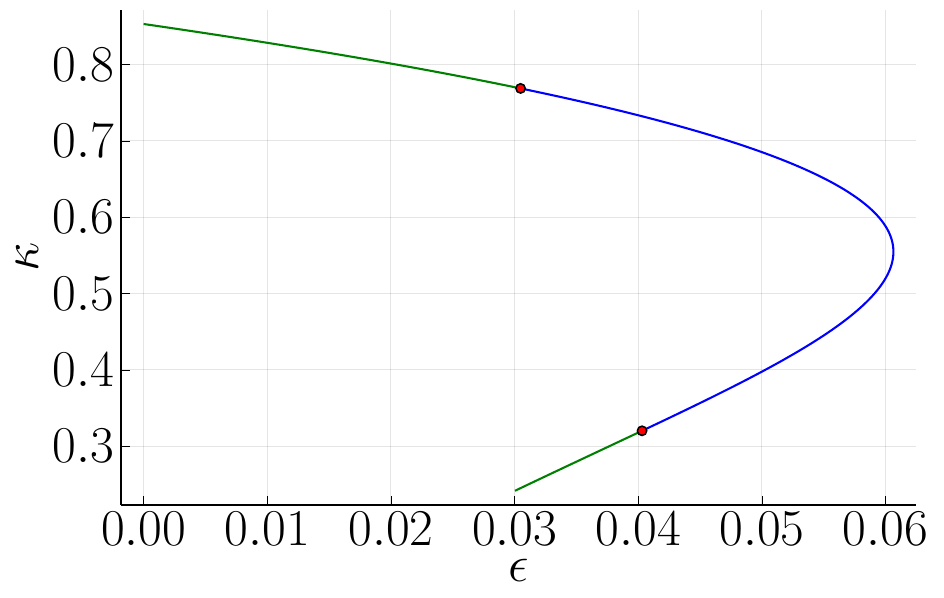}
  \caption{The branch is split into three parts, top/bottom (green)
    turning (blue).}
  \label{fig:CGL-example-parts}
\end{figure}

At the beginning of the branch it was possible to parameterize the
curve using \(\epsilon\). Once we reach the turning point of the
branch this is no longer possible, and a different parameterization is
required. For this, we split the curve into three separate parts, one
top part where the parameterization is done in \(\epsilon\), one
turning part where the parameterization is done in \(\kappa\) and one
bottom part where the parameterization is again done in \(\epsilon\).
The top, turning and bottom parts of the branch are shown in
Figure~\ref{fig:CGL-example-parts}. In practice we take the parts so
that they have a bit of an overlap.

\begin{remark}
  It would be possible to parameterize the whole curve by \(\kappa\).
  However, the parameterization in \(\epsilon\) performs
  computationally better near the beginning and end of the curve, and
  splitting the curve in these three parts is therefore beneficial.
  Exactly where the splitting is done is not so important, there are
  large regions where the parameterization in \(\epsilon\) and
  \(\kappa\) have more or less equal performance. In our case we make
  the top split roughly halfway between the beginning of the curve and
  the turning point, and the bottom split halfway between the turning
  point and \(\epsilon = 0.02\).
\end{remark}

For the top and bottom parts the procedure is exactly as described
above. For the turning part we have to flip the role of \(\epsilon\)
and \(\kappa\), instead of splitting it into intervals in
\(\epsilon\), we split it in \(\kappa\) as
\begin{equation*}
  \bigcup_{i = 1}^{N} \inter{\kappa}_{i}.
\end{equation*}
When applying the interval Newton method on \(G\), we now fix
\(\inter{\kappa}_i\) and look for a zero in \(\mu\), \(\gamma\) and \(\epsilon\)
instead. The main difference this gives is that we have to compute the
Jacobian of \(G\) with respect to \(\mu\), \(\gamma\) and
\(\epsilon\), instead of with respect to \(\mu\), \(\gamma\) and
\(\kappa\). Other than that, the procedure is the same.

With this, we get three continuous curves, one for the top part, one
for the turning part and one for the bottom part. The only thing that
remains is to prove that these three curves can be joined together
into one curve. In this case we cannot use the approach described
above, that requires the curves to be parameterized by the same
variable. Instead, we will prove that the curves can be connected by
finding a point that is proved to be contained in both curves, in
which case the two curves have to be joined together.

\begin{figure}
  \centering
  \begin{subfigure}[t]{0.45\textwidth}
    \includegraphics[width=\textwidth]{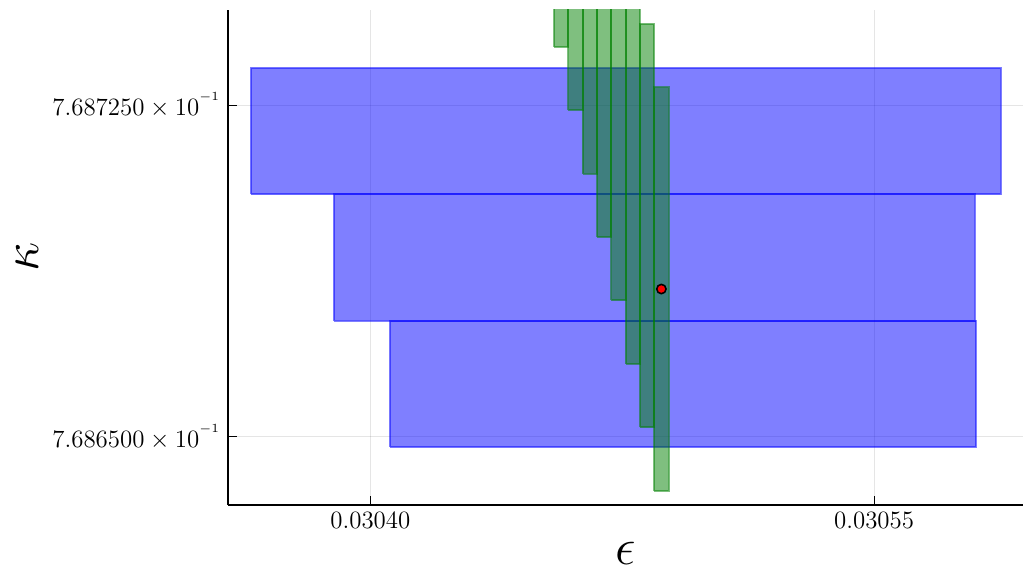}
    \caption{}
    \label{fig:CGL-example-top-turn}
  \end{subfigure}
  \hspace{0.05\textwidth}
  \begin{subfigure}[t]{0.45\textwidth}
    \includegraphics[width=\textwidth]{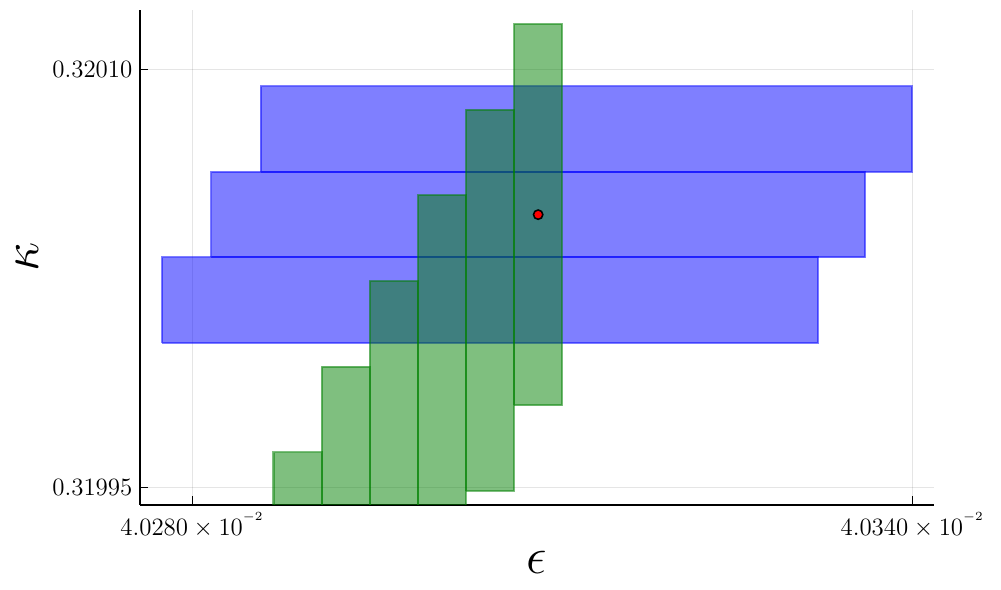}
    \caption{}
    \label{fig:CGL-example-turn-bottom}
  \end{subfigure}
  \caption{In the left figure, the red point lies on both the top and
    turning parts of the curve. In the right figure, it lies on the
    turning and bottom parts. The green boxes indicate the enclosures
    of uniqueness for the top and bottom parts of the curve, while the
    blue boxes show the enclosures of uniqueness for the turning
    part.}
  \label{fig:CGL-example-joining}
\end{figure}

This process of joining the top part and the turning part, is depicted
in Figure~\ref{fig:CGL-example-top-turn}. The green boxes represent
the uniqueness for the top part and the blue boxes the uniqueness for
the turning part, the red point is a very accurately computed
enclosure of one solution. The red point lies within the boxes of
uniqueness for the top part, and hence has to be part of the top
curve. Since it also lies within the box of uniqueness for the turning
part, it has to also be part of the turning curve. That means that
these two curves share a common point, and can hence be joined
together. For connecting the turning part and the bottom part the
procedure is the same, and depicted in
Figure~\ref{fig:CGL-example-turn-bottom}.

Verifying the whole branch is computationally costly and in general
the cost increases further along the branch.
Table~\ref{table:CGL-example} gives some statistics about the computational
costs for the top, turning and bottom parts of the branch, the number
of subintervals that they are split into, the computational time to
verify the existence of the curve and to count the number of critical
points. The bottom part takes the longest time, and, compared to its
length, it has the largest number of subintervals. In fact, to avoid
an excessive computational cost we stop the verification at an earlier
stage. The numerical approximation goes on until
\(\epsilon \approx 0.02\), see
Figure~\ref{fig:CGL-example-approximation}, whereas we stop the
verification at \(\epsilon \approx 0.03\), see
Figure~\ref{fig:CGL-example-parts}.

\begin{table}[ht]
  \centering
  \begin{tabular}{lcccc}
    \toprule
    & Top & Turn & Bottom\\
    \midrule
    Subintervals & 6730 & 21248 & 7641\\
    Subintervals & 6157 & 21216 & 7082\\
    Runtime existence of curve (seconds) & 297 & 1787 & 701\\
    Runtime existence of curve (core hours) & 11 & 64 & 25\\
    Runtime critical points (seconds) & 166 & 27175 & 5039\\
    Runtime critical points (core hours) & 6 & 966 & 179\\
    \bottomrule
  \end{tabular}
  \caption{Data about the number of subintervals and the computational
    time for the top, turning and bottom parts of the first branch in
    \textbf{Case I}. The computations were done on a \(128\) core
    machine, and runtime in core hours is the runtime in seconds
    multiplied by \(\frac{128}{3600}\).}
  \label{table:CGL-example}
\end{table}

\subsection{Other cases}
\label{sec:other-cases-cgl}
In \textbf{Case I} the verification of the remaining 7 branches
follows exactly the same approach as in the above example. Some of
them, however, come at a significantly higher computational cost, in
particular for verifying the number of critical points. For the
branches \(j = 7\) and \(j = 8\) we only verify the top part of the
branch, as the verification of the turning and bottom parts becomes
prohibitively costly. For most of the other branches we also stop the
verification before reaching \(\epsilon \approx 0.02\) (which is how
far the numerical approximations are computed), again to avoid an
excessive computational time.

In \textbf{Case II} the situation is even more delicate. Here the
beginning of the branches is the hardest to compute, but the
other parts also require significant effort (this is discussed in a
bit more detail in Section~\ref{sec:pointwise-results-cgl}). For this
reason, we only verify some segments of the branches. Other than the
significantly higher computational cost, the procedure is generally
the same as in the above example. For \(j = 3\) one has to switch to
the parameterization in \(\epsilon\) before it makes the turn in
\(\kappa\). For \(j = 5\) one would have to split it into more than
three parts to be able to handle the larger number of turns, though we
do not attempt to verify that part.

\begin{figure}[htbp]
  \centering
  \begin{subfigure}[t]{0.32\textwidth}
    \includegraphics[width=\textwidth]{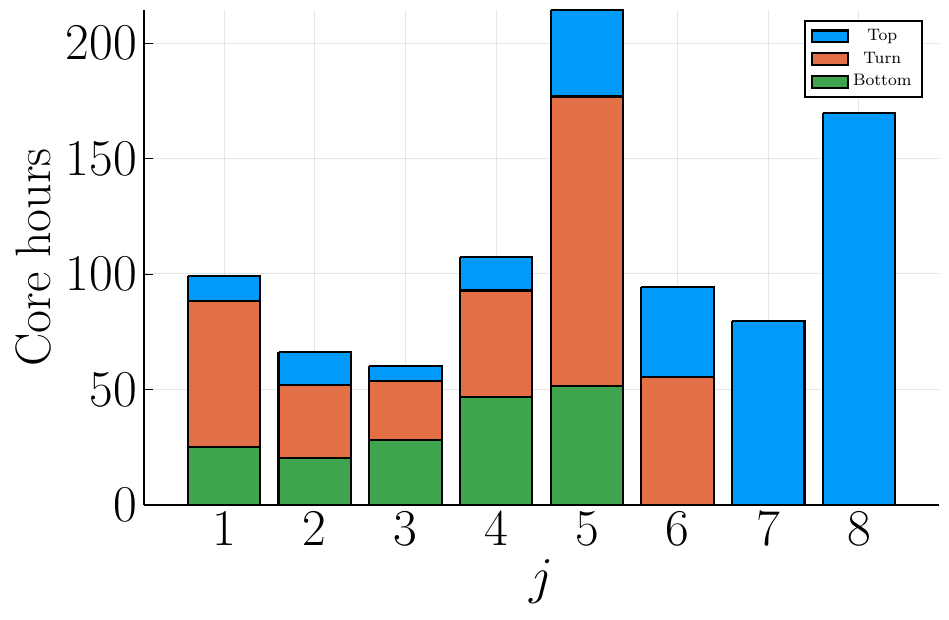}
    \caption{Proving existence (\textbf{Case I}).}
    \label{fig:CGL-runtime-d=1}
  \end{subfigure}
  \hfill
  \begin{subfigure}[t]{0.32\textwidth}
    \includegraphics[width=\textwidth]{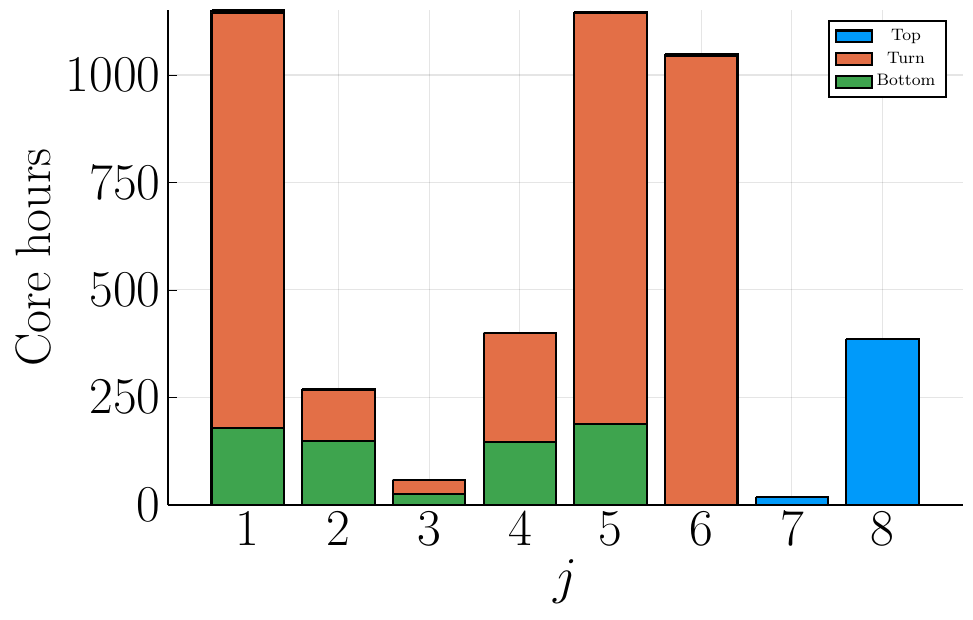}
    \caption{Critical points (\textbf{Case I}).}
    \label{fig:CGL-runtime-critical-points-d=1}
  \end{subfigure}
  \hfill
  \begin{subfigure}[t]{0.32\textwidth}
    \includegraphics[width=\textwidth]{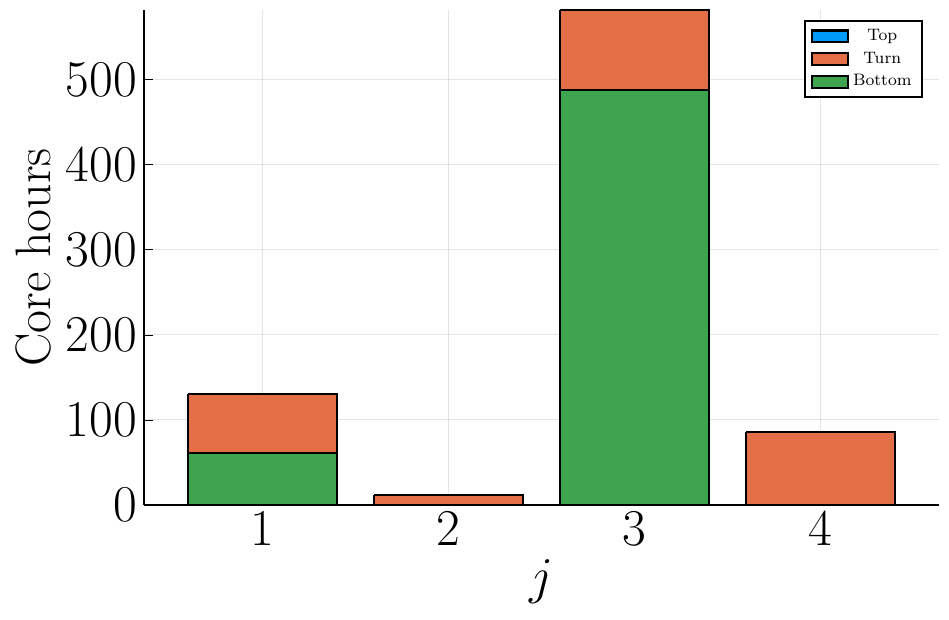}
    \caption{Proving existence (\textbf{Case II}).}
    \label{fig:CGL-runtime-d=3}
  \end{subfigure}
  \caption{Runtimes for \textbf{Case I} (8 branches) and \textbf{Case II} (4 branches).}
  \label{fig:overall-runtime}
\end{figure}

The computational times for proving the existence of the 8 branches in
\textbf{Case I} are given in Figure~\ref{fig:CGL-runtime-d=1}, and the
time to count the number of critical points for \(|Q|\) along the
branches is given in Figure~\ref{fig:CGL-runtime-critical-points-d=1}.
For \textbf{Case II} the times for proving existence are given in
Figure~\ref{fig:CGL-runtime-d=3}. The runtimes are reported in terms
of core hours; this is the wall clock time for running the
computations multiplied by the number of CPU cores that were used (in
this case 128).

\begin{remark}[Computational time]
  The portion of the branches covered was scaled to allow all
  computations to run comfortably over a weekend. With additional
  computational time, one could cover larger portions of the branches.
  Alternatively, the computations could be significantly shortened by,
  for example, only covering the top parts of the branches in
  \textbf{Case I}.

  Note that only the computation of the branches incurs a high
  computational cost. The proofs for the NLS equation in
  Section~\ref{sec:proof-nls} run in a matter of minutes on a regular
  laptop.
\end{remark}

\section{Solution at infinity}
\label{sec:solution-infinity}
In this section, we study solutions to the equation
\begin{equation}\label{eq:Q-infty-equation}
  (1 - i\epsilon)\left(Q'' + \frac{d - 1}{\xi}Q'\right) + i\kappa\xi Q'
  + i \frac{\kappa}{\sigma}Q - \omega Q + (1 + i\delta)|Q|^{2\sigma}Q = 0,
\end{equation}
satisfying the condition
\(Q(\xi) \sim \xi^{-\frac{1}{\sigma} - i\frac{\omega}{\kappa}}\) as
\(\xi \to \infty\). We will see that there exists a manifold of
solutions, which can be parameterized by \(\gamma \in \mathbb{C}\).

To perform the shooting method discussed in
Section~\ref{sec:shooting-method}, for a fixed \(\xi_{1} > 1\), we
need to be able to compute enclosures of \(Q(\xi_{1})\) and
\(Q'(\xi_{1})\), as well as the derivatives with respect to the
parameters \(\gamma\), \(\kappa\) and \(\epsilon\). Since the
non-linearity \(|Q|^{2\sigma}Q\) is not holomorphic, we need to
separate derivatives in the real and imaginary directions. For the
derivatives with respect to \(\gamma\) we use the notation
\(Q_{\real\gamma}\) and \(Q_{\imag\gamma}\) for the derivatives in the
real and imaginary directions respectively. Since \(\kappa\) and
\(\epsilon\) are real, there is only one direction to consider, and we
write \(Q_{\kappa}\) and \(Q_{\epsilon}\) for the corresponding
derivatives. To count the number of critical points of \(|Q|\), we
also need to prove that the solution is monotone for sufficiently
large \(\xi\).

To compute the solution, we follow a similar approach to that
in~\cite{Plech2001}. The main idea is to use the variation of
parameters applied to the linearized equation to construct an operator
whose fixed points are solutions to the equation. All the necessary
bounds are then derived from this fixed point equation.

The associated linear equation is given by
\begin{equation}\label{eq:Q-infty-linear}
  (1 - i\epsilon)\left(Q'' + \frac{d - 1}{\xi}Q'\right) + i\kappa\xi Q' + i \frac{\kappa}{\sigma}Q - \omega Q = 0.
\end{equation}
Making the change of variables
\begin{equation}
  \label{eq:abc}
  a = \frac{1}{2}\left(\frac{1}{\sigma} + i \frac{\omega}{\kappa}\right),\quad
  b = \frac{d}{2},\quad
  c = \frac{-i \kappa}{2(1 - i\epsilon)},\quad
  z = c\xi^{2},
\end{equation}
this can be transformed into the well-studied Kummer's equation,
\begin{equation*}
  z \frac{d^{2}w}{dz^{2}} + (b - z)\frac{dw}{dz} - aw = 0.
\end{equation*}
One solution of Kummer's equation is the confluent hypergeometric
function \(U(a, b, z)\), with another, linearly independent, solution
given by \(V(a, b, z) = e^{z}U(b - a, b, -z)\). This gives us two
linearly independent solutions \(P(\xi)\) and \(E(\xi)\)
of~\eqref{eq:Q-infty-linear}, given by
\begin{equation*}
  P(\xi) = U(a, b, c\xi^{2}) \quad\text{and}\quad
  E(\xi) = e^{c\xi^{2}}U(b - a, b, -c\xi^{2}).
\end{equation*}
We also note that the associated Wronskian is given by
\begin{equation*}
  W(\xi) = P(\xi)E'(\xi) - P'(\xi)E(\xi) = 2ce^{\sign(\imag(c)) \pi i (b - a)}\xi \left(c\xi^{2}\right)^{-b}e^{c\xi^{2}}.
\end{equation*}
Of the functions \(P(\xi)\) and \(E(\xi)\), only \(P(\xi)\) has the
right asymptotic behavior at infinity. More precisely, \(P(\xi)\)
asymptotically behaves like
\(P(\xi) \sim \xi^{-\frac{1}{\sigma} - i \frac{\omega}{\kappa}}\) (see
Lemma~\ref{lemma:U}), whereas \(E(\xi)\) behaves like
\(E(\xi) \sim e^{c\xi^{2}}\xi^{\frac{1}{\sigma} - d +
  i\frac{\omega}{\kappa}}\). For \(\epsilon > 0\) we have
\(\real c > 0\), so \(E(\xi)\) is exponentially increasing. However,
even for \(\epsilon = 0\), the large oscillations coming from the
exponential factor stop \(E(\xi)\) from having finite energy.

To go from solutions to the linear equation to solutions
to~\eqref{eq:Q-infty-equation}, we follow the method of variation of
parameters and look for solutions of the form
\(Q(\xi) = c_{1}(\xi)P(\xi) + c_{2}(\xi)E(\xi)\), where
\begin{equation*}
  c_{1}'(\xi) = \frac{1 + i\delta}{1 - i \epsilon}E(\xi)W(\xi)^{-1}|Q(\xi)|^{2\sigma}Q(\xi)
  \quad\text{and}\quad
  c_{2}'(\xi) = -\frac{1 + i\delta}{1 - i \epsilon}P(\xi)W(\xi)^{-1}|Q(\xi)|^{2\sigma}Q(\xi).
\end{equation*}
To get the right asymptotic behavior at infinity, we require that
\(c_{2}(\xi)\) goes to zero rapidly as \(\xi \to \infty\), to
counteract the exponential factor in \(E(\xi)\). In contrast,
\(c_{1}\) should have a nonzero limit at infinity. For this reason, it
is natural to use the representations
\begin{align*}
  c_{1}(\xi) &= \gamma + \int_{\xi_{1}}^{\xi}\frac{1 + i\delta}{1 - i \epsilon}E(\eta)W(\eta)^{-1}|Q(\eta)|^{2\sigma}Q(\eta)\,d\eta,\\
  c_{2}(\xi) &= \int_{\xi}^{\infty}\frac{1 + i\delta}{1 - i \epsilon}P(\eta)W(\eta)^{-1}|Q(\eta)|^{2\sigma}Q(\eta)\,d\eta.
\end{align*}
The vanishing interval of integration for \(c_{2}\) ensures that it
goes to zero in the limit. Using the notation
\begin{equation*}
  J_{P}(\xi) = \frac{1 + i\delta}{1 - i \epsilon}P(\xi)W(\xi)^{-1},\quad
  J_{E}(\xi) = \frac{1 + i\delta}{1 - i \epsilon}E(\xi)W(\xi)^{-1}
\end{equation*}
and
\begin{equation*}
  I_{P}(Q, \xi) = \int_{\xi}^{\infty} J_{P}(\eta)|Q(\eta)|^{2\sigma}Q(\eta)\,d\eta,\quad
  I_{E}(Q, \xi) = \int_{\xi_{1}}^{\xi} J_{E}(\eta)|Q(\eta)|^{2\sigma}Q(\eta)\,d\eta,
\end{equation*}
this gives us the equation
\begin{equation}\label{eq:Q-infty-fixed-point-equation}
  Q(\xi) = \gamma P(\xi) + P(\xi)I_{E}(Q, \xi) + E(\xi)I_{P}(Q, \xi)
\end{equation}
for \(Q\). Specifically, if we let \(T\) be the operator
\begin{equation}\label{eq:T}
  T(Q)(\xi) = \gamma P(\xi) + P(\xi)I_{E}(Q, \xi) + E(\xi)I_{P}(Q, \xi),
\end{equation}
then fixed points of this operator correspond to solutions
of~\eqref{eq:Q-infty-equation}.

Our strategy for computing enclosures of \(Q\) is to study the
operator \(T\) and the fixed point
equation~\eqref{eq:Q-infty-fixed-point-equation}. We split this into
two parts:
\begin{enumerate}
\item First, we prove the existence of a (locally unique) fixed point
  of \(T\) by proving that it is a contraction in a ball in a weighted
  function space. This gives us both the existence of a solution, as
  well as initial bounds for its norm.
\item Next, we compute refined enclosures of the fixed point using a
  bootstrapping approach applied to the fixed point
  equation~\eqref{eq:Q-infty-fixed-point-equation}. For the shooting
  argument to be successful, these enclosures need to be sufficiently
  tight.
\end{enumerate}
Both of these steps require precise control over the solutions \(P\)
and \(E\) to the linear equation. It is then relatively
straightforward to prove the existence of a fixed point. Getting
sufficiently refined enclosures for the shooting argument to succeed
does, however, require significantly more work. The main culprit is
getting good control over the integral \(I_{P}(Q, \xi)\), whose
integrand has exponentially fast oscillations coming from \(J_{P}\).
To handle this, we integrate by parts multiple times, giving us
progressively better control. This does, however, require the use of
higher order derivatives of \(P\), \(E\) and \(Q\). A significant
amount of the work is therefore dedicated to also handling these
higher order derivatives.

To control the derivatives with respect to \(\real\gamma\),
\(\imag\gamma\), \(\kappa\) and \(\epsilon\), we directly
differentiate the fixed-point equation. To get tight enclosures, we
then follow a similar approach to that for enclosing \(Q\) and \(Q'\).
In particular, this requires precise control over the derivatives of
\(P\) and \(E\) with respect to the parameters \(\kappa\) and
\(\epsilon\).

We proceed as follows:
\begin{enumerate}
\item In Section~\ref{sec:function-bounds}, we give asymptotic bounds
  for functions related to \(P\) and \(E\). This includes both higher
  order derivatives with respect to \(\xi\) (required for the
  integration by parts of \(I_{P}\)) and derivatives with respect to
  \(\kappa\) and \(\epsilon\).
\item Section~\ref{sec:fixed-point} is devoted to studying the
  operator \(T\), and we give explicit conditions for the existence of
  a fixed point.
\item In Section~\ref{sec:fixed-point-enclosures}, we discuss how to
  use the fixed point equation~\eqref{eq:Q-infty-fixed-point-equation}
  to compute tight enclosures of \(Q\) and its derivatives.
\item Finally, Section~\ref{sec:verify-monotonicity} goes through the
  extra details required to prove monotonicity of \(Q\) for large
  values of \(\xi\).
\end{enumerate}

Lastly, let us comment on some of the notation and other conventions
used in this section. In general, we use \('\) to denote
differentiation with respect to the argument, e.g.\ \(P'(\xi)\),
and subscripts to denote differentiation with respect to parameters,
e.g.\ \(J_{P,\kappa}(\xi)\). For the most part, we do not make explicit
the functions' dependence on parameters such as \(\kappa\) and
\(\epsilon\). In the same way, we generally drop the \(Q\) argument
for \(I_{P}\) and \(I_{E}\) defined above, i.e.\ we prefer to write
just \(I_{P}(\xi)\) instead of \(I_{P}(Q, \xi)\). The sections contain
a large number of asymptotic bounds with explicitly given constants.
For example, from Lemma~\ref{lemma:bounds-list} we have the bound
\begin{equation*}
  |P(\xi)| \leq C_{P}\xi^{-\frac{1}{\sigma}},
\end{equation*}
valid for \(\xi \geq \xi_{1}\) and where \(C_{P}\) is an explicitly
given constant. In general, these explicit constants will depend on
parameters such as \(\kappa\) and \(\epsilon\), as well as the value
of \(\xi_{1}\). However, to keep the notation brief, we do not make
this dependence explicit. We use the variable \(C\) for all of these
types of constants, with a subscript to denote which function it is
associated with. For constants associated with derivatives with
respect to parameters, we generally prefer to put the parameter in the
top level of the subscript. For example, we use the notation
\(C_{P,\epsilon}\), rather than \(C_{P_{\epsilon}}\), for a constant
associated with \(P_{\epsilon}\). In some cases, there are multiple
constants appearing in the bound for a function; in this case, we opt
to use an integer subscript to enumerate them. An example is the bound
for \(Q''\), which depends on the two constants \(C_{Q'',1}\) and
\(C_{Q'',2}\).

\subsection{Asymptotic bounds for \(P\) and \(E\)}
\label{sec:function-bounds}
In this subsection we give asymptotic bounds for \(P\) and \(E\) as
well as several related functions.

We start by giving bounds for the confluent hypergeometric function
\(U\), all other bounds are then derived from these bounds. Note that
differentiation of \(U\) with respect to the argument gives us another
confluent hypergeometric function, given by
\begin{equation}
  \label{eq:U-dz}
  U^{(n)}(a, b, z) = (-1)^{n}U(a + n, b + n, z)(a)_{n},
\end{equation}
where \((a)_{n} = a(a + 1)(a + 2)\cdots(a + n - 1)\) is the rising
factorial. As a consequence, having bounds for \(U\) also allows us to
get bounds for derivatives with respect to \(z\). The following lemma
gives the asymptotic expansion of \(U\) as well as an explicit bound
for the remainder, both of which are available in the literature.

\begin{lemma}
  \label{lemma:U}
  Let \(a, b, z \in \mathbb{C}\) with \(|\imag z| > |b - 2a|\) or
  \(|\real z| > |b - 2a|\). For \(n \in \mathbb{Z}_{\geq 0}\) we have
  \begin{equation*}
    U(a, b, z) = \left(\sum_{k = 0}^{n - 1} \frac{(a)_{k}(a - b + 1)_{k}}{k!(-z)^{k}} + R_{U}(a, b, n, z)z^{-n}\right)z^{-a},
  \end{equation*}
  where
  \begin{equation*}
    \left|R_{U}\!\left(a,b,n,z\right)\right| \leq
    C_{R_{U}}(a, b, n, z) =
    \left|\frac{\left(a\right)_{n} \left(a - b + 1\right)_{n}}{n!}\right|
    \frac{2 \sqrt{1 + \frac{1}{2} \pi n}}{1 - \sigma(a, b, z)}
    \exp\!\left(\frac{\pi \rho(a, b, z)}{\left(1 - \sigma(a, b, z)\right) \left|z\right|}\right)
  \end{equation*}
  with
  \begin{equation*}
    \sigma(a, b, z) = \frac{\left|b - 2 a\right|}{\left|z\right|},\quad
    \rho(a, b, z) = \left|{a}^{2} - a b + \frac{b}{2}\right| + \frac{\sigma(a, b, z) \left(1 + \frac{\sigma(a, b, z)}{4}\right)}{{\left(1 - \sigma(a, b, z)\right)}^{2}}.
  \end{equation*}

  In particular, if \(z_1 \in \mathbb{C}\) satisfies
  \(|\imag z_1| > |b - 2a|\) or \(|\real z_1| > |b - 2a|\), then for
  all \(z \in \mathbb{C}\) such that \(\arg(z) = \arg(z_{1})\) and
  \(|z| \geq |z_{1}|\), we have
  \begin{equation*}
    C_{R_{U}}(a, b, n, z) \leq C_{R_{U}}(a, b, n, z_{1})
  \end{equation*}
  and, moreover,
  \begin{equation*}
    |U(a, b, z)| \leq C_{U}(a,b,n,z_{1})|z^{-a}|,
  \end{equation*}
  with \(C_{U}(a,b,n,z_{1})\) given by
  \begin{equation*}
    C_{U}(a,b,n,z_{1}) = \sum_{k = 0}^{n - 1} \left|\frac{(a)_{k}(a - b + 1)_{k}}{k!z_1^{k}}\right| +
    C_{R_{U}}(a, b, n, z_{1})|z_{1}^{-n}|.
  \end{equation*}
\end{lemma}

\begin{proof}
  The first statement follows from~\cite[Sec.~13.7]{NIST:DLMF} (see
  also entries~\cite{fungrim-c8fcc7,fungrim-d1b3b5,fungrim-461a54}
  from Fungrim~\cite{Johansson2020fungrim}).

  For the second statement, the bound
  \(C_{R_{U}}(a, b, n, z) \leq C_{R_{U}}(a, b, n, z_{1})\) follows
  immediately from the fact that \(C_{R_{U}}\) is decreasing in
  \(|z|\). Since the condition \(\arg(z) = \arg(z_{1})\) implies that
  \(z\) also satisfies \(|\imag z| > |b - 2a|\) or
  \(|\real z| > |b - 2a|\), the expression for the constant \(C_{U}\)
  then follows from bounding the asymptotic expansion for
  \(U(a, b, z)\) termwise and applying the monotonicity of the error
  bound.
\end{proof}

For our purposes we also need bounds for derivatives with respect to
\(a\), i.e. \(U_{a}\). In this case there are no bounds for the
remainder terms directly available in the literature. Instead, we
adapt the approach in~\cite{Olver1991b} to cover also this case.

\begin{lemma}
  \label{lemma:U-a}
  Let \(a, b, z \in \mathbb{C}\), satisfying \(|\imag z| > |b - 2a|\)
  or \(|\real z| > |b - 2a|\) as well as \(0 < \real(a) < \real(b)\)
  and \(|\arg(z)| < \pi\). For \(n \in \mathbb{Z}_{\geq 0}\) such that
  \(\real(a - b + n + 1) > 0\) we have
  \begin{align*}
    U_{a}(a, b, z)
    &= \Bigg(
      -\sum_{k = 0}^{n - 1}\frac{(a)_{k}(a - b + 1)_{k}}{k!(-z)^{k}}
      + \log(z)^{-1}\sum_{k = 0}^{n - 1}\frac{\partial}{\partial a}\left(
      \frac{(a)_{k}(a - b + 1)_{k}}{k!(-z)^{k}}
      \right)\\
    &\qquad+ \left(
      \left(1 - \frac{\Gamma'(a)}{\Gamma(a)}\log(z)^{-1}\right)R_{U}(a, b, n, z)
      + R_{U,1}(a, b, z) + R_{U,2}(a, b, z)
      \right)z^{-n}
      \Bigg)\log(z)z^{-a}.
  \end{align*}
  Here \(R_{U}\) is as in Lemma~\ref{lemma:U},
  \begin{multline*}
    \left|R_{U,1}\!\left(a,b,n,z\right)\right| \leq
    C_{R_{U,1}}(a, b, n, z) =
    \frac{\rho_{\gamma}C_{\chi}(a, b, n)e^{-\gamma\imag(a)}}{|\Gamma(a)|}\\
    \left(
      C_{1}(a, n) + (|\gamma|C_{1}(a, n) + C_{2}(a, n))|\log(|z|)|^{-1}
    \right)
  \end{multline*}
  and
  \begin{equation*}
    \left|R_{U,2}\!\left(a,b,n,z\right)\right| \leq
    C_{R_{U,2}}(a, b, n, z) =
    \frac{\rho_{\gamma}C_{\chi,a}(a, b, n)C_{1}(a, n)e^{-\gamma\imag(a)}}{|\Gamma(a)|}|\log(|z|)|^{-1},
  \end{equation*}
  with \(\gamma = \arg(z)\),
  \begin{align*}
    C_{\chi}(a, b, n)
    &= \frac{|\sin(\pi(a - b + 1))|}{\pi}
      \frac{\Gamma(-\real(a - b))\Gamma(\real(a - b + n + 1))}{n!},\\
    C_{\chi,a}(a, b, n)
    &= |\cos(\pi(a - b + 1))|
      \frac{\Gamma(-\real(a - b))\Gamma(\real(a - b + n + 1))}{n!}\\
    &\qquad+ \frac{|\sin(\pi(a - b + 1))|}{\pi}
      \left(
      \frac{1}{\real(a - b)^{2}}
      + \frac{\Gamma(-\real(a - b))\Gamma(\real(a - b + n))}{(n - 1)!}
      \right),\\
    C_{1}(a, n)
    &= \Gamma(\real(a + n)),\\
    C_{2}(a, n)
    &= \frac{1}{\real(a + n)^{2}} + \Gamma(\real(a + n + 1))
  \end{align*}
  and
  \begin{equation*}
    \rho_\gamma =
    \begin{cases}
      1 & \text{if } |\gamma| \leq \frac{\pi}{2}, \\
      \sin(\pi - |\gamma|)^{-1} & \text{if } \frac{\pi}{2} < |\gamma| < \pi.
    \end{cases}
  \end{equation*}

  In particular, if \(z_1 \in \mathbb{C}\), with \(|z_{1}| > 1\),
  satisfies the requirements on \(z\) above, then for all
  \(z \in \mathbb{C}\) such that \(\arg(z) = \arg(z_{1})\) and
  \(|z| \geq |z_1|\), we have
  \begin{equation*}
    |U_a(a, b, z)| \leq C_{U,a}(a,b,n,z_{1})|\log(z)z^{-a}|
  \end{equation*}
  with \(C_{U,a}(a,b,n,z_{1})\) given by
  \begin{align*}
    C_{U,a}(a,b,n,z_{1})
    &= \sum_{k = 0}^{n - 1}\left|\frac{(a)_{k}(a - b + 1)_{k}}{k!(-z_{1})^{k}}\right|
      + |\log(z_{1})|^{-1}\sum_{k = 0}^{n - 1}\left|\frac{\partial}{\partial a}\left(
      \frac{(a)_{k}(a - b + 1)_{k}}{k!(-z_{1})^{k}}
      \right)\right|\\
    &\qquad+ \left(
      \left(1 + \left|\frac{\Gamma'(a)}{\Gamma(a)}\log(z_{1})^{-1}\right|\right)R_{U}(a, b, n, z_{1})
      + R_{U,1}(a, b, z_{1}) + R_{U,2}(a, b, z_{1})
      \right)|z_{1}|^{-n}.
  \end{align*}
\end{lemma}

\begin{proof}
  We follow the methodology established by Olver~\cite{Olver1991b},
  specifically the techniques detailed in Sections 1 and 2. Note that
  they prove stronger results than we intend to do here; we therefore
  only need part of the paper.

  The idea is to take the asymptotic expansion of \(U(a, b, z)\) from
  Lemma~\ref{lemma:U} above and differentiate it termwise. However, to
  get a bound for the remainder term we first need an explicit formula
  for it. The explicit formula we use is~\cite[Equation
  2.5]{Olver1991b} (note that their parameter \(b\) corresponds to
  \(a - b + 1\) in our case). For completeness, we include a brief
  derivation.

  The starting point is the integral representation
  \begin{equation*}
    U(a, b, z) = \frac{1}{\Gamma(a)} \int_{0}^{\infty}e^{-zt}t^{a - 1}(1 + t)^{-(a - b + 1)}\,dt.
  \end{equation*}
  However, this integral only converges when \(\real(a) > 0\) and
  \(\real(z) > 0\). The requirement on \(a\) is satisfied by
  assumption, but not the one on \(z\). However, by analytic
  continuation one can, for an arbitrary angle
  \(\gamma \in (-\pi, \pi)\), get the representation
  \begin{equation*}
    U(a, b, z) = \frac{1}{\Gamma(a)} \int_{0}^{\infty e^{-i\gamma}}e^{-zt}t^{a - 1}(1 + t)^{-(a - b + 1)}\,dt,
  \end{equation*}
  which is valid for all \(z\) such that \(\real(ze^{-i\gamma}) > 0\),
  still requiring \(\real(a) > 0\).

  The next step is to use the Taylor expansion
  \begin{equation*}
    (1 + t)^{-(a - b + 1)} =
    \sum_{k = 0}^{n - 1}(-1)^{k}\frac{(a - b + 1)_{k}}{k!}t^{k} + \chi_{n}(a, b, t).
  \end{equation*}
  Inserting this into the integral gives us
  \begin{equation*}
    U(a, b, z) = \sum_{k = 0}^{n - 1}\frac{(a - b + 1)_{k}}{k!}
    \frac{(-1)^{k}}{\Gamma(a)}
    \int_{0}^{\infty e^{-i\gamma}}e^{-zt}t^{k + a - 1}\,dt
    + \frac{1}{\Gamma(a)} \int_{0}^{\infty e^{-i\gamma}}e^{-zt}t^{a - 1}\chi_{n}(a, b, t)\,dt.
  \end{equation*}
  The identity
  \begin{equation*}
    \frac{(-1)^{k}}{\Gamma(a)}
    \int_{0}^{\infty e^{-i\gamma}}e^{-zt}t^{k + a - 1}\,dt
    = \frac{(a)_{k}}{(-z)^{k}}z^{-a}
  \end{equation*}
  gives us the asymptotic expansion in Lemma~\ref{lemma:U}. The
  remainder term \(R_{U}(a, b, z)\) is then given by
  \begin{equation*}
    R_{U}(a, b, z) = \frac{z^{n + a}}{\Gamma(a)} \int_{0}^{\infty e^{-i\gamma}}e^{-zt}t^{a - 1}\chi_{n}(a, b, t)\,dt.
  \end{equation*}

  We are interested in bounding
  \(\frac{\partial}{\partial a}R_U(a, b, z)\). Differentiating the
  integral representation with respect to \(a\) gives us three terms:
  one arising from the \(z^{n + a} / \Gamma(a)\) factor, and two from
  the derivatives of \(t^{a - 1}\) and \(\chi_{n}(a, b, t)\) in the
  integrand. If we let
  \begin{equation*}
    R_{U,1}(a, b, z) = \log(z)^{-1}\frac{z^{n + a}}{\Gamma(a)} \int_{0}^{\infty e^{-i\gamma}}e^{-zt}\log(t)t^{a - 1}\chi_{n}(a, b, t)\,dt
  \end{equation*}
  and
  \begin{equation*}
    R_{U,2}(a, b, z) = \log(z)^{-1}\frac{z^{n + a}}{\Gamma(a)} \int_{0}^{\infty e^{-i\gamma}}e^{-zt}t^{a - 1}\frac{\partial}{\partial a}\chi_{n}(a, b, t)\,dt
  \end{equation*}
  denote the two derivatives associated with \(t^{a - 1}\) and
  \(\chi_{n}(a, b, t)\) respectively, then we can write the full
  derivative as
  \begin{equation*}
    \frac{\partial}{\partial a}R_{U}(a, b, z)
    = \left(\log(z) - \frac{\Gamma'(a)}{\Gamma(a)}\right)R_{U}(a, b, z) + R_{U,1}(a, b, z) + R_{U,2}(a, b, z).
  \end{equation*}

  Before bounding \(R_{U,1}\) and \(R_{U,2}\), let us take a closer
  look at \(\chi_{n}(a, b, t)\). We start by noting that
  \begin{equation*}
    \chi_{n}(a, b, t) = \frac{t^{n}}{2\pi i}\int_{\mathcal{L}} \frac{(1 + w)^{-(a - b + 1)}}{w^{n}}(w - t)\,dw,
  \end{equation*}
  where \(\mathcal{L}\) is a simple closed contour enclosing \(w = 0\)
  and \(w = t\), but not \(w = -1\). This loop can be deformed into
  the one depicted in~\cite[Figure 1]{Olver1991b}. If
  \(\real(a) < \real(b)\), which holds by assumption, then it can be
  collapsed on to the interval \((-\infty, -1]\). For this we need to
  take into account the branch cut along \((-\infty, -1]\), which adds
  a factor
  \begin{equation*}
    e^{-i(a - b + 1)\pi} - e^{-i(a - b + 1)(-\pi)} = -2i\sin(\pi(a - b + 1)).
  \end{equation*}
  This gives us
  \begin{equation*}
    \chi_{n}(a, b, t) = -\frac{\sin(\pi(a - b + 1))}{\pi}t^{n}\int_{-1}^{-\infty} \frac{(1 + w)^{-(a - b + 1)}}{w^{n}(w - t)}\,dw,
  \end{equation*}
  which after a change of variables from \(w\) to \(-w\) gives us
  \begin{equation*}
    \chi_{n}(a, b, t) = (-1)^{n}\frac{\sin(\pi(a - b + 1))}{\pi}t^{n}\int_{1}^{\infty} \frac{(w - 1)^{-(a - b + 1)}}{w^{n}(w + t)}\,dw.
  \end{equation*}
  We have
  \begin{equation*}
    |\chi_{n}(a, b, t)| \leq \frac{|\sin(\pi(a - b + 1))|}{\pi}|t|^{n}\int_{1}^{\infty} \frac{(w - 1)^{-\real(a - b + 1)}}{w^{n + 1}|1 + t / w|}\,dw.
  \end{equation*}
  Some trigonometry gives us that the factor \(|1 + t / w|^{-1}\) is
  uniformly bounded by
  \begin{equation*}
    \rho_\gamma =
    \begin{cases}
      1 & \text{if } |\gamma| \leq \frac{\pi}{2}, \\
      \sin(\pi - |\gamma|)^{-1} & \text{if } \frac{\pi}{2} < |\gamma| < \pi.
    \end{cases}
  \end{equation*}
  This gives us \(|\chi_{n}(a, b, t)| \leq C_{\chi}(a, b, n) |t|^n\),
  with
  \begin{align*}
    C_{\chi}(a, b, n)
    &= \frac{\rho_{\gamma}|\sin(\pi(a - b + 1))|}{\pi}
      \int_{1}^{\infty} \frac{(w - 1)^{-\real(a - b + 1)}}{w^{n + 1}}\,dw\\
    &= \frac{\rho_{\gamma}|\sin(\pi(a - b + 1))|}{\pi}
      \frac{\Gamma(-\real(a - b))\Gamma(\real(a - b + n + 1))}{n!}.
  \end{align*}
  Furthermore
  \begin{multline*}
    \frac{\partial}{\partial a}\chi_{n}(a, b, t)
    = (-1)^{n}\cos(\pi(a - b + 1))t^{n}
    \int_{1}^{\infty} \frac{(w - 1)^{-(a - b + 1)}}{w^{n}(w + t)}\,dw\\
    - (-1)^{n}\frac{\sin(\pi(a - b + 1))}{\pi}t^{n}
    \int_{1}^{\infty} \log(w - 1)\frac{(w - 1)^{-(a - b + 1)}}{w^{n}(w + t)}\,dw,
  \end{multline*}
  giving us
  \begin{equation*}
    \left|\frac{\partial}{\partial a}\chi_{n}(a, b, t)\right|
    \leq C_{\chi,a}(a, b, n)|t|^{n}
  \end{equation*}
  with
  \begin{multline*}
    C_{\chi,a}(a, b, n) = \rho_{\gamma}|\cos(\pi(a - b + 1))|
    \int_{1}^{\infty} \frac{(w - 1)^{-(a - b + 1)}}{w^{n + 1}}\,dw\\
    + \frac{\rho_{\gamma}|\sin(\pi(a - b + 1))|}{\pi}
    \int_{1}^{\infty} |\log(w - 1)|\frac{(w - 1)^{-(a - b + 1)}}{w^{n + 1}}\,dw.
  \end{multline*}
  The last integral can be bounded as
  \begin{align*}
    \int_{1}^{\infty} |\log(w - 1)|\frac{(w - 1)^{-(a - b + 1)}}{w^{n + 1}}\,dw
    &\leq \int_{1}^{2} (-\log(w - 1))(w - 1)^{-(a - b + 1)}\,dw
      + \int_{1}^{\infty} \frac{(w - 1)^{-(a - b + 1)}}{w^{n}}\,dw\\
    &= \frac{1}{\real(a - b)^{2}}
      + \frac{\Gamma(-\real(a - b))\Gamma(\real(a - b + n))}{(n - 1)!}.
  \end{align*}

  We now give a bound for \(R_{U,1}\). We have
  \begin{equation*}
    |R_{U,1}(a, b, z)|
    \leq C_{\chi}(a, b, n)\left|\log(z)^{-1}\frac{z^{n + a}}{\Gamma(a)}\right| \int_{0}^{\infty e^{-i\gamma}}\left|
      e^{-zt}\log(t)t^{n + a - 1}
    \right|\,dt.
  \end{equation*}
  If we replace \(t\) by \(\tau e^{-i\gamma}\) we get
  \begin{equation*}
    |R_{U,1}(a, b, z)|
    \leq C_{\chi}(a, b, n)\left|\log(z)^{-1}\frac{z^{n + a}}{\Gamma(a)}\right| \int_{0}^{\infty}
    e^{-\real(ze^{-i\gamma})\tau}|\log(\tau e^{-i\gamma})|\tau^{\real(n + a - 1)}\,d\tau.
  \end{equation*}
  By taking \(\gamma = \arg(z)\) we get \(\real(ze^{-i\gamma}) = |z|\), giving us
  \begin{equation*}
    |R_{U,1}(a, b, z)|
    \leq C_{\chi}(a, b, n)\left|\log(z)^{-1}\frac{z^{n + a}}{\Gamma(a)}\right| \int_{0}^{\infty}
    e^{-|z|\tau}|\log(\tau e^{-i\gamma})|\tau^{\real(n + a - 1)}\,d\tau.
  \end{equation*}
  The change of variables \(s = |z|\tau\) gives us
  \begin{equation*}
    \int_{0}^{\infty} e^{-|z|\tau}|\log(\tau e^{-i\gamma})|\tau^{\real(n + a - 1)}\,d\tau
    = |z|^{-\real(n +  a)}\int_{0}^{\infty} e^{-s}|\log(s|z|^{-1}e^{-i\gamma})|s^{\real(n + a - 1)}\,ds.
  \end{equation*}
  We have
  \begin{equation*}
    |\log(s|z|^{-1}e^{-i\gamma})| \leq |\log(|z|)| + |\log(s)| + |\gamma|,
  \end{equation*}
  giving us
  \begin{align*}
    \int_{0}^{\infty} e^{-s}|\log(s|z|^{-1}e^{-i\gamma})|s^{\real(n + a - 1)}\,ds
    &= |\log(|z|)|\int_{0}^{\infty} e^{-s}s^{\real(n + a - 1)}\,ds\\
    &\qquad+ \int_{0}^{\infty} e^{-s}|\log(s)|s^{\real(n + a - 1)}\,ds\\
    &\qquad+ |\gamma|\int_{0}^{\infty} e^{-s}s^{\real(n + a - 1)}\,ds.
  \end{align*}
  Let
  \begin{equation*}
    C_{1}(a, n) = \int_{0}^{\infty} e^{-s}s^{\real(n + a - 1)}\,ds,\quad
    C_{2}(a, n) = \int_{0}^{\infty} e^{-s}|\log(s)|s^{\real(n + a - 1)}\,ds.
  \end{equation*}
  We get \(C_{1}(a, n) = \Gamma(\real(a + n))\) and
  \begin{align*}
    C_{2}(a, n)
    &= \int_{0}^{1} e^{-s}(-\log(s))s^{\real(n + a - 1)}\,ds
      + \int_{1}^{\infty} e^{-s}\log(s)s^{\real(n + a - 1)}\,ds\\
    &\leq \int_{0}^{1} (-\log(s))s^{\real(n + a - 1)}\,ds
      + \int_{1}^{\infty} e^{-s}s^{\real(n + a)}\,ds\\
    &\leq \frac{1}{\real(a + n)^{2}} + \Gamma(\real(a + n + 1))
  \end{align*}

  Combining all of this gives us
  \begin{multline*}
    |R_{U,1}(a, b, z)|
    \leq C_{\chi}(a, b, n)\left|\log(z)^{-1}\frac{z^{n + a}}{\Gamma(a)}\right||\log(|z|)||z|^{-\real(n + a)}\\
    \left(
      C_{1}(a, n) + (|\gamma|C_{1}(a, n) + C_{2}(a, n))|\log(|z|)|^{-1}
    \right).
  \end{multline*}
  Using that \(\left|\frac{\log(|z|)}{\log(z)}\right| \leq 1\) and
  \begin{equation*}
    \frac{|z^{n + a}|}{|z|^{\real(n + a)}}
    = \frac{|z|^{\real(n + a)}e^{-\arg(z)\imag(n + a)}}{|z|^{\real(n + a)}}
    = e^{-\arg(z)\imag(a)},
  \end{equation*}
  we get
  \begin{equation*}
    |R_{U,1}(a, b, z)|
    \leq \frac{C_{\chi}(a, b, n)}{|\Gamma(a)|}e^{-\arg(z)\imag(a)}
    \left(
      C_{1}(a, n) + (|\gamma|C_{1}(a, n) + C_{2}(a, n))|\log(|z|)|^{-1}
    \right).
  \end{equation*}

  Next we give bounds for \(R_{U,2}\). We have
  \begin{equation*}
    |R_{U,2}(a, b, z)| \leq \left|\log(z)^{-1}\frac{z^{n + a}}{\Gamma(a)}\right|
    \int_{0}^{\infty e^{-i\gamma}}\left|e^{-zt}t^{a - 1}\frac{\partial}{\partial a}\chi_{n}(a, b, t)\right|\,dt.
  \end{equation*}
  Adding the bound for \(\frac{\partial}{\partial a}\chi_{n}(a, b, t)\)
  gives us
  \begin{equation*}
    |R_{U,2}(a, b, z)| \leq C_{\chi,a}(a, b, n)\left|\log(z)^{-1}\frac{z^{n + a}}{\Gamma(a)}\right|
    \int_{0}^{\infty e^{-i\gamma}}\left|e^{-zt}t^{n + a - 1}\right|\,dt.
  \end{equation*}
  As for \(R_{U,1}\), we replace \(t\) by \(\tau e^{-i\gamma}\)
  to get
  \begin{equation*}
    |R_{U,2}(a, b, z)| \leq C_{\chi,a}(a, b, n)\left|\log(z)^{-1}\frac{z^{n + a}}{\Gamma(a)}\right|
    \int_{0}^{\infty}e^{-\real(ze^{-i\gamma})\tau}\tau^{\real(n + a - 1)}\,d\tau.
  \end{equation*}
  Taking \(\gamma = \arg(z)\) and making the change of variables \(s
  = |z|\tau\) then gives
  \begin{equation*}
    |z|^{-\real(n + a)}\int_{0}^{\infty}e^{-s}s^{\real(n + a - 1)}\,ds.
  \end{equation*}
  This gives us
  \begin{equation*}
    |R_{U,2}(a, b, z)|
    \leq C_{\chi,a}(a, b, n)C_{1}(a, n)\left|\log(z)^{-1}\frac{z^{n + a}}{\Gamma(a)}\right||z|^{-\real(n + a)}
    \leq \frac{C_{\chi,a}(a, b, n)C_{1}(a, n)}{|\Gamma(a)|}e^{-\gamma\imag(a)}|\log(z)|^{-1}.
  \end{equation*}

  The final statement of the lemma follows in the same way as for
  Lemma~\ref{lemma:U}.
\end{proof}

These two lemmas will be the basis for the asymptotic bounds for \(P\)
and \(E\) and their derivatives. Before giving these bounds, let us
however introduce a bit more notation that will be helpful for the
handling of derivatives. To begin with, if we let
\begin{equation*}
  B_{W} = -\frac{1 + i\delta}{i\kappa}e^{-\sign(\imag(c))\pi i(b - a)}c^{b},
\end{equation*}
then \(J_{P}\) and \(J_{E}\) can be written as
\begin{align*}
  J_{P}(\xi) = B_{W}P(\xi)e^{-c\xi^{2}}\xi^{d - 1}
  \quad\text{and}\quad
  J_{E}(\xi) = B_{W}E(\xi)e^{-c\xi^{2}}\xi^{d - 1}.
\end{align*}
When enclosing \(I_{P,\kappa}\) and \(I_{P,\epsilon}\) we will have to
work with the functions \(J_{P,\kappa}\) and \(J_{P,\epsilon}\). For
this, it is convenient to write them in the form
\begin{equation*}
  J_{P,\kappa}(\xi) = D(\xi)e^{-c\xi^{2}}\xi^{d + 1}
  \quad\text{and}\quad
  J_{P,\epsilon}(\xi) = H(\xi)e^{-c\xi^{2}}\xi^{d + 1},
\end{equation*}
with \(D\) and \(H\) given by
\begin{align}
  \label{eq:D}
  D(\xi) &= P(\xi)(-c_{\kappa}B_{W} + B_{W,\kappa}\xi^{-2}) + P_{\kappa}(\xi)B_{W}\xi^{-2},\\
  \label{eq:H}
  H(\xi) &= P(\xi)(-c_{\epsilon}B_{W} + B_{W,\epsilon}\xi^{-2}) + P_{\epsilon}(\xi)B_{W}\xi^{-2}.
\end{align}

We are now ready to give the asymptotic bounds. These include bounds
for \(P\) and \(E\) as well as their derivatives with respect to
\(\xi\), \(\kappa\) and \(\epsilon\). For derivatives with respect to
\(\xi\) we go up to order three, which is related to how many
iterations of integration by parts are performed when enclosing
\(I_{P}\). We also give separate bounds for \(J_{P}\), \(J_{E}\),
\(D\), and \(H\) as well as some of their derivatives, since this
simplifies many of the later computations.

All of these bounds are relatively straightforward consequences of the
bounds for \(U\) and \(U_{a}\) from the two above lemmas. However,
assembling all of the explicit expressions is a rather tedious effort;
we therefore leave most of the details to
Appendix~\ref{sec:function-bounds-proofs}. It is worth pointing out
that the constants, of course, depend on the parameters \(\kappa\),
\(\omega\), \(\epsilon\) and \(\delta\) as well as the value of
\(\xi_{1}\), even if this dependence is not explicit in our notation.
Furthermore, the requirements of Lemmas~\ref{lemma:U} and
~\ref{lemma:U-a} give us conditions that need to be satisfied for
these bounds to be valid. We do not explicitly specify these conditions
here, instead, they are verified when these constants are actually
computed in the code associated with the computer-assisted proof.

\begin{lemma}
  \label{lemma:bounds-list}
  Let \(\xi_{1} > 1\). Under suitable conditions on the parameters
  \(\kappa\), \(\omega\), \(\epsilon\) and \(\delta\), we have, for
  \(\xi \geq \xi_{1}\), the bounds
  \begin{alignat*}{2}
    |P^{(m)}(\xi)| &\leq C_{P^{(m)}}\xi^{-\frac{1}{\sigma} - m}&\quad\text{for}\quad 0 \leq m \leq 3,\\
    |P_{\kappa}^{(m)}(\xi)| &\leq C_{P^{(m)},\kappa}\log(\xi)\xi^{-\frac{1}{\sigma} - m}&\quad\text{for}\quad 0 \leq m \leq 2,\\
    |P_{\epsilon}^{(m)}(\xi)| &\leq C_{P^{(m)},\epsilon}\xi^{-\frac{1}{\sigma} - m}&\quad\text{for}\quad 0 \leq m \leq 2,\displaybreak[3]\\
    |E^{(m)}(\xi)| &\leq C_{E^{(m)}}e^{\real(c)\xi^{2}}\xi^{\frac{1}{\sigma} - d + m}&\quad\text{for}\quad 0 \leq m \leq 3,\\
    |E_{\kappa}^{(m)}(\xi)| &\leq C_{E^{(m)},\kappa}e^{\real(c)\xi^{2}}\xi^{\frac{1}{\sigma} - d + 2 + m}&\quad\text{for}\quad 0 \leq m \leq 1,\\
    |E_{\epsilon}^{(m)}(\xi)| &\leq C_{E^{(m)},\epsilon}e^{\real(c)\xi^{2}}\xi^{\frac{1}{\sigma} - d + 2 + m}&\quad\text{for}\quad 0 \leq m \leq 1.
  \end{alignat*}
  Similarly, we have
  \begin{align*}
    |J_{P}(\xi)| &\leq C_{J_{P}}e^{-\real(c)\xi^{2}}\xi^{-\frac{1}{\sigma} + d - 1},\\
    |J_{P,\kappa}(\xi)| &\leq C_{J_{P},\kappa}e^{-\real(c)\xi^{2}}\xi^{-\frac{1}{\sigma} + d + 1},\\
    |J_{P,\epsilon}(\xi)| &\leq C_{J_{P},\epsilon}e^{-\real(c)\xi^{2}}\xi^{-\frac{1}{\sigma} + d + 1},\displaybreak[3]\\
    |J_{E}(\xi)| &\leq C_{J_{E}}\xi^{\frac{1}{\sigma} - 1},\\
    |J_{E,\kappa}(\xi)| &\leq C_{J_{E},\kappa}\log(\xi)\xi^{\frac{1}{\sigma} - 1},\\
    |J_{E,\epsilon}(\xi)| &\leq C_{J_{E},\epsilon}\xi^{\frac{1}{\sigma} - 1},
  \end{align*}
  as well as
  \begin{alignat*}{2}
    |D^{(m)}(\xi)| &\leq C_{D^{(m)}}\xi^{-\frac{1}{\sigma} - m}&\quad\text{for}\quad 0 \leq m \leq 3,\\
    |H^{(m)}(\xi)| &\leq C_{H^{(m)}}\xi^{-\frac{1}{\sigma} - m}&\quad\text{for}\quad 0 \leq m \leq 3.
  \end{alignat*}
\end{lemma}

\begin{proof}
  The bounds for \(P\) and \(E\) and their derivatives with respect to
  \(\xi\) are given in Lemma~\ref{lemma:C-P-E}, the derivatives with
  respect to \(\kappa\) in Lemma~\ref{lemma:C-P-E-dkappa} and the
  derivatives with respect to \(\epsilon\) in
  Lemma~\ref{lemma:C-P-E-depsilon}. The bounds for \(J_{P}\) and
  \(J_{E}\) are given in Lemma~\ref{lemma:C-J-P-J-E} and their
  derivatives with respect to \(\kappa\) and \(\epsilon\) in
  Lemma~\ref{lemma:C-J-P-kappa-J-E-kappa}. Finally \(D\) and \(H\) are
  handled in Lemma~\ref{lemma:C-D-H}.

  Note that in all cases, the constants have explicit expressions
  given in terms of the bounds for \(U\) and \(U_{a}\) from
  Lemmas~\ref{lemma:U} and ~\ref{lemma:U-a}. For the derivatives with
  respect to \(\xi\), similar bounds would of course hold for higher
  order derivatives. However, the orders included in the lemma are the
  only ones for which we give explicit expressions for the bounds.
\end{proof}

\subsection{Existence of fixed point}
\label{sec:fixed-point}
In this subsection we give explicit conditions for the existence of a
fixed point of the operator
\begin{equation*}
  T(Q)(\xi) = \gamma P(\xi) + P(\xi)I_{E}(\xi) + E(\xi)I_{P}(\xi)
\end{equation*}
from~\eqref{eq:T}. The fixed-point argument is performed in the Banach
space \(\Space_\normv\), with \(\normv \geq 0\), of continuous
functions \(Q : [\xi_1, \infty) \to \mathbb{C}\) for which the norm
\begin{equation*}
  \|Q\|_{\normv} = \sup_{\xi \geq \xi_{1}} \xi^{\frac{1}{\sigma} - \normv}|Q(\xi)|
\end{equation*}
is finite. While the inclusion of \(\normv\) in the norm is not
strictly necessary for the fixed-point argument itself, it is required
to bound the derivatives with respect to \(\kappa\) and \(\epsilon\),
where \(\log\) terms appear in the leading order.

To begin, we provide the following lemma that gives us bounds for
\(I_{P}\) and \(I_{E}\) in terms of the norm \(\|\cdot\|_{\normv}\).

\begin{lemma}\label{lemma:I_P-I_E}
  Assume that \(\xi_{1} > 1\),
  \((2\sigma + 1)\normv - \frac{2}{\sigma} + d - 2 < 0\) and
  \((2\sigma + 1)\normv - 2 < 0\). Then, for \(\xi \geq \xi_{1}\), we
  have the following bounds
  \begin{align*}
    |I_{P}(\xi)| &\leq C_{I_{P}}\|Q\|_{\normv}^{2\sigma + 1}e^{-\real(c)\xi^{2}}\xi^{(2\sigma + 1)\normv - \frac{2}{\sigma} + d - 2},\\
    |I_{E}(\xi)| &\leq C_{I_{E}}\|Q\|_{\normv}^{2\sigma + 1}\xi_{1}^{(2\sigma + 1)\normv - 2},
  \end{align*}
  with
  \begin{align*}
    C_{I_{P}} &= \min\left(\frac{C_{J_{P}}}{|(2\sigma + 1)\normv - \frac{2}{\sigma} + d - 2|}, \frac{C_{J_{P}}}{2\real(c)}\xi_{1}^{-2}\right),\\
    C_{I_{E}} &= \frac{C_{J_{E}}}{|(2\sigma + 1)\normv - 2|}.
  \end{align*}
  The second term in the minimum for \(C_{I_{P}}\), is understood to
  be finite only for \(\real(c) > 0\).
\end{lemma}

\begin{proof}
  Let us start by noting that
  \begin{equation*}
    |Q(\eta)|^{2\sigma + 1}
    \leq \|Q\|_{\normv}^{2\sigma + 1}\eta^{(2\sigma + 1)\left(\normv - \frac{1}{\sigma}\right)}
    = \|Q\|_{\normv}^{2\sigma + 1}\eta^{(2\sigma + 1)\normv - \frac{1}{\sigma} - 2}.
  \end{equation*}
  Combining this with the bound
  \(|J_{P}(\eta)| \leq
  C_{J_{P}}e^{-\real(c)\eta^{2}}\eta^{-\frac{1}{\sigma} + d - 1}\)
  from Lemma~\ref{lemma:bounds-list} gives us
  \begin{equation*}
    |I_{P}(\xi)|
    \leq \int_{\xi}^{\infty} |J_{P}(\eta)||Q(\eta)|^{2\sigma + 1}\,d\eta
    \leq C_{J_{P}}\|Q\|_{\normv}^{2\sigma + 1}
    \int_{\xi}^{\infty} e^{-\real(c)\eta^{2}}\eta^{(2\sigma + 1)\normv - \frac{2}{\sigma} + d - 3}\,d\eta.
  \end{equation*}
  The bound \(|J_{E}(\xi)| \leq C_{J_{E}}\xi^{\frac{1}{\sigma} - 1}\)
  similarly gives us
  \begin{equation*}
    |I_{E}(\xi)|
    \leq \int_{\xi_{1}}^{\xi}|J_{E}(\eta)||Q(\eta)|^{2\sigma + 1}\,d\eta
    \leq C_{J_{E}}\|Q\|_{\normv}^{2\sigma + 1}
    \int_{\xi_{1}}^{\xi} \eta^{(2\sigma + 1)\normv - 3}\,d\eta.
  \end{equation*}

  For \(\real(c) \geq 0\), the first integral can be bounded as
  \begin{equation*}
    \begin{split}
      \int_{\xi}^{\infty} e^{-\real(c)\eta^{2}}\eta^{(2\sigma + 1)\normv - \frac{2}{\sigma} + d - 3}\,d\eta
      &\leq e^{-\real(c)\xi^{2}} \int_{\xi}^{\infty} \eta^{(2\sigma + 1)\normv - \frac{2}{\sigma} + d - 3}\,d\eta\\
      &= \frac{1}{|(2\sigma + 1)\normv - \frac{2}{\sigma} + d - 2|}e^{-\real(c)\xi^{2}}\xi^{(2\sigma + 1)\normv - \frac{2}{\sigma} + d - 2}.
    \end{split}
  \end{equation*}
  In case \(\real(c) > 0\) we can get a better estimate by first
  integrating by parts. In this case we have
  \begin{equation*}
    \begin{split}
      \int_{\xi}^{\infty} e^{-\real(c)\eta^{2}}\eta^{(2\sigma + 1)\normv - \frac{2}{\sigma} + d - 3}\,d\eta
      &= \frac{1}{2\real(c)}e^{-\real(c)\xi^{2}}\xi^{(2\sigma + 1)\normv - \frac{2}{\sigma} + d - 4}\\
      &\quad+ \frac{1}{2\real(c)}\int_{\xi}^{\infty} e^{-\real(c)\eta^{2}}\frac{d}{d\eta}\left(\eta^{(2\sigma + 1)\normv - \frac{2}{\sigma} + d - 3}\right)\,d\eta.
    \end{split}
  \end{equation*}
  Since
  \(\frac{d}{d\eta}\left(\eta^{(2\sigma + 1)\normv - \frac{2}{\sigma}
      + d - 3}\right) < 0\), the last integral is negative and we get
  an upper bound by removing it. This gives us the upper bound
  \begin{equation*}
    \begin{split}
      \int_{\xi}^{\infty} e^{-\real(c)\eta^{2}} \eta^{(2\sigma + 1)\normv - \frac{2}{\sigma} + d - 3}\,d\eta
      &\leq \frac{1}{2\real(c)}e^{-\real(c)\xi^{2}} \xi^{(2\sigma + 1)\normv - \frac{2}{\sigma} + d - 4}\\
      &\leq \frac{1}{2\real(c)}\xi_{1}^{-2}e^{-\real(c)\xi^{2}} \xi^{(2\sigma + 1)\normv - \frac{2}{\sigma} + d - 2}
    \end{split}
  \end{equation*}
  To avoid a bound which degenerates as \(\real(c)\) approaches zero,
  we take the minimum of this bound and the above one.

  For the second integral we immediately get the bound
  \begin{equation*}
    \int_{\xi_{1}}^{\xi} \eta^{(2\sigma + 1)\normv - 3}\,d\eta
    \leq \frac{1}{|(2\sigma + 1)\normv - 2|}\xi_{1}^{(2\sigma + 1)\normv - 2}.
  \end{equation*}
\end{proof}

Consequently we obtain the following version of~\cite[Lemma
A.2]{Plech2001} with explicit constants.

\begin{lemma}\label{lemma:fixed-point-bounds}
  Under the assumptions of Lemma~\ref{lemma:I_P-I_E}, the operator
  \(T\) defines a continuous mapping
  \(T : \Space_\normv \to \Space_\normv\). Moreover
  \begin{equation}
    \label{eq:T-1}
    \|T(Q)\|_{\normv} \leq C_{P}|\gamma|\xi_{1}^{-\normv} + C_{T}\xi_{1}^{2\sigma \normv - 2}\|Q\|_{\normv}^{2\sigma + 1}
  \end{equation}
  and
  \begin{equation}
    \label{eq:T-2}
    \|T(u) - T(v)\|_{\normv} \leq M_{\sigma}C_{T}\xi_{1}^{-2 + 2\sigma \normv}\|u - v\|_{\normv}(\|u\|_{\normv}^{2\sigma} + \|v\|_{\normv}^{2\sigma}),
  \end{equation}
  for all \(u, v \in \Space_\normv\). Here
  \begin{equation*}
    C_{T} = C_{P}C_{I_{E}} + C_{E}C_{I_{P}}
  \end{equation*}
  and \(M_{\sigma}\) is as in Lemma~\ref{lemma:M}.
\end{lemma}

\begin{proof}
  To prove~\eqref{eq:T-1}, we need to bound
  \(\xi^{\frac{1}{\sigma} - \normv}|T(Q)(\xi)|\) for
  \(\xi \geq \xi_{1}\). Bounding it termwise we get
  \begin{equation*}
    \xi^{\frac{1}{\sigma} - \normv}|T(Q)(\xi)| \leq
    \xi^{\frac{1}{\sigma} - \normv}|\gamma||P(\xi)|
    + \xi^{\frac{1}{\sigma} - \normv}|P(\xi)||I_{E}(\xi)|
    + \xi^{\frac{1}{\sigma} - \normv}|E(\xi)||I_{P}(\xi)|.
  \end{equation*}
  From Lemma~\ref{lemma:bounds-list} we get for the first term
  \begin{equation*}
    \xi^{\frac{1}{\sigma} - \normv}|\gamma||P(\xi)| \leq C_{P}|\gamma|\xi^{-\normv} \leq C_{P}|\gamma|\xi_{1}^{-\normv}.
  \end{equation*}
  Combining Lemmas~\ref{lemma:bounds-list} and~\ref{lemma:I_P-I_E} we
  have for the second and third term
  \begin{equation*}
    \xi^{\frac{1}{\sigma} - \normv}|P(\xi)||I_{E}(\xi)|
    \leq C_{P}C_{I_{E}}\|Q\|_{\normv}^{2\sigma + 1}\xi^{-\normv}\xi_{1}^{(2\sigma + 1)\normv - 2}
    \leq C_{P}C_{I_{E}}\|Q\|_{\normv}^{2\sigma + 1}\xi_{1}^{2\sigma\normv - 2}
  \end{equation*}
  and
  \begin{equation*}
    \xi^{\frac{1}{\sigma} - \normv}|E(\xi)||I_{P}(\xi)|
    \leq C_{E}C_{I_{P}}\|Q\|_{\normv}^{2\sigma + 1}\xi^{2\sigma\normv - 2}
    \leq C_{E}C_{I_{P}}\|Q\|_{\normv}^{2\sigma + 1}\xi_{1}^{2\sigma\normv - 2}
  \end{equation*}
  respectively. Combining these three bounds gives us the result.

  For \eqref{eq:T-2} we get
  \begin{equation*}
    \begin{split}
      \xi^{\frac{1}{\sigma} - \normv}|T(u)(\xi) - T(v)(\xi)|
      &\leq \xi^{\frac{1}{\sigma} - \normv}|P(\xi)|
        \int_{\xi_{1}}^{\xi}|J_{E}(\eta)|\left||u(\eta)|^{2\sigma}u(\eta) - |v(\eta)|^{2\sigma}v(\eta)\right|\,d\eta\\
      &\quad+ \xi^{\frac{1}{\sigma} - \normv}|E(\xi)|
        \int_{\xi}^{\infty}|J_{P}(\eta)|\left||u(\eta)|^{2\sigma}u(\eta) - |v(\eta)|^{2\sigma}v(\eta)\right|\,d\eta.
    \end{split}
  \end{equation*}
  By Lemma~\ref{lemma:M} we have the bound
  \begin{equation*}
    \left||u(\eta)|^{2\sigma}u(\eta) - |v(\eta)|^{2\sigma}v(\eta)\right|
    \leq M_{\sigma} |u(\eta) - v(\eta)|(|u(\eta)|^{2\sigma} + |v(\eta)|^{2\sigma}),
  \end{equation*}
  giving us
  \begin{equation*}
    \left||u(\eta)|^{2\sigma}u(\eta) - |v(\eta)|^{2\sigma}v(\eta)\right|
    \leq M_{\sigma} \|u - v\|_{\normv}(\|u\|_{\normv}^{2\sigma} + \|v\|_{\normv}^{2\sigma})\eta^{(2\sigma + 1)\normv - \frac{1}{\sigma} - 2}.
  \end{equation*}
  From which we get
  \begin{multline*}
    \xi^{\frac{1}{\sigma} - \normv}|P(\xi)|
    \int_{\xi_{1}}^{\xi}|J_{E}(\eta)|\left||u(\eta)|^{2\sigma}u(\eta) - |v(\eta)|^{2\sigma}v(\eta)\right|\,d\eta\\
    \leq M_{\sigma} \|u - v\|_{\normv}(\|u\|_{\normv}^{2\sigma} + \|v\|_{\normv}^{2\sigma})
    \xi^{\frac{1}{\sigma} - \normv}|P(\xi)|
    \int_{\xi_{1}}^{\xi}|J_{E}(\eta)|\eta^{(2\sigma + 1)\normv - \frac{1}{\sigma} - 2}\,d\eta
  \end{multline*}
  and
  \begin{multline*}
    \xi^{\frac{1}{\sigma} - \normv}|E(\xi)|
    \int_{\xi}^{\infty}|J_{P}(\eta)|\left||u(\eta)|^{2\sigma}u(\eta) - |v(\eta)|^{2\sigma}v(\eta)\right|\,d\eta\\
    \leq M_{\sigma} \|u - v\|_{\normv}(\|u\|_{\normv}^{2\sigma} + \|v\|_{\normv}^{2\sigma})
    \xi^{\frac{1}{\sigma} - \normv}|E(\xi)|
    \int_{\xi}^{\infty}|J_{P}(\eta)|\eta^{(2\sigma + 1)\normv - \frac{1}{\sigma} - 2}\,d\eta
  \end{multline*}
  These integrals can be bounded in the same way as in
  Lemma~\ref{lemma:I_P-I_E}, giving us the result.

  The continuity follows in the same way as in~\cite[Proposition
  A.2]{Plech2001}.
\end{proof}

The estimates \eqref{eq:T-1} and \eqref{eq:T-2} show that \(T\) is a
contraction of the ball
\(B_\rho = \{Q \in \Space_\normv : \|Q\|_\normv \leq \rho\}\) into
itself if \(\rho\) is such that
\begin{align}
  \label{eq:T-ineq-1}
  C_{P}|\gamma|\xi_{1}^{-\normv} + C_{T}\xi_{1}^{2\sigma \normv - 2}\rho^{2\sigma + 1} &\leq \rho,\\
  \label{eq:T-ineq-2}
  2M_{\sigma}C_{T}\xi_{1}^{2\sigma \normv - 2}\rho^{2\sigma} &< 1.
\end{align}
From the Banach fixed point theorem we then get the existence of a
fixed point. More precisely we have the following result,
see~\cite[Proposition A.3]{Plech2001}.

\begin{proposition}\label{prop:fixed-point}
  Under the assumptions of Lemma~\ref{lemma:fixed-point-bounds} and if
  the two inequalities \eqref{eq:T-ineq-1} and \eqref{eq:T-ineq-2} are
  satisfied, then the operator \(T\) has a unique fixed point \(Q\) in
  the ball \(B_{\rho}\).
\end{proposition}

\subsection{Enclosures for fixed point}
\label{sec:fixed-point-enclosures}
The above proposition allows us to prove the existence of a solution
\(Q\) satisfying \(T(Q) = Q\) and gives us bounds on its norm. It
does, however, not directly give us an enclosure for \(Q(\xi_{1})\),
other than what can be deduced from the bound of the norm. The goal of
this section is to develop bounds that allow us to compute accurate
enclosures of \(Q\), as well as its derivatives, at \(\xi_{1}\).

The fixed-point equation satisfied by \(Q\), combined with the bounds
for \(I_{P}\) and \(I_{E}\) from Lemma~\ref{lemma:I_P-I_E}, already
allows for the computation of an improved enclosure for \(Q(\xi_{1})\).
From the fixed-point equation we have
\begin{equation*}
  Q(\xi_{1}) = \gamma P(\xi_{1}) + P(\xi_{1})I_{E}(\xi_{1}) + E(\xi_{1})I_{P}(\xi_{1}).
\end{equation*}
In the specific case when \(\xi = \xi_{1}\), we have
\(I_{E}(\xi_{1}) = 0\), hence giving us
\begin{equation*}
  Q(\xi_{1}) = \gamma P(\xi_{1}) + E(\xi_{1})I_{P}(\xi_{1}).
\end{equation*}
One can directly compute tight enclosures for \(P(\xi_{1})\) and
\(E(\xi_{1})\). What remains is enclosing \(I_{P}(\xi_{1})\), which
can be done using Lemma~\ref{lemma:I_P-I_E} and the bounds for the
norm of \(Q\) from Proposition~\ref{prop:fixed-point}. While this does
give a significantly improved enclosure of \(Q(\xi_{1})\) compared to
what can be deduced from the bound of the norm alone, it is still not
good enough to prove the existence of a zero for the function
\(G\)~\eqref{eq:G}. To get better enclosures we will refine the
estimates for \(I_{P}\) through integration by parts.

To compute enclosures of the various derivatives of \(Q\) we directly
differentiate the fixed-point equation, giving us equations for the
derivatives. These equations will consist of derivatives of \(P\) and
\(E\), as well as \(I_{P}\) and \(I_{E}\). Enclosures are then
computed similarly to those for \(Q\), by computing initial bounds for the
norms and combining this with refined enclosures for the derivatives
of \(I_{P}\) and \(I_{E}\).

We split the work of this part into four subsections:
\begin{enumerate}
\item In Section~\ref{sec:I_P-I_E-direct-bounds} we give various
  bounds for \(I_{E}\) and \(I_{P}\) and their derivatives. The
  approach is similar to that of Lemma~\ref{lemma:I_P-I_E}, and follows
  from direct bounds of the involved functions.
\item In Section~\ref{sec:I_P-integration-by-parts-bounds} we give
  refined bounds for \(I_{P}\) and its derivatives based on integration
  by parts.
\item In Section~\ref{sec:bounds-on-norms} we derive bounds for the
  norms of the various derivatives of \(Q\). These are based on the
  bound for the norm of \(Q\) from the fixed-point argument and the
  bounds for derivatives of \(I_{P}\) and \(I_{E}\) from the earlier
  subsection.
\item In Section~\ref{sec:enclosures-xi-1} we describe how to compute
  tight enclosures at \(\xi = \xi_{1}\) using the derived bounds.
\end{enumerate}

Before continuing, let us start by specifying the various quantities
we will be bounding. Recall that for evaluating the function \(G\) and
its Jacobian we need to compute enclosures for \(Q\) and \(Q'\), as
well as their derivatives with respect to \(\real\gamma\),
\(\imag\gamma\), \(\kappa\) and \(\epsilon\). To apply the refined
integration by parts bounds for \(I_{P}\) we also need bounds for
higher order derivatives in \(\xi\), namely \(Q''\) and \(Q'''\). In
the end, we are thus interested in the functions
\begin{equation*}
  Q,\ Q',\ Q'',\ Q''',\ Q_{\real\gamma},\ Q_{\real\gamma}',\ Q_{\imag\gamma},\ Q_{\imag\gamma}',\ Q_{\kappa},\ Q_{\kappa}',\ Q_{\epsilon} \text{ and } Q_{\epsilon}'.
\end{equation*}
From the fixed-point equation for \(Q\) we can derive equations
satisfied by the derivatives. For derivatives with respect to \(\xi\)
we note the following identities, which follow from
\(I_{E}'(\xi) = J_{E}(\xi)|Q(\xi)|^{2\sigma}Q(\xi)\) and
\(I_{P}'(\xi) = -J_{P}(\xi)|Q(\xi)|^{2\sigma}Q(\xi)\), together with
the Wronskian \(W = PE' - P'E\):
\begin{align*}
  P(\xi)I_{E}'(\xi) + E(\xi)I_{P}'(\xi)
  &= 0,\\
  P'(\xi)I_{E}'(\xi) + E'(\xi)I_{P}'(\xi)
  &= -\frac{1 + i\delta}{1 - i\epsilon}|Q(\xi)|^{2\sigma}Q(\xi),\\
  P''(\xi)I_{E}'(\xi) + E''(\xi)I_{P}'(\xi)
  &= -\frac{1 + i\delta}{1 - i\epsilon}\frac{W'(\xi)}{W(\xi)}|Q(\xi)|^{2\sigma}Q(\xi)\\
  &= \frac{1 + i\delta}{1 - i\epsilon}((d - 1)\xi^{-1} - 2c\xi)|Q(\xi)|^{2\sigma}Q(\xi).
\end{align*}
This gives us
\begin{align*}
  Q(\xi)
  &= \gamma P(\xi)
    + P(\xi)I_{E}(\xi) + E(\xi)I_{P}(\xi),\\
  Q'(\xi)
  &= \gamma P'(\xi)
    + P'(\xi)I_{E}(\xi) + E'(\xi)I_{P}(\xi),\\
  Q''(\xi)
  &= \gamma P''(\xi)
    + P''(\xi)I_{E}(\xi)
    + E''(\xi)I_{P}(\xi)
    - \frac{1 + i\delta}{1 - i\epsilon}|Q(\xi)|^{2\sigma}Q(\xi),\\
  Q'''(\xi)
  &= \gamma P'''(\xi)
    + P'''(\xi)I_{E}(\xi)
    + E'''(\xi)I_{P}(\xi)
    + \frac{1 + i\delta}{1 - i\epsilon}((d - 1)\xi^{-1} - 2c\xi)|Q(\xi)|^{2\sigma}Q(\xi)\\
    &\quad- \frac{1 + i\delta}{1 - i\epsilon}\frac{\partial}{\partial \xi}\left(|Q(\xi)|^{2\sigma}Q(\xi)\right).
\end{align*}
For the derivatives with respect to the parameters we get
\begin{align*}
  Q_{\real\gamma}(\xi)
  &= P(\xi)
    + P(\xi)I_{E,\real\gamma}(\xi)
    + E(\xi)I_{P,\real\gamma}(\xi),\\
  Q_{\real\gamma}'(\xi)
  &= P'(\xi)
    + P'(\xi)I_{E,\real\gamma}(\xi)
    + E'(\xi)I_{P,\real\gamma}(\xi),\\
  Q_{\imag\gamma}(\xi)
  &= iP(\xi)
    + P(\xi)I_{E,\imag\gamma}(\xi)
    + E(\xi)I_{P,\imag\gamma}(\xi),\\
  Q_{\imag\gamma}'(\xi)
  &= iP'(\xi)
    + P'(\xi)I_{E,\imag\gamma}(\xi)
    + E'(\xi)I_{P,\imag\gamma}(\xi),\\
  Q_{\kappa}(\xi)
  &= \gamma P_{\kappa}(\xi)
    + P_{\kappa}(\xi)I_{E}(\xi)
    + P(\xi)I_{E,\kappa}(\xi)
    + E_{\kappa}(\xi)I_{P}(\xi)
    + E(\xi)I_{P,\kappa}(\xi),\\
  Q_{\kappa}'(\xi)
  &= \gamma P_{\kappa}'(\xi)
    + P_{\kappa}'(\xi)I_{E}(\xi)
    + P'(\xi)I_{E,\kappa}(\xi)
    + E_{\kappa}'(\xi)I_{P}(\xi)
    + E'(\xi)I_{P,\kappa}(\xi),\\
  Q_{\epsilon}(\xi)
  &= \gamma P_{\epsilon}(\xi)
    + P_{\epsilon}(\xi)I_{E}(\xi)
    + P(\xi)I_{E,\epsilon}(\xi)
    + E_{\epsilon}(\xi)I_{P}(\xi)
    + E(\xi)I_{P,\epsilon}(\xi),\\
  Q_{\epsilon}'(\xi)
  &= \gamma P_{\epsilon}'(\xi)
    + P_{\epsilon}'(\xi)I_{E}(\xi)
    + P'(\xi)I_{E,\epsilon}(\xi)
    + E_{\epsilon}'(\xi)I_{P}(\xi)
    + E'(\xi)I_{P,\epsilon}(\xi).
\end{align*}
Note that all of the various derivatives of \(P\) and \(E\) appearing
in these expressions are covered by Lemma~\ref{lemma:bounds-list}. The
goal of Sections~\ref{sec:I_P-I_E-direct-bounds}
and~\ref{sec:I_P-integration-by-parts-bounds} is to derive bounds for
all of the involved derivatives of \(I_{P}\) and \(I_{E}\).

\subsubsection{Direct bounds for \(I_{P}\) and \(I_{E}\)}
\label{sec:I_P-I_E-direct-bounds}
In this section we give bounds for \(I_{P}\) and \(I_{E}\), as well as
their derivatives with respect to the parameters, based on direct
bounds for the involved functions and integrands (as opposed to
through integration by parts).

The bounds we derive will depend on a priori bounds for the norms of
\(Q\) and its derivatives. In particular, we will need bounds for
derivatives of \(|Q|^{2\sigma}Q\) in terms of bounds for the norms of
the derivatives of \(Q\). For this reason, the following lemma is
useful.
\begin{lemma}\label{lemma:norm-bounds}
  For \(\xi \geq \xi_{1}\) we have the bounds
  \begin{align*}
    \left||Q(\xi)|^{2\sigma}Q(\xi)\right|
    &\leq \|Q\|_{\normv}^{2\sigma + 1}\xi^{(2\sigma + 1)\normv - \frac{1}{\sigma} - 2},\\
    \left|\frac{\partial}{\partial \xi}\left(|Q(\xi)|^{2\sigma}Q(\xi)\right)\right|
    &\leq (2\sigma + 1)\|Q\|_{\normv}^{2\sigma}\|Q'\|_{\normv}\xi^{(2\sigma + 1)\normv - \frac{1}{\sigma} - 2},\\
    \left|\frac{\partial^{2}}{\partial \xi^{2}}\left(|Q(\xi)|^{2\sigma}Q(\xi)\right)\right|
    &\leq (2\sigma + 1)(
      2\sigma\|Q'\|_{\normv}^{2}
      + \|Q\|_{\normv}\|Q''\|_{\normv}
      )\|Q\|_{\normv}^{2\sigma - 1}\xi^{(2\sigma + 1)\normv - \frac{1}{\sigma} - 2},\\
    \left|\frac{\partial^{3}}{\partial \xi^{3}}\left(|Q(\xi)|^{2\sigma}Q(\xi)\right)\right|
    &\leq (2\sigma + 1)(
      2\sigma(2\sigma - 1)\|Q'\|_{\normv}^{3}
      + 6\sigma\|Q\|_{\normv}\|Q'\|_{\normv}\|Q''\|_{\normv}\\
    &\qquad+ \|Q\|_{\normv}^{2}\|Q'''\|_{\normv}
      )\|Q\|_{\normv}^{2\sigma - 2}\xi^{(2\sigma + 1)\normv - \frac{1}{\sigma} - 2},
  \end{align*}
  as well as
  \begin{align*}
    \left|\frac{\partial}{\partial \real\gamma}\left(|Q(\xi)|^{2\sigma}Q(\xi)\right)\right|
    &\leq (2\sigma + 1)\|Q\|_{\normv}^{2\sigma}\|Q_{\real\gamma}\|_{\normv}\xi^{(2\sigma + 1)\normv - \frac{1}{\sigma} - 2},\\
    \left|\frac{\partial}{\partial \real\gamma}\frac{\partial}{\partial \xi}\left(|Q(\xi)|^{2\sigma}Q(\xi)\right)\right|
    &\leq (2\sigma + 1)(
      2\sigma\|Q'\|_{\normv}\|Q_{\real\gamma}\|_{\normv}
      + \|Q\|_{\normv}\|Q_{\real\gamma}'\|_{\normv}
      )\|Q\|_{\normv}^{2\sigma - 1}\xi^{(2\sigma + 1)\normv - \frac{1}{\sigma} - 2},
  \end{align*}
  with identical bounds for the derivatives with respect to
  \(\imag\gamma\), \(\kappa\) and \(\epsilon\).
\end{lemma}

\begin{proof}
  The bounds for the first order derivatives follow from
  \begin{equation*}
    \frac{\partial}{\partial t}\left(|Q|^{2\sigma}Q\right)
    = (\sigma + 1)|Q|^{2\sigma}\frac{\partial Q}{\partial t}
    + \sigma|Q|^{2\sigma-2}Q^{2}\overline{\left(\frac{\partial Q}{\partial t}\right)},
  \end{equation*}
  valid for any real direction \(t\), in particular for
  \(t \in \{\xi,\ \real\gamma,\ \imag\gamma,\ \kappa,\ \epsilon\}\).
  Bounds for remaining derivatives follow by differentiating this
  identity further.
\end{proof}

For \(I_{E}\) and its derivatives with respect to \(\gamma\),
\(\kappa\) and \(\epsilon\) there are no oscillations to take into
account, and direct bounds suffice to get good enough estimates. For
\(I_{E}\), this means that the bounds in Lemma~\ref{lemma:I_P-I_E}
suffice. For the derivatives with respect to the parameters we have
the following lemma.
\begin{lemma}\label{lemma:I-E-bounds}
  If \((2\sigma + 1)\normv - 2 < 0\), then for
  \(\xi \geq \xi_{1} \geq 1\) we have
  \begin{align*}
    |I_{E,\real\gamma}(\xi)|
    &\leq (2\sigma + 1)C_{I_{E}}\|Q\|_{\normv}^{2\sigma}
      \|Q_{\real\gamma}\|_\normv\xi_{1}^{(2\sigma + 1)\normv - 2},\\
    |I_{E,\imag\gamma}(\xi)|
    &\leq (2\sigma + 1)C_{I_{E}}\|Q\|_{\normv}^{2\sigma}
      \|Q_{\imag\gamma}\|_\normv\xi_{1}^{(2\sigma + 1)\normv - 2},\\
    \left|I_{E,\kappa}(\xi)\right|
    &\leq (
      C_{I_{E},\kappa}\|Q\|_{\normv}
      + (2\sigma + 1)C_{I_{E}}\left\|Q_{\kappa}\right\|_{\normv}
      )\|Q\|_{\normv}^{2\sigma}\xi_{1}^{(2\sigma + 1)\normv - 2},\\
    \left|I_{E,\epsilon}(\xi)\right|
    &\leq (
      C_{I_{E},\epsilon}\|Q\|_{\normv}
      + (2\sigma + 1)C_{I_{E}}\left\|Q_{\epsilon}\right\|_{\normv}
      )\|Q\|_{\normv}^{2\sigma}\xi_{1}^{(2\sigma + 1)\normv - 2}.
  \end{align*}
  Here \(C_{I_{E}}\) is as in Lemma~\ref{lemma:I_P-I_E} and
  \begin{align*}
    C_{I_{E},\kappa}
    &= C_{J_{E},\kappa}\left(
      \frac{1}{((2\sigma + 1)\normv - 2)^2}
      + \frac{\log(\xi_{1})}{|(2\sigma + 1)\normv - 2|}
      \right),\\
    C_{I_{E},\epsilon}
    &= \frac{C_{J_{E},\epsilon}}{|(2\sigma + 1)\normv - 2|}.
  \end{align*}
\end{lemma}

\begin{proof}
  Let us start by bounding \(I_{E,\kappa}\), for which we have
  \begin{equation*}
    I_{E,\kappa}(\xi)
    = \int_{\xi_{1}}^{\xi}J_{E,\kappa}(\eta)|Q(\eta)|^{2\sigma}Q(\eta)\,d\eta
    + \int_{\xi_{1}}^{\xi}J_{E}(\eta)
    \frac{\partial}{\partial \kappa}(|Q(\eta)|^{2\sigma}Q(\eta))\,d\eta.
  \end{equation*}
  From Lemma~\ref{lemma:bounds-list} we have
  \(|J_{E}(\xi)| \leq C_{J_{E}}\xi^{\frac{1}{\sigma} - 1}\) and
  \(|J_{E,\kappa}(\xi)| \leq
  C_{J_{E},\kappa}\log(\xi)\xi^{\frac{1}{\sigma} - 1}\), combined with
  the bounds for \(Q\) from Lemma~\ref{lemma:norm-bounds} this gives
  us
  \begin{equation*}
    |I_{E,\kappa}(\xi)|
    \leq
    C_{J_{E},\kappa}\|Q\|_{\normv}^{2\sigma + 1}\int_{\xi_{1}}^{\xi}\log(\eta)\eta^{(2\sigma + 1)\normv - 3}\,d\eta
    + (2\sigma + 1)C_{J_{E}}\|Q\|_{\normv}^{2\sigma}\|Q_{\kappa}\|_{\normv}\int_{\xi_{1}}^{\xi}\eta^{\frac{1}{\sigma} - 3}\,d\eta.
  \end{equation*}
  The second integral is bounded by
  \(\frac{\xi_{1}^{(2\sigma + 1)\normv -2}}{|(2\sigma + 1)\normv -
    2|}\). The first integral we bound as
  \begin{equation*}
    \int_{\xi_{1}}^{\xi} \log(\eta)\eta^{(2\sigma + 1)\normv - 3}\,d\eta
    \leq \int_{\xi_{1}}^{\infty} \log(\eta)\eta^{(2\sigma + 1)\normv - 3}\,d\eta
    = \left(
      \frac{1}{((2\sigma + 1)\normv - 2)^{2}}
      + \frac{\log(\xi_{1})}{|(2\sigma + 1)\normv - 2|}
    \right)\xi_{1}^{(2\sigma + 1)\normv - 2}.
  \end{equation*}
  Combining these two bounds gives us the result for \(I_{E,\kappa}\).

  The bounds for \(I_{E,\real\gamma}\), \(I_{E,\imag\gamma}\) and
  \(I_{E,\epsilon}\) are similar. In the case of \(I_{E,\real\gamma}\)
  and \(I_{E,\imag\gamma}\) there is only one term to bound, since
  \(J_{E}\) does not depend on \(\gamma\). For \(I_{E,\epsilon}\), the
  bound
  \(|J_{E,\epsilon}(\xi)| \leq C_{J_{E},\epsilon}\xi^{\frac{1}{\sigma}
    - 1}\) means that there is no logarithmic factor, but otherwise
  the approach is the same.
\end{proof}

Finally, we get to the bounds for \(I_{P}\) and its derivatives with
respect to the parameters. For these derivatives we get
\begin{align*}
  I_{P,\real\gamma}(\xi)
  &= \int_{\xi}^{\infty}J_{P}(\eta)
    \frac{\partial}{\partial \real\gamma}(|Q(\eta)|^{2\sigma}Q(\eta))\,d\eta,\\
  I_{P,\imag\gamma}(\xi)
  &= \int_{\xi}^{\infty}J_{P}(\eta)
    \frac{\partial}{\partial \imag\gamma}(|Q(\eta)|^{2\sigma}Q(\eta))\,d\eta,\\
  I_{P,\kappa}(\xi)
  &= \int_{\xi}^{\infty}J_{P,\kappa}(\eta)|Q(\eta)|^{2\sigma}Q(\eta)\,d\eta
    + \int_{\xi}^{\infty}J_{P}(\eta)
    \frac{\partial}{\partial \kappa}(|Q(\eta)|^{2\sigma}Q(\eta))\,d\eta\\
  &= I_{P,\kappa,1}(\xi) + I_{P,\kappa,2}(\xi),\\
  I_{P,\epsilon}(\xi)
  &= \int_{\xi}^{\infty}J_{P,\epsilon}(\eta)|Q(\eta)|^{2\sigma}Q(\eta)\,d\eta
    + \int_{\xi}^{\infty}J_{P}(\eta)
    \frac{\partial}{\partial \epsilon}(|Q(\eta)|^{2\sigma}Q(\eta))\,d\eta\\
  &= I_{P,\epsilon,1}(\xi) + I_{P,\epsilon,2}(\xi).
\end{align*}
For \(I_{P,\real\gamma}\), \(I_{P,\imag\gamma}\), \(I_{P,\kappa,2}\)
and \(I_{P,\epsilon,2}\), it suffices to give bounds similar to those
for \(I_{P}\) in Lemma~\ref{lemma:I_P-I_E}.

\begin{lemma}\label{lemma:I-P-dgamma-dkappa-depsilon}
  Assume that \(\xi_{1} > 1\) and that
  \((2\sigma + 1)\normv - \frac{2}{\sigma} + d - 2 < 0\). Then, for
  \(\xi \geq \xi_{1}\) we have the bounds
  \begin{align*}
    |I_{P,\real\gamma}(\xi)|
    &\leq (2\sigma + 1)C_{I_{P}}\|Q\|_{\normv}^{2\sigma}\|Q_{\real\gamma}\|_\normv
      e^{-\real(c)\xi^{2}}\xi^{(2\sigma + 1)\normv - \frac{2}{\sigma} + d - 2},\\
    |I_{P,\imag\gamma}(\xi)|
    &\leq (2\sigma + 1)C_{I_{P}}\|Q\|_{\normv}^{2\sigma}\|Q_{\imag\gamma}\|_\normv
      e^{-\real(c)\xi^{2}}\xi^{(2\sigma + 1)\normv - \frac{2}{\sigma} + d - 2},\\
    \left|I_{P,\kappa,2}(\xi)\right|
    &\leq (2\sigma + 1)C_{I_{P}}\|Q\|_{\normv}^{2\sigma}\|Q_{\kappa}\|_{\normv}
      e^{-\real(c)\xi^{2}}\xi^{(2\sigma + 1)\normv - \frac{2}{\sigma} + d - 2},\\
    \left|I_{P,\epsilon,2}(\xi)\right|
    &\leq (2\sigma + 1)C_{I_{P}}\|Q\|_{\normv}^{2\sigma}\|Q_{\epsilon}\|_{\normv}
      e^{-\real(c)\xi^{2}}\xi^{(2\sigma + 1)\normv - \frac{2}{\sigma} + d - 2}
  \end{align*}
  with the constant \(C_{I_{P}}\) as in Lemma~\ref{lemma:I_P-I_E}.
\end{lemma}

\begin{proof}
  This follows in the same way as the bound for \(I_{P}\) in
  Lemma~\ref{lemma:I_P-I_E}, except that we use the bounds for the
  derivatives of \(|Q(\eta)|^{2\sigma}Q(\eta)\) with respect to the
  parameters from Lemma~\ref{lemma:norm-bounds}.
\end{proof}

In case of \(I_{P,\kappa,1}\) and \(I_{P,\epsilon,1}\), these type of
direct bounds are not possible when \(\epsilon = 0\). For
\(I_{P,\kappa,1}\) (the situation is similar for \(I_{P,\epsilon,1}\))
we would have the bound
\begin{equation*}
  \left|\int_{\xi}^{\infty}J_{P,\kappa}(\eta)|Q(\eta)|^{2\sigma}Q(\eta)\,d\eta\right|
  \leq C_{J_{P},\kappa}\|Q\|^{2\sigma + 1}\int_{\xi}^{\infty}e^{-\real(c)\eta^{2}}\eta^{(2\sigma + 1)\normv - \frac{2}{\sigma} + d - 1}\,d\eta.
\end{equation*}
In the supercritical case we have \(\frac{2}{\sigma} + d > 0\), so
this is only integrable when \(\real(c) > 0\), corresponding to
\(\epsilon > 0\). To get a convergent integral also when
\(\epsilon = 0\) we therefore need to take into account the
oscillations coming from \(J_{P,\kappa}\). This is done through
integration by parts in the following section.

\subsubsection{Bounds for \(I_{P}\) through integration by parts}
\label{sec:I_P-integration-by-parts-bounds}

To get refined bounds for \(I_{P}\) as well as bounds for
\(I_{P,\kappa,1}\) and \(I_{P,\epsilon,1}\) we need to account for the
exponential oscillations in the integrand, which is done through
integration by parts. In addition, when computing enclosures at
\(\xi = \xi_{1}\) in Section~\ref{sec:enclosures-xi-1} the expressions
coming from the integration by parts are useful for computing tighter
enclosures. For this reason we include formulas also for
\(I_{P,\real\gamma}\), \(I_{P,\imag\gamma}\), \(I_{P,\kappa,2}\) and
\(I_{P,\epsilon,2}\), even if they are not needed for the uniform
bounds. The idea is to use the identity
\begin{equation*}
  e^{-c\eta^{2}} = -\frac{1}{2c}\eta^{-1}\frac{\partial}{\partial \eta} e^{-c\eta^{2}}.
\end{equation*}
on the exponential factors appearing in the integrands. In case of
\(I_{P}\) we want to integrate by parts up to three times, giving us
the following lemma.

\begin{lemma}\label{lemma:I_P-expansion}
  Integrating by parts \(n\) times, we can, for \(1 \leq n \leq 3\),
  write \(I_P\) as
  \begin{equation}
    \label{eq:I_P-expansion}
    I_{P}(\xi) =
    B_{W}\left(
      \sum_{k = 1}^{n} \frac{1}{(2c)^{k}} I_{P,k}(\xi)
      + \frac{1}{(2c)^{n}}\hat{I}_{P,n + 1}(\xi)
    \right),
  \end{equation}
  where
  \begin{align*}
    I_{P,1}(\xi)
    &= e^{-c\xi^{2}}P(\xi)\xi^{d - 2}|Q(\xi)|^{2\sigma}Q(\xi),\\
    I_{P,2}(\xi)
    &= e^{-c\xi^{2}}\xi^{-1}
      \frac{\partial}{\partial \xi}\left(P(\xi)\xi^{d - 2}|Q(\xi)|^{2\sigma}Q(\xi)\right),\\
    I_{P,3}(\xi)
    &= e^{-c\xi^{2}}\xi^{-1}
      \frac{\partial}{\partial \xi}\left(\xi^{-1}\frac{\partial}{\partial \xi}\left(P(\xi)\xi^{d - 2}|Q(\xi)|^{2\sigma}Q(\xi)\right)\right)
  \end{align*}
  and
  \begin{align*}
    \hat{I}_{P,2}(\xi)
    &= \int_{\xi}^{\infty}e^{-c\eta^{2}}
      \frac{\partial}{\partial \eta}\left(P(\eta)\eta^{d - 2}|Q(\eta)|^{2\sigma}Q(\eta)\right)\,d\eta,\\
    \hat{I}_{P,3}(\xi)
    &= \int_{\xi}^{\infty}
      e^{-c\eta^{2}}\frac{\partial}{\partial \eta}\left(
      \eta^{-1}
      \frac{\partial}{\partial \eta}\left(P(\eta)\eta^{d - 2}|Q(\eta)|^{2\sigma}Q(\eta)\right)
      \right)\,d\eta,\\
    \hat{I}_{P,4}(\xi)
    &= \int_{\xi}^{\infty}
      e^{-c\eta^{2}}\frac{\partial}{\partial \eta}\left(
      \eta^{-1}
      \frac{\partial}{\partial \eta}\left(
      \eta^{-1}
      \frac{\partial}{\partial \eta}\left(P(\eta)\eta^{d - 2}|Q(\eta)|^{2\sigma}Q(\eta)\right)
      \right)\right)\,d\eta.
  \end{align*}
\end{lemma}

\begin{proof}
  Follows directly through iterative integration by parts, using that
  \(J_{P}(\eta) = B_{W}P(\eta)e^{-c\eta^{2}}\eta^{d - 1}\) and the
  identity
  \(e^{-c\eta^{2}} = -\frac{1}{2c}\eta^{-1}\frac{\partial}{\partial
    \eta} e^{-c\eta^{2}}\).
\end{proof}

For \(I_{P,\real\gamma}\), \(I_{P,\imag\gamma}\), \(I_{P,\kappa,2}\)
and \(I_{P,\epsilon,2}\) we only need one iteration of integration by
parts.

\begin{lemma}\label{lemma:I_P-derivatives-expansion-1}
  Integrating by parts we have, for \(n = 1\),
  \begin{align*}
    I_{P,\real\gamma}(\xi)
    &= B_{W}\left(
      \sum_{k = 1}^{n} \frac{1}{(2c)^{k}} I_{P,\real\gamma,k}(\xi)
      + \frac{1}{(2c)^{n}}\hat{I}_{P,\real\gamma,n + 1}(\xi)
      \right),\\
    I_{P,\imag\gamma}(\xi)
    &= B_{W}\left(
      \sum_{k = 1}^{n} \frac{1}{(2c)^{k}} I_{P,\imag\gamma,k}(\xi)
      + \frac{1}{(2c)^{n}}\hat{I}_{P,\imag\gamma,n + 1}(\xi)
      \right),\displaybreak[3]\\
    I_{P,\kappa,2}(\xi)
    &= B_{W}\left(
      \sum_{k = 1}^{n} \frac{1}{(2c)^{k}} I_{P,\kappa,2,k}(\xi)
      + \frac{1}{(2c)^{n}}\hat{I}_{P,\kappa,2,n + 1}(\xi)
      \right),\\
    I_{P,\epsilon,2}(\xi)
    &= B_{W}\left(
      \sum_{k = 1}^{n} \frac{1}{(2c)^{k}} I_{P,\epsilon,2,k}(\xi)
      + \frac{1}{(2c)^{n}}\hat{I}_{P,\epsilon,2,n + 1}(\xi)
      \right).
  \end{align*}
  The individual terms are the same as for \(I_{P}\) in
  Lemma~\ref{lemma:I_P-expansion}, except that one replaces
  \(|Q(\eta)|^{2\sigma}Q(\eta)\) by its derivative with respect to
  \(\real\gamma\), \(\imag\gamma\), \(\kappa\) or \(\epsilon\)
  respectively.
\end{lemma}

\begin{proof}
  Follows in the same way as Lemma~\ref{lemma:I_P-expansion}.
\end{proof}

For \(I_{P,\kappa,1}\) and \(I_{P,\epsilon,1}\) we need two iterations.

\begin{lemma}\label{lemma:I_P-derivatives-expansion-2}
  Integrating by parts twice we have, for \(n = 2\),
  \begin{align*}
    I_{P,\kappa,1}(\xi)
    &= \sum_{k = 1}^{n} \frac{1}{(2c)^{k}} I_{P,\kappa,1,k}(\xi)
      + \frac{1}{(2c)^{n}}\hat{I}_{P,\kappa,1,n + 1}(\xi),\\
    I_{P,\epsilon,1}(\xi)
    &= \sum_{k = 1}^{n} \frac{1}{(2c)^{k}} I_{P,\epsilon,1,k}(\xi)
      + \frac{1}{(2c)^{n}}\hat{I}_{P,\epsilon,1,n + 1}(\xi).
  \end{align*}
  Here
  \begin{align*}
    I_{P,\kappa,1,1}(\xi)
    &= e^{-c\xi^{2}}D(\xi)\xi^{d}|Q(\xi)|^{2\sigma}Q(\xi),\\
    I_{P,\kappa,1,2}(\xi)
    &= e^{-c\xi^{2}}\xi^{-1}
      \frac{\partial}{\partial \xi}\left(D(\xi)\xi^{d}|Q(\xi)|^{2\sigma}Q(\xi)\right),\\
    I_{P,\epsilon,1,1}(\xi)
    &= e^{-c\xi^{2}}H(\xi)\xi^{d}|Q(\xi)|^{2\sigma}Q(\xi),\\
    I_{P,\epsilon,1,2}(\xi)
    &= e^{-c\xi^{2}}\xi^{-1}
      \frac{\partial}{\partial \xi}\left(H(\xi)\xi^{d}|Q(\xi)|^{2\sigma}Q(\xi)\right)
  \end{align*}
  and
  \begin{align*}
    \hat{I}_{P,\kappa,1,3}(\xi)
    &= \int_{\xi}^{\infty}
      e^{-c\eta^{2}}\frac{\partial}{\partial \eta}\left(
      \eta^{-1}
      \frac{\partial}{\partial \eta}\left(D(\eta)\eta^{d}|Q(\eta)|^{2\sigma}Q(\eta)\right)
      \right)\,d\eta,\\
    \hat{I}_{P,\epsilon,1,3}(\xi)
    &= \int_{\xi}^{\infty}
      e^{-c\eta^{2}}\frac{\partial}{\partial \eta}\left(
      \eta^{-1}
      \frac{\partial}{\partial \eta}\left(H(\eta)\eta^{d}|Q(\eta)|^{2\sigma}Q(\eta)\right)
      \right)\,d\eta.
  \end{align*}
\end{lemma}

\begin{proof}
  Follows in the same way as Lemma~\ref{lemma:I_P-expansion}, using
  that \(J_{P,\kappa}(\xi) = e^{-c\xi^{2}}D(\xi)\xi^{d + 1}\) and
  \(J_{P,\epsilon}(\xi) = e^{-c\xi^{2}}H(\xi)\xi^{d + 1}\).
\end{proof}

Let us start by bounding the remainder terms in the above expansions,
i.e.\ \(\hat{I}_{P,2}\) and similar. Starting with \(I_{P}\) we have
the following bounds.

\begin{lemma}\label{lemma:I_P-remainder-bounds}
  Let \(\alpha = (2\sigma + 1)\normv - \frac{2}{\sigma} + d\) and
  assume that \(\alpha - 2 < 0\). For \(\xi \geq \xi_{1} > 1\) we have
  the bounds
  \begin{align*}
    |\hat{I}_{P,2}(\xi)|
    &\leq \left(
      \frac{C_{P'} + |d - 2|C_{P}}{|\alpha - 4|}\|Q\|_{\normv}^{2\sigma + 1}\xi^{-1}
      + \frac{(2\sigma + 1)C_{P}}{|\alpha - 3|}\|Q\|_{\normv}^{2\sigma}\|Q'\|_{\normv}
      \right)e^{-\real(c)\xi^{2}}\xi^{\alpha - 3}
  \end{align*}
  and
  \begin{align*}
    |\hat{I}_{P,3}(\xi)|
    &\leq \Bigg(
      \frac{C_{P''} + |2d - 5|C_{P'} + |(d - 2)(d - 4)|C_{P}}{|\alpha - 6|}\|Q\|_{\normv}^{2\sigma + 1}
      \xi^{-2}\\
    &\qquad + (2\sigma + 1)\frac{2C_{P'} + |2d - 5|C_{P}}{|\alpha - 5|}\|Q\|_{\normv}^{2\sigma}\|Q'\|_{\normv}
      \xi^{-1}\\
    &\qquad + \frac{(2\sigma + 1)C_{P}}{|\alpha - 4|}(2\sigma\|Q'\|_{\normv}^{2} + \|Q\|_{\normv}\|Q''\|_{\normv})\|Q\|_{\normv}^{2\sigma - 1}
      \Bigg)e^{-\real(c)\xi^{2}}\xi^{\alpha - 4},
  \end{align*}
  as well as
  \begin{align*}
    |\hat{I}_{P,4}(\xi)|
    &\leq \Bigg(
      \frac{C_{P'''} + |3d - 9|C_{P''} + |3d^{2} - 21d + 33|C_{P'} + |(d - 2)(d - 4)(d - 6)|C_{P}}{|\alpha - 8|}
      \|Q\|_{\normv}^{2\sigma + 1}\xi^{-3}\\
    &\qquad + (2\sigma + 1)\frac{3C_{P''} + |6d - 18|C_{P'} + |3d^{2} - 21d + 33|C_{P}}{|\alpha - 7|}
      \|Q\|_{\normv}^{2\sigma}\|Q'\|_{\normv}\xi^{-2}\\
    &\qquad + (2\sigma + 1)\frac{3C_{P'} + |3d - 9|C_{P}}{|\alpha - 6|}
      (2\sigma\|Q'\|_{\normv}^{2} + \|Q\|_{\normv}\|Q''\|_{\normv})\|Q\|_{\normv}^{2\sigma - 1}
      \xi^{-1}\\
    &\qquad + \frac{(2\sigma + 1)C_{P}}{|\alpha - 5|}
      (
      2\sigma(2\sigma - 1)\|Q'\|_{\normv}^{3}
      + 6\sigma\|Q\|_{\normv}\|Q'\|_{\normv}\|Q''\|_{\normv}
      + \|Q\|_{\normv}^{2}\|Q'''\|_{\normv}
      )\|Q\|_{\normv}^{2\sigma - 2}
      \Bigg)e^{-\real(c)\xi^{2}}\xi^{\alpha - 5}.
  \end{align*}
\end{lemma}

\begin{proof}
  For \(\hat{I}_{P,2}\) we have
  \begin{align*}
    \hat{I}_{P,2}(\xi)
    &= \int_{\xi}^{\infty}e^{-c\eta^{2}}
      P'(\eta)\eta^{d - 2}|Q(\eta)|^{2\sigma}Q(\eta)
      \,d\eta\\
    &\quad + (d - 2)\int_{\xi}^{\infty}e^{-c\eta^{2}}
      P(\eta)\eta^{d - 3}|Q(\eta)|^{2\sigma}Q(\eta)
      \,d\eta\\
    &\quad + \int_{\xi}^{\infty}e^{-c\eta^{2}}
      P(\eta)\eta^{d - 2}\frac{\partial}{\partial \eta}\left(|Q(\eta)|^{2\sigma}Q(\eta)\right)
      \,d\eta,
  \end{align*}
  The result follows by using the bounds from
  Lemmas~\ref{lemma:bounds-list} and~\ref{lemma:norm-bounds} and
  grouping the terms appropriately.

  For \(\hat{I}_{P,3}\) and \(\hat{I}_{P,4}\) we have
  \begin{align*}
    \hat{I}_{P,3}(\xi)
    &= \int_{\xi}^{\infty}e^{-c\eta^{2}}
      P''(\eta)\eta^{d - 3}|Q(\eta)|^{2\sigma}Q(\eta)
      \,d\eta\\
    &\quad + (2d - 5)\int_{\xi}^{\infty}e^{-c\eta^{2}}
      P'(\eta)\eta^{d - 4}|Q(\eta)|^{2\sigma}Q(\eta)
      \,d\eta\\
    &\quad + 2\int_{\xi}^{\infty}e^{-c\eta^{2}}
      P'(\eta)\eta^{d - 3}\frac{\partial}{\partial \eta}\left(|Q(\eta)|^{2\sigma}Q(\eta)\right)
      \,d\eta\\
    &\quad + (d - 2)(d - 4)\int_{\xi}^{\infty}e^{-c\eta^{2}}
      P(\eta)\eta^{d - 5}|Q(\eta)|^{2\sigma}Q(\eta)
      \,d\eta\\
    &\quad + (2d - 5)\int_{\xi}^{\infty}e^{-c\eta^{2}}
      P(\eta)\eta^{d - 4}\frac{\partial}{\partial \eta}\left(|Q(\eta)|^{2\sigma}Q(\eta)\right)
      \,d\eta\\
    &\quad + \int_{\xi}^{\infty}e^{-c\eta^{2}}
      P(\eta)\eta^{d - 3}\frac{\partial^{2}}{\partial \eta^{2}}\left(|Q(\eta)|^{2\sigma}Q(\eta)\right)
      \,d\eta
  \end{align*}
  and
  \begin{align*}
    \hat{I}_{P,4}(\xi)
    &= \int_{\xi}^{\infty}e^{-c\eta^{2}}
      P'''(\eta)\eta^{d - 4}|Q(\eta)|^{2\sigma}Q(\eta)
      \,d\eta\\
    &\quad + (3d - 9)\int_{\xi}^{\infty}e^{-c\eta^{2}}
      P''(\eta)\eta^{d - 5}|Q(\eta)|^{2\sigma}Q(\eta)
      \,d\eta\\
    &\quad + 3\int_{\xi}^{\infty}e^{-c\eta^{2}}
      P''(\eta)\eta^{d - 4}\frac{\partial}{\partial \eta}\left(|Q(\eta)|^{2\sigma}Q(\eta)\right)
      \,d\eta\\
    &\quad + (3d^{2} - 21d + 33)\int_{\xi}^{\infty}e^{-c\eta^{2}}
      P'(\eta)\eta^{d - 6}|Q(\eta)|^{2\sigma}Q(\eta)
      \,d\eta\\
    &\quad + (6d - 18)\int_{\xi}^{\infty}e^{-c\eta^{2}}
      P'(\eta)\eta^{d - 5}\frac{\partial}{\partial \eta}\left(|Q(\eta)|^{2\sigma}Q(\eta)\right)
      \,d\eta\\
    &\quad + 3\int_{\xi}^{\infty}e^{-c\eta^{2}}
      P'(\eta)\eta^{d - 4}\frac{\partial^{2}}{\partial \eta^{2}}\left(|Q(\eta)|^{2\sigma}Q(\eta)\right)
      \,d\eta\\
    &\quad + (d - 2)(d - 4)(d - 6)\int_{\xi}^{\infty}e^{-c\eta^{2}}
      P(\eta)\eta^{d - 7}|Q(\eta)|^{2\sigma}Q(\eta)
      \,d\eta\\
    &\quad + (3d^{2} - 21d + 33)\int_{\xi}^{\infty}e^{-c\eta^{2}}
      P(\eta)\eta^{d - 6}\frac{\partial}{\partial \eta}\left(|Q(\eta)|^{2\sigma}Q(\eta)\right)
      \,d\eta\\
    &\quad + (3d - 9)\int_{\xi}^{\infty}e^{-c\eta^{2}}
      P(\eta)\eta^{d - 5}\frac{\partial^{2}}{\partial \eta^{2}}\left(|Q(\eta)|^{2\sigma}Q(\eta)\right)
      \,d\eta\\
    &\quad + \int_{\xi}^{\infty}e^{-c\eta^{2}}
      P(\eta)\eta^{d - 4}\frac{\partial^{3}}{\partial \eta^{3}}\left(|Q(\eta)|^{2\sigma}Q(\eta)\right)
      \,d\eta.
  \end{align*}
  Similarly to \(\hat{I}_{P,2}\) these can be bounded termwise and the
  result follows by grouping similar terms together.
\end{proof}

For \(\hat{I}_{P,\real\gamma,2}\), \(\hat{I}_{P,\imag\gamma,2}\),
\(\hat{I}_{P,\kappa,2,2}\) and \(\hat{I}_{P,\epsilon,2,2}\) the bounds
are similar to those of \(\hat{I}_{P,2}\). We omit the proofs since
they follow in the same way as Lemma~\ref{lemma:I_P-remainder-bounds}.

\begin{lemma}\label{lemma:I_P-remainder-bounds-1}
  Let \(\alpha = (2\sigma + 1)\normv - \frac{2}{\sigma} + d\). Assume
  that \(\xi_{1} > 1\) and \(\alpha - 2 < 0\). For
  \(\xi \geq \xi_{1} > 1\) we have the bounds
  \begin{align*}
    |\hat{I}_{P,\real\gamma,2}(\xi)|
    &\leq \Bigg(
      (2\sigma + 1)\frac{C_{P'} + |d - 2|C_{P}}{|\alpha - 4|}
      \|Q\|_{\normv}^{2\sigma}\|Q_{\real\gamma}\|_{\normv}\xi^{-1}\\
    &\qquad + \frac{(2\sigma + 1)C_{P}}{|\alpha - 3|}
      (2\sigma\|Q'\|_{\normv}\|Q_{\real\gamma}\|_{\normv} + \|Q\|_{\normv}\|Q_{\real\gamma}'\|_{\normv})\|Q\|_{\normv}^{2\sigma - 1}
      \Bigg)e^{-\real(c)\xi^{2}}\xi^{\alpha - 3}.
  \end{align*}
  The functions \(\hat{I}_{P,\imag\gamma,2}\),
  \(\hat{I}_{P,\kappa,2,2}\) and \(\hat{I}_{P,\epsilon,2,2}\) are
  bounded by the same expression, replacing the derivative in
  \(\real\gamma\) with \(\imag\gamma\), \(\kappa\) and \(\epsilon\)
  respectively.
\end{lemma}

For \(\hat{I}_{P,\kappa,1,3}\) and \(\hat{I}_{P,\epsilon,1,3}\) we
have the following lemma. Again, we omit the proof since it is very
similar to that for the bound of \(\hat{I}_{P,3}\) in
Lemma~\ref{lemma:I_P-remainder-bounds}.

\begin{lemma}\label{lemma:I_P-remainder-bounds-2}
  Let \(\alpha = (2\sigma + 1)\normv - \frac{2}{\sigma} + d\). Assume
  that \(\xi_{1} > 1\) and \(\alpha - 2 < 0\). For
  \(\xi \geq \xi_{1} > 1\) we have the bounds
  \begin{align*}
    |\hat{I}_{P,\kappa,1,3}(\xi)|
    &\leq \Bigg(
      \frac{C_{D''} + |2d - 1|C_{D'} + |d(d - 2)|C_{D}}{|\alpha - 4|}
      \|Q\|_{\normv}^{2\sigma + 1}\xi^{-2}\\
    &\quad + (2\sigma + 1)\frac{2C_{D'} + |2d - 1|C_{D}}{|\alpha - 3|}
      \|Q\|_{\normv}^{2\sigma}\|Q'\|_{\normv}\xi^{-1}\\
    &\quad + \frac{(2\sigma + 1)C_{D}}{|\alpha - 2|}
      (2\sigma\|Q'\|_{\normv}^{2} + \|Q\|_{\normv}\|Q''\|_{\normv})\|Q\|_{\normv}^{2\sigma - 1}
      \Bigg)e^{-\real(c)\xi^{2}}\xi^{\alpha - 2}.
  \end{align*}
  For \(\hat{I}_{P,\epsilon,1,3}\) the bound is the same, except \(D\)
  is replaced with \(H\).
\end{lemma}

Let us next give the refined bounds for \(I_{P}\). We give two
separate bounds, one associated with one iteration of integration by
parts and the other with two. The first version is used to bootstrap
bounds for the norm of \(Q''\), which can then be used in the second
version.

\begin{lemma}\label{lemma:I-P-refined}
  Let \(\alpha = (2\sigma + 1)\normv - \frac{2}{\sigma} + d\). Assume
  that \(\xi_{1} > 1\) and \(\alpha - 2 < 0\). Then, for
  \(\xi \geq \xi_{1}\), we have the two bounds
  \begin{align*}
    |I_{P}(\xi)|
    &\leq (
      C_{I_{P},1,1}\|Q\|_{\normv}\xi^{-1}
      + C_{I_{P},1,2}\|Q'\|_{\normv}
      )\|Q\|_{\normv}^{2\sigma}e^{-\real(c)\xi^{2}}
      \xi^{\alpha - 3},\\
    |I_{P}(\xi)|
    &\leq \Big(
    C_{I_{P},2,1}\|Q\|_{\normv}^{2}
    + C_{I_{P},2,2}\|Q\|_{\normv}\|Q'\|_{\normv}\xi^{-1}\\
    &\qquad+ C_{I_{P},2,3}\|Q'\|_{\normv}^{2}
    + C_{I_{P},2,4}\|Q\|_{\normv}\|Q''\|_{\normv}
    \Big)\|Q\|_{\normv}^{2\sigma - 1}e^{-\real(c)\xi^{2}}
    \xi^{\alpha - 4}.
  \end{align*}
  The constants are given by
  \begin{align*}
    C_{I_{P},1,1}
    &= \frac{|B_{W}|}{|2c|}\left(
      C_{P}
      + \frac{C_{P'}}{|\alpha - 4|}
      + \frac{|d - 2|C_{P}}{|\alpha - 4|}
      \right),\\
    C_{I_{P},1,2}
    &= \frac{|B_{W}|}{|2c|}
      \frac{(2\sigma + 1)C_{P}}{|\alpha - 3|}
  \end{align*}
  and
  \begin{align*}
    C_{I_{P},2,1}
    &= \frac{|B_{W}|}{|2c|}C_{P} + \frac{|B_{W}|}{|4c^{2}|}\Bigg(
      C_{P'}
      + |d - 2|C_{P}
      + \frac{C_{P''} + |2d - 5|C_{P'} + |(d - 2)(d - 4)|C_{P}}{|\alpha - 6|}
      \Bigg)\xi_{1}^{-2},\\
    C_{I_{P},2,2}
    &= \frac{|B_{W}|}{|4c^{2}|}(2\sigma + 1)\left(
      C_{P}
      + \frac{2C_{P'} + |2d - 5|C_{P}}{|\alpha - 5|}
      \right),\\
    C_{I_{P},2,3}
    &= \frac{|B_{W}|}{|4c^{2}|}\frac{(2\sigma + 1)2\sigma C_{P}}{|\alpha - 4|},\\
    C_{I_{P},2,4}
    &= \frac{|B_{W}|}{|4c^{2}|}\frac{(2\sigma + 1)C_{P}}{|\alpha - 4|}.
  \end{align*}
\end{lemma}

\begin{proof}
  From Lemma~\ref{lemma:I_P-expansion} we have that \(I_{P}\) can be
  written as both
  \begin{equation*}
    I_{P}(\xi)
    = B_{W}\left(
      \frac{1}{2c}I_{P,1}(\xi)
      + \frac{1}{2c}\hat{I}_{P,2}(\xi)
    \right)
  \end{equation*}
  and
  \begin{equation*}
    I_{P}(\xi)
    = B_{W}\left(
      \frac{1}{2c}I_{P,1}(\xi)
      + \frac{1}{(2c)^{2}}I_{P,2}(\xi)
      + \frac{1}{(2c)^{2}}\hat{I}_{P,3}(\xi)
    \right),
  \end{equation*}
  corresponding to \(n = 1\) and \(n = 2\) respectively. We have
  bounds for \(\hat{I}_{P,2}\) and \(\hat{I}_{P,3}\) from
  Lemma~\ref{lemma:I_P-remainder-bounds}. For the other two terms we
  have, using Lemmas~\ref{lemma:bounds-list}
  and~\ref{lemma:norm-bounds},
  \begin{equation*}
    |I_{P,1}(\xi)| \leq C_{P}\|Q\|_{\normv}^{2\sigma + 1}e^{-\real(c)\xi^{2}}\xi^{\alpha - 4}
  \end{equation*}
  and
  \begin{equation*}
    |I_{P,2}(\xi)| \leq (
    (C_{P'} + |d - 2|C_{P})\|Q\|_{\normv}^{2\sigma + 1}\xi^{-1}
    + (2\sigma + 1)C_{P}\|Q\|_{\normv}^{2\sigma}\|Q'\|_{\normv}
    )e^{-\real(c)\xi^{2}}\xi^{\alpha - 5}.
  \end{equation*}
  Grouping the terms based on the derivatives of \(Q\) appearing gives
  us the result.
\end{proof}

In case of \(I_{P,\kappa,1}\) and \(I_{P,\epsilon,1}\), we have the
following lemma.

\begin{lemma}\label{lemma:I-P-dkappa-depsilon-1}
  Let \(\alpha = (2\sigma + 1)\normv - \frac{2}{\sigma} + d\). Assume
  that \(\xi_{1} > 1\) and \(\alpha - 2 < 0\). Then, for
  \(\xi \geq \xi_{1}\), we have the bounds
  \begin{align*}
    \left|I_{P,\kappa,1}(\xi)\right|
    &\leq \Big(
      C_{I_{P},\kappa,1,1}\|Q\|_{\normv}^{2}
      + C_{I_{P},\kappa,1,2}\|Q\|_{\normv}\|Q'\|_{\normv}\xi^{-1}\\
    &\qquad+ C_{I_{P},\kappa,1,3}\|Q'\|_{\normv}^{2}
      + C_{I_{P},\kappa,1,4}\|Q\|_{\normv}\|Q''\|_{\normv}
      \Big)\|Q\|_{\normv}^{2\sigma - 1}e^{-\real(c)\xi^{2}}\xi^{\alpha - 2},
  \end{align*}
  and
  \begin{align*}
    \left|I_{P,\epsilon,1}(\xi)\right|
    &\leq \Big(
      C_{I_{P},\epsilon,1,1}\|Q\|_{\normv}^{2}
      + C_{I_{P},\epsilon,1,2}\|Q\|_{\normv}\|Q'\|_{\normv}\xi^{-1}\\
    &\qquad+ C_{I_{P},\epsilon,1,3}\|Q'\|_{\normv}^{2}
      + C_{I_{P},\epsilon,1,4}\|Q\|_{\normv}\|Q''\|_{\normv}
      \Big)\|Q\|_{\normv}^{2\sigma - 1}e^{-\real(c)\xi^{2}}\xi^{\alpha - 2}.
  \end{align*}
  Here the constants for \(I_{P,\kappa,1}\) are given by
  \begin{align*}
    C_{I_{P},\kappa,1,1}
    &= \frac{C_{D}}{|2c|} + \frac{1}{|2c|^{2}}\left(C_{D'} + dC_{D} + \frac{
      C_{D''} + |2d - 1|C_{D'} + d|d - 2|C_{D}
      }{
      |\alpha - 4|
      }\right)\xi_{1}^{-2},\\
    C_{I_{P},\kappa,1,2}
    &= \frac{(2\sigma + 1)}{|2c|^{2}}\left(
      C_{D} + \frac{
      2C_{D'} + |2d - 1|C_{D}
      }{
      |\alpha - 3|
      }
      \right),\\
    C_{I_{P},\kappa,1,3}
    &= \frac{1}{|2c|^{2}}\frac{(2\sigma + 1)2\sigma C_{D}}{|\alpha - 2|},\\
    C_{I_{P},\kappa,1,4}
    &= \frac{1}{|2c|^{2}}\frac{(2\sigma + 1)C_{D}}{|\alpha - 2|}.
  \end{align*}
  The constants for \(I_{P,\epsilon,1}\) are the same, except \(D\)
  replaced with \(H\).
\end{lemma}

\begin{proof}
  We give the proof for \(I_{P,\kappa,1}\), the result for
  \(I_{P,\epsilon,1}\) follows in the same way.

  From Lemma~\ref{lemma:I_P-derivatives-expansion-2} we have
  \begin{equation}
    I_{P,\kappa,1}(\xi) =
    \frac{1}{2c}I_{P,\kappa,1,1}(\xi)
    + \frac{1}{(2c)^{2}}I_{P,\kappa,1,2}(\xi)
    + \frac{1}{(2c)^{2}}\hat{I}_{P,\kappa,1,3}(\xi).
  \end{equation}

  For \(I_{P,\kappa,1,1}(\xi)\) and \(I_{P,\kappa,1,2}(\xi)\) we get,
  using the bounds for \(D\) from Lemma~\ref{lemma:bounds-list} and
  for the norms from Lemma~\ref{lemma:norm-bounds},
  \begin{align*}
    |I_{P,\kappa,1,1}(\xi)|
    &\leq C_{D}e^{-\real(c)\xi^{2}}
    \|Q\|_{\normv}^{2\sigma + 1}\xi^{\alpha - 2},\\
    I_{P,\kappa,1,2}(\xi)
    &= e^{-c\xi^{2}}\Big(
      D'(\xi)\xi^{d - 1}|Q(\xi)|^{2\sigma}Q(\xi)\\
    &\qquad+ d D(\xi)\xi^{d - 2}|Q(\xi)|^{2\sigma}Q(\xi)
      + D(\xi)\xi^{d - 1}\frac{\partial}{\partial \xi}\left(|Q(\xi)|^{2\sigma}Q(\xi)\right)
      \Big)\\
    &\leq \left(
      (C_{D'} + dC_{D})\|Q\|_{\normv}\xi^{-1}
      + (2\sigma + 1)C_{D}\|Q'\|_{\normv}
      \right)
    \|Q\|_{\normv}^{2\sigma}e^{-\real(c)\xi^{2}}\xi^{\alpha - 3}.
  \end{align*}
  Combined with the bounds for \(\hat{I}_{P,\kappa,1,3}\) from
  Lemma~\ref{lemma:I_P-remainder-bounds-2}, this gives us the result.
\end{proof}

\subsubsection{Bounds on norms}
\label{sec:bounds-on-norms}
In this section we give bounds for the norms of derivatives of \(Q\).
These are based on the bound for \(\|Q\|_{\normv}\) from the
fixed-point argument in Proposition~\ref{prop:fixed-point} as well as
the bounds for \(I_{P}\), \(I_{E}\) and their derivatives from
Sections~\ref{sec:I_P-I_E-direct-bounds}
and~\ref{sec:I_P-integration-by-parts-bounds}.

We start by bounding the norms of \(Q'\), \(Q''\) and \(Q'''\).
\begin{lemma}\label{lemma:norm-dQ}
  Under the assumptions of Lemma~\ref{lemma:fixed-point-bounds}, we
  have
  \begin{align*}
    \left\|Q'\right\|_{\normv}
    &\leq C_{P'}|\gamma|\xi_{1}^{-\normv - 1}
      + C_{Q'}\|Q\|_{\normv}^{2\sigma + 1}\xi_{1}^{2\sigma\normv - 1},\\
    \left\|Q''\right\|_{\normv}
    &\leq C_{P''}|\gamma|\xi_{1}^{-\normv - 2}
      + \left(
      C_{Q'',1}\|Q\|_{\normv}\xi_{1}^{-1}
      + C_{Q'',2}\|Q'\|_{\normv}
      \right)\|Q\|_{\normv}^{2\sigma}\xi_{1}^{2\sigma\normv - 1},\\
    \left\|Q'''\right\|_{\normv}
    &\leq C_{P'''}|\gamma|\xi_{1}^{-\normv - 3}
      + \Big(
      C_{Q''',1}\|Q\|_{\normv}^{2}
      + C_{Q''',2}\|Q\|_{\normv}\|Q'\|_{\normv}\xi_{1}^{-1}\\
    &\qquad+ C_{Q''',3}\|Q'\|_{\normv}^{2}
      + C_{Q''',4}\|Q\|_{\normv}\|Q''\|_{\normv}
      \Big)\|Q\|_{\normv}^{2\sigma - 1}\xi_{1}^{2\sigma\normv - 1}.
  \end{align*}
  The constants are given by
  \begin{align*}
    C_{Q'}
    &= C_{E'}C_{I_{P}}
      + C_{P'}C_{I_{E}}\xi_{1}^{-2},\\
    C_{Q'',1}
    &= C_{E''}C_{I_{P},1,1}
      + \left|\frac{1 + i\delta}{1 - i\epsilon}\right|
      + C_{P''}C_{I_{E}}\xi_{1}^{-2},\\
    C_{Q'',2}
    &= C_{E''}C_{I_{P},1,2},\\
    C_{Q''',1}
    &= C_{E'''}C_{I_{P},2,1}
      + \left|\frac{1 + i\delta}{1 - i\epsilon}\right|(2|c| + (d - 1)\xi_{1}^{-2})
      + C_{P'''}C_{I_{E}}\xi_{1}^{-4},\\
    C_{Q''',2}
    &= C_{E'''}C_{I_{P},2,2}
      + (2\sigma + 1)\left|\frac{1 + i\delta}{1 - i\epsilon}\right|,\\
    C_{Q''',3}
    &= C_{E'''}C_{I_{P},2,3},\\
    C_{Q''',4}
    &= C_{E'''}C_{I_{P},2,4}.
  \end{align*}
\end{lemma}

\begin{proof}
  For \(Q'\) we have
  \begin{equation*}
    Q'(\xi) = \gamma P'(\xi) + P'(\xi)I_{E}(\xi) + E'(\xi)I_{P}(\xi).
  \end{equation*}
  Multiplying by \(\xi^{\frac{1}{\sigma} - \normv}\) and bounding it
  termwise, using Lemma~\ref{lemma:bounds-list} for \(P\) and \(E\)
  and Lemma~\ref{lemma:I_P-I_E} for \(I_{P}\) and \(I_{E}\) we get
  \begin{align*}
    \xi^{\frac{1}{\sigma} - \normv}|\gamma P'(\xi)|
    &\leq C_{P'}|\gamma|\xi_{1}^{-\normv - 1},\\
    \xi^{\frac{1}{\sigma} - \normv}|P'(\xi)I_{E}(\xi)|
    &\leq C_{P'}C_{I_{E}}
      \|Q\|_{\normv}^{2\sigma + 1}\xi_{1}^{2\sigma\normv - 3},\\
    \xi^{\frac{1}{\sigma} - \normv}|E'(\xi)I_{P}(\xi)|
    &\leq C_{E'}C_{I_{P}}
      \|Q\|_{\normv}^{2\sigma + 1}\xi_{1}^{2\sigma\normv - 1}.
  \end{align*}
  The sum gives a bound for
  \(\xi^{\frac{1}{\sigma} - \normv}|Q'(\xi)|\) for \(\xi \geq \xi_1\),
  which in turn gives us a bound for \(\|Q'\|_\normv\).

  Similarly, for \(Q''\) we have
  \begin{equation*}
    Q''(\xi) = \gamma P''(\xi)
    + P''(\xi)I_{E}(\xi)
    + E''(\xi)I_{P}(\xi)
    - \frac{1 + i\delta}{1 - i\epsilon}|Q(\xi)|^{2\sigma}Q(\xi).
  \end{equation*}
  For the individual terms we get, this time using the first bound in
  Lemma~\ref{lemma:I-P-refined} for \(I_{P}\),
  \begin{align*}
    \xi^{\frac{1}{\sigma} - \normv}|\gamma P''(\xi)|
    &\leq C_{P''}|\gamma|\xi_{1}^{-\normv - 2},\\
    \xi^{\frac{1}{\sigma} -\normv}|P''(\xi)I_{E}(\xi)|
    &\leq C_{P''}C_{I_{E}}\|Q\|_{\normv}^{2\sigma + 1}\xi_{1}^{2\sigma\normv - 4},\\
    \xi^{\frac{1}{\sigma} - \normv}|E''(\xi)I_{P}(\xi)|
    &\leq C_{E''}
      (C_{I_{P},1,1}\|Q\|_{\normv}\xi_{1}^{-1} + C_{I_{P},1,2}\|Q'\|_{\normv})\|Q\|_{\normv}^{2\sigma}\xi_{1}^{2\sigma\normv - 1},\\
    \xi^{\frac{1}{\sigma} - \normv}\left|\frac{1 + i\delta}{1 - i\epsilon}|Q(\xi)|^{2\sigma}Q(\xi)\right|
    &\leq \left|\frac{1 + i\delta}{1 - i\epsilon}\right|\|Q\|_{\normv}^{2\sigma + 1}\xi_{1}^{2\sigma\normv - 2}.
  \end{align*}
  Collecting the terms with \(\|Q\|_\normv^{2\sigma + 1}\) gives us
  \(C_{Q'',1}\) and collecting the ones with
  \(\|Q\|_\normv^{2\sigma}\|Q'\|_\normv\) gives us \(C_{Q'',2}\).

  For \(Q'''\) we have
  \begin{multline*}
    Q'''(\xi) = \gamma P'''(\xi)
    + P'''(\xi)I_{E}(\xi)
    + E'''(\xi)I_{P}(\xi)
    + \frac{1 + i\delta}{1 - i\epsilon}((d - 1)\xi^{-1} - 2c\xi)|Q(\xi)|^{2\sigma}Q(\xi)\\
    - \frac{1 + i\delta}{1 - i\epsilon}\frac{\partial}{\partial \xi}\left(|Q(\xi)|^{2\sigma}Q(\xi)\right).
  \end{multline*}
  For the individual terms we get, this time using the second bound in
  Lemma~\ref{lemma:I-P-refined} for \(I_{P}\), and
  Lemma~\ref{lemma:norm-bounds} for the derivative of the
  non-linearity,
  \begin{align*}
    \xi^{\frac{1}{\sigma} - \normv}|\gamma P'''(\xi)|
    &\leq C_{P'''}|\gamma|\xi_{1}^{-\normv - 3},\\
    \xi^{\frac{1}{\sigma} -\normv}|P'''(\xi)I_{E}(\xi)|
    &\leq C_{P'''}C_{I_{E}}\|Q\|_{\normv}^{2\sigma + 1}\xi_{1}^{2\sigma\normv - 5},\\
    \xi^{\frac{1}{\sigma} - \normv}|E'''(\xi)I_{P}(\xi)|
    &\leq C_{E'''}
      (
      C_{I_{P},2,1}\|Q\|_{\normv}^{2}
      + C_{I_{P},2,2}\|Q\|_{\normv}\|Q'\|_{\normv}\xi_{1}^{-1}\\
    &\quad+ C_{I_{P},2,3}\|Q'\|_{\normv}^{2}
      + C_{I_{P},2,4}\|Q\|_{\normv}\|Q''\|_{\normv}
      )\|Q\|_{\normv}^{2\sigma - 1}\xi_{1}^{2\sigma\normv - 1},\\
    \xi^{\frac{1}{\sigma} - \normv}\left|\frac{1 + i\delta}{1 - i\epsilon}((d - 1)\xi^{-1} - 2c\xi)|Q(\xi)|^{2\sigma}Q(\xi)\right|
    &\leq \left|\frac{1 + i\delta}{1 - i\epsilon}\right|
      (2|c| + (d - 1)\xi_{1}^{-2})\|Q\|_{\normv}^{2\sigma + 1}\xi_{1}^{2\sigma\normv - 1},\\
    \xi^{\frac{1}{\sigma} - \normv}\left|\frac{1 + i\delta}{1 - i\epsilon}\frac{\partial}{\partial \xi}\left(|Q(\xi)|^{2\sigma}Q(\xi)\right)\right|
    &\leq (2\sigma + 1)\left|\frac{1 + i\delta}{1 - i\epsilon}\right|\|Q\|_{\normv}^{2\sigma}\|Q'\|_{\normv}\xi_{1}^{2\sigma\normv - 2}.
  \end{align*}
  Collecting the terms with \(\|Q\|_\normv^{2\sigma + 1}\) gives us
  \(C_{Q''',1}\), the ones with
  \(\|Q\|_\normv^{2\sigma}\|Q'\|_\normv\) gives us \(C_{Q''',2}\), the
  ones with \(\|Q\|_\normv^{2\sigma - 1}\|Q'\|_\normv^2\) gives us
  \(C_{Q''',3}\) and the ones with
  \(\|Q\|_\normv^{2\sigma}\|Q''\|_\normv\) gives us \(C_{Q''',4}\).
\end{proof}

In case of \(Q_{\real\gamma}\), \(Q_{\imag\gamma}\), \(Q_{\kappa}\)
and \(Q_{\epsilon}\), all bounds follow a similar strategy:
differentiating the fixed-point equation with respect to the parameter
and solving for the associated derivative of \(Q\).

For \(Q_{\real\gamma}\) and \(Q_{\imag\gamma}\) we get the following
lemma.
\begin{lemma}\label{lemma:norm-Q-dgamma}
  Under the assumptions of Lemma~\ref{lemma:fixed-point-bounds}, we
  have
  \begin{equation*}
    \left\|Q_{\real\gamma}\right\|_{\normv}, \left\|Q_{\imag\gamma}\right\|_{\normv}
    \leq \frac{C_{P}\xi_{1}^{-\normv}}{1 - (2\sigma + 1)C_{T}\xi_{1}^{2\sigma\normv - 2}\|Q\|_{\normv}^{2\sigma}},
  \end{equation*}
  as long as the denominator is positive.
\end{lemma}

\begin{proof}
  Following the same approach as in
  Lemma~\ref{lemma:fixed-point-bounds}, for \(Q_{\real\gamma}\) we get
  \begin{equation*}
    \|Q_{\real\gamma}\|_{\normv}
    \leq C_{P}\xi_{1}^{-\normv} + (2\sigma + 1)C_{T}\xi_{1}^{-2 + 2\sigma\normv}\|Q\|_{\normv}^{2\sigma}\|Q_{\real\gamma}\|_{\normv}.
  \end{equation*}
  Solving for \(\|Q_{\real\gamma}\|_{\normv}\) gives us the result, as
  long as the factor in front of \(\|Q_{\real\gamma}\|_{\normv}\) is
  positive. The proof for \(Q_{\imag\gamma}\) is identical.
\end{proof}

For \(Q_{\kappa}\), the approach is the same, but slightly more work
is required due to the more involved bounds. Note that in this case
we have to assume that \(\normv > 0\), instead of just
\(\normv \geq 0\).
\begin{lemma}\label{lemma:norm-Q-dkappa}
  Under the assumptions of Lemma~\ref{lemma:I_P-I_E} as well as
  \(\normv > 0\), we have
  \begin{align*}
    \left\|Q_{\kappa}\right\|_{\normv}
    &\leq \frac{1}{
      1
      - C_{Q_{\kappa},6}\|Q\|_{\normv}^{2\sigma}
      }\Big(
      C_{Q_{\kappa},1}|\gamma|\\
    &\qquad+ \Big(
      C_{Q_{\kappa},2}\|Q\|_{\normv}^{2}
      + C_{Q_{\kappa},3}\|Q\|_{\normv}\|Q'\|_{\normv}
      + C_{Q_{\kappa},4}\|Q'\|_{\normv}^{2}
      + C_{Q_{\kappa},5}\|Q\|_{\normv}\|Q''\|_{\normv}
      \Big)\|Q\|_{\normv}^{2\sigma - 1}
      \Big),
  \end{align*}
  as long as the denominator is positive. The constants are given by
  \begin{align*}
    C_{Q_{\kappa},1}
    &= C_{P,\kappa}\frac{e^{-1}}{\normv}\\
    C_{Q_{\kappa},2}
    &= \left(
      C_{P,\kappa}C_{I_{E}}\frac{e^{-1}}{\normv}\xi_{1}^{\normv}
    + C_{P}C_{I_{E},\kappa}
    + C_{E,\kappa}C_{I_{P},1,1}
    + C_{E}C_{I_{P},\kappa,1,1}\right)\xi_{1}^{2\sigma\normv - 2},\\
    C_{Q_{\kappa},3}
    &= (C_{E,\kappa}C_{I_{P},1,2}
      + C_{E}C_{I_{P},\kappa,1,2}\xi_{1}^{-2}
      )\xi_{1}^{2\sigma\normv - 1},\\
    C_{Q_{\kappa},4}
    &= C_{E}C_{I_{P},\kappa,1,3}\xi_{1}^{2\sigma\normv - 2},\\
    C_{Q_{\kappa},5}
    &= C_{E}C_{I_{P},\kappa,1,4}\xi_{1}^{2\sigma\normv - 2},\\
    C_{Q_{\kappa},6}
    &= (2\sigma + 1)(
      C_{P}C_{I_{E}}
      + C_{E}C_{I_{P}}
      )\xi_{1}^{2\sigma\normv - 2}
  \end{align*}
\end{lemma}

\begin{proof}
  For \(Q_\kappa\) we have
  \begin{equation*}
    Q_{\kappa}(\xi)
    = \gamma P_{\kappa}(\xi)
    + P_{\kappa}(\xi)I_{E}(\xi)
    + P(\xi)I_{E,\kappa}(\xi)
    + E_{\kappa}(\xi)I_{P}(\xi)
    + E(\xi)I_{P,\kappa}(\xi).
  \end{equation*}
  Similarly to Lemma~\ref{lemma:norm-dQ} we multiply with
  \(\xi^{\frac{1}{\sigma} - \normv}\), take the absolute value and
  split into terms. Using Lemma~\ref{lemma:bounds-list} for \(P\) and
  \(E\), Lemma~\ref{lemma:I_P-I_E} for \(I_{E}\), the first bound in
  Lemma~\ref{lemma:I-P-refined} for \(I_{P}\),
  Lemma~\ref{lemma:I-E-bounds} for \(I_{E,\kappa}\) and
  Lemmas~\ref{lemma:I-P-dgamma-dkappa-depsilon}
  and~\ref{lemma:I-P-dkappa-depsilon-1} for \(I_{P,\kappa}\), we get
  the bounds
  \begin{align*}
    \xi^{\frac{1}{\sigma} - \normv}\left|\gamma P_{\kappa}(\xi)\right|
    &\leq C_{P,\kappa}|\gamma|\log(\xi)\xi^{-\normv},\\
    \xi^{\frac{1}{\sigma} - \normv}\left|P_{\kappa}(\xi)I_{E}(\xi)\right|
    &\leq C_{P,\kappa}C_{I_{E}}
      \|Q\|_{\normv}^{2\sigma + 1}\log(\xi)\xi^{-\normv}\xi_{1}^{(2\sigma + 1)\normv - 2},\\
    \xi^{\frac{1}{\sigma} - \normv}\left|P(\xi)I_{E,\kappa}(\xi)\right|
    &\leq C_{P}
      (
      C_{I_{E},\kappa}\|Q\|_{\normv}
      + (2\sigma + 1)C_{I_{E}}\left\|Q_{\kappa}\right\|_{\normv}
      )\|Q\|_{\normv}^{2\sigma}\xi^{-\normv}\xi_{1}^{(2\sigma + 1)\normv - 2},\\
    \xi^{\frac{1}{\sigma} - \normv}\left|E_{\kappa}(\xi)I_{P}(\xi)\right|
    &\leq C_{E,\kappa}
      (C_{I_{P},1,1}\|Q\|_{\normv}\xi^{-1} + C_{I_{P},1,2}\|Q'\|_{\normv})\|Q\|_{\normv}^{2\sigma}
      \xi^{2\sigma\normv - 1},\\
    \xi^{\frac{1}{\sigma} - \normv}\left|E(\xi)I_{P,\kappa}(\xi)\right|
    &\leq C_{E}\Big(C_{I_{P},\kappa,1,1}\|Q\|_{\normv}^{2}
    + C_{I_{P},\kappa,1,2}\|Q\|_{\normv}\|Q'\|_{\normv}\xi^{-1}
      + C_{I_{P},\kappa,1,3}\|Q'\|_{\normv}^{2}\\
    &\qquad+ C_{I_{P},\kappa,1,4}\|Q\|_{\normv}\|Q''\|_{\normv}
      + (2\sigma + 1)C_{I_{P}}\|Q\|_{\normv}\|Q_{\kappa}\|_{\normv}\\
    &\quad\Big)\|Q\|_{\normv}^{2\sigma - 1}\xi^{2\sigma\normv - 2}.
  \end{align*}
  To get bounds that are uniform in \(\xi\) we note that for
  \(\normv > 0\), \(\log(\xi)\xi^{-\normv}\) attains its maximum at
  \(\xi = e^{1 / \normv}\), where it is given by
  \(\frac{e^{-1}}{\normv}\). The rest of the terms are bounded by
  their value at \(\xi_{1}\). Combining these uniform bounds and
  collecting all parts with \(\|Q_\kappa\|_\normv\), we have
  \begin{align*}
    \|Q_{\kappa}\|_{\normv}
    &\leq C_{P,\kappa}|\gamma|\frac{e^{-1}}{\normv}\\
    &\quad+ C_{P,\kappa}C_{I_{E}}
      \|Q\|_{\normv}^{2\sigma + 1}\frac{e^{-1}}{\normv}\xi_{1}^{(2\sigma + 1)\normv - 2}\\
    &\quad+ C_{P}C_{I_{E},\kappa}
      \|Q\|_{\normv}^{2\sigma + 1}\xi_{1}^{2\sigma\normv - 2}\\
    &\quad+ C_{E,\kappa}
      (C_{I_{P},1,1}\|Q\|_{\normv}\xi_{1}^{-1} + C_{I_{P},1,2}\|Q'\|_{\normv})\|Q\|_{\normv}^{2\sigma}
      \xi_{1}^{2\sigma\normv - 1}\\
    &\quad+ C_{E}\Big(C_{I_{P},\kappa,1,1}\|Q\|_{\normv}^{2}
    + C_{I_{P},\kappa,1,2}\|Q\|_{\normv}\|Q'\|_{\normv}\xi_{1}^{-1}
      + C_{I_{P},\kappa,1,3}\|Q'\|_{\normv}^{2}\\
    &\qquad+ C_{I_{P},\kappa,1,4}\|Q\|_{\normv}\|Q''\|_{\normv}
      \Big)\|Q\|_{\normv}^{2\sigma - 1}\xi_{1}^{2\sigma\normv - 2}\\
    &\quad+ (2\sigma + 1)\left(
      C_{P}C_{I_{E}}
      + C_{E}C_{I_{P}}
      \right)
      \|Q\|_{\normv}^{2\sigma}\|Q_{\kappa}\|_{\normv}\xi_{1}^{2\sigma\normv - 2}.
  \end{align*}
  We can now solve for \(\|Q_{\kappa}\|_{\normv}\), giving us
  the result.
\end{proof}

The bound for \(Q_{\epsilon}\) is similar to that for \(Q_{\kappa}\).
\begin{lemma}\label{lemma:norm-Q-depsilon}
  Under the assumptions of Lemma~\ref{lemma:I_P-I_E} as well as
  \(\normv > 0\), we have
  \begin{align*}
    \left\|Q_{\epsilon}\right\|_{\normv}
    &\leq \frac{1}{
      1
      - C_{Q_{\epsilon},6}\|Q\|_{\normv}^{2\sigma}
      }\Big(
      C_{Q_{\epsilon},1}|\gamma|\\
    &\qquad+ \left(
      C_{Q_{\epsilon},2}\|Q\|_{\normv}^{2}
      + C_{Q_{\epsilon},3}\|Q\|_{\normv}\|Q'\|_{\normv}
      + C_{Q_{\epsilon},4}\|Q'\|_{\normv}^{2}
      + C_{Q_{\epsilon},5}\|Q\|_{\normv}\|Q''\|_{\normv}
      \right)\|Q\|_{\normv}^{2\sigma - 1}
      \Big),
  \end{align*}
  as long as the denominator is positive. The constants are given by
  \begin{align*}
    C_{Q_{\epsilon},1}
    &= C_{P,\epsilon}\xi_{1}^{-\normv}\\
    C_{Q_{\epsilon},2}
    &= \left(
      C_{P,\epsilon}C_{I_{E}}
    + C_{P}C_{I_{E},\epsilon}
    + C_{E,\epsilon}C_{I_{P},1,1}
    + C_{E}C_{I_{P},\epsilon,1,1}\right)\xi_{1}^{2\sigma\normv - 2},\\
    C_{Q_{\epsilon},3}
    &= (C_{E,\epsilon}C_{I_{P},1,2}
      + C_{E}C_{I_{P},\epsilon,1,2}\xi_{1}^{-2})\xi_{1}^{2\sigma\normv - 1},\\
    C_{Q_{\epsilon},4}
    &= C_{E}C_{I_{P},\epsilon,1,3}\xi_{1}^{2\sigma\normv - 2},\\
    C_{Q_{\epsilon},5}
    &= C_{E}C_{I_{P},\epsilon,1,4}\xi_{1}^{2\sigma\normv - 2},\\
    C_{Q_{\epsilon},6}
    &= (2\sigma + 1)(
      C_{P}C_{I_{E}}
      + C_{E}C_{I_{P}}
      )\xi_{1}^{2\sigma\normv - 2}
  \end{align*}
\end{lemma}
\begin{proof}
  For \(Q_\epsilon\) we have
  \begin{equation*}
    Q_{\epsilon}(\xi)
    = \gamma P_{\epsilon}(\xi)
    + P_{\epsilon}(\xi)I_{E}(\xi)
    + P(\xi)I_{E,\epsilon}(\xi)
    + E_{\epsilon}(\xi)I_{P}(\xi)
    + E(\xi)I_{P,\epsilon}(\xi).
  \end{equation*}
  Following the same approach as in Lemma~\ref{lemma:norm-Q-dkappa} we
  get the termwise bounds
  \begin{align*}
    \xi^{\frac{1}{\sigma} - \normv}\left|\gamma P_{\epsilon}(\xi)\right|
    &\leq C_{P,\epsilon}|\gamma|\xi^{-\normv},\\
    \xi^{\frac{1}{\sigma} - \normv}\left|P_{\epsilon}(\xi)I_{E}(\xi)\right|
    &\leq C_{P,\epsilon}C_{I_{E}}
      \|Q\|_{\normv}^{2\sigma + 1}\xi^{-\normv}\xi_{1}^{(2\sigma + 1)\normv - 2},\\
    \xi^{\frac{1}{\sigma} - \normv}\left|P(\xi)I_{E,\epsilon}(\xi)\right|
    &\leq C_{P}
      (
      C_{I_{E},\epsilon}\|Q\|_{\normv}
      + (2\sigma + 1)C_{I_{E}}\left\|Q_{\epsilon}\right\|_{\normv}
      )\|Q\|_{\normv}^{2\sigma}\xi^{-\normv}\xi_{1}^{(2\sigma + 1)\normv - 2},\\
    \xi^{\frac{1}{\sigma} - \normv}\left|E_{\epsilon}(\xi)I_{P}(\xi)\right|
    &\leq C_{E,\epsilon}
      (C_{I_{P},1,1}\|Q\|_{\normv}\xi^{-1} + C_{I_{P},1,2}\|Q'\|_{\normv})\|Q\|_{\normv}^{2\sigma}
      \xi^{2\sigma\normv - 1},\\
    \xi^{\frac{1}{\sigma} - \normv}\left|E(\xi)I_{P,\epsilon}(\xi)\right|
    &\leq C_{E}\Big(C_{I_{P},\epsilon,1,1}\|Q\|_{\normv}^{2}
    + C_{I_{P},\epsilon,1,2}\|Q\|_{\normv}\|Q'\|_{\normv}\xi^{-1}
      + C_{I_{P},\epsilon,1,3}\|Q'\|_{\normv}^{2}\\
    &\qquad+ C_{I_{P},\epsilon,1,4}\|Q\|_{\normv}\|Q''\|_{\normv}
      + (2\sigma + 1)C_{I_{P}}\|Q\|_{\normv}\|Q_{\epsilon}\|_{\normv}\\
    &\quad\Big)\|Q\|_{\normv}^{2\sigma - 1}\xi^{2\sigma\normv - 2}.
  \end{align*}
  In this case, all the terms are bounded by their value at
  \(\xi_{1}\), and the result follows using the same approach as in
  Lemma~\ref{lemma:norm-Q-dkappa}.
\end{proof}

What remains is bounding \(Q_{\real\gamma}'\), \(Q_{\imag\gamma}'\),
\(Q_{\kappa}'\) and \(Q_{\epsilon}'\). The bounds for
\(Q_{\real\gamma}'\) and \(Q_{\imag\gamma}'\) are straightforward, and
given in the following lemma.
\begin{lemma}\label{lemma:norm-Q-dgamma-dxi}
  Under the assumptions of Lemma~\ref{lemma:I_P-I_E}, we have
  \begin{align*}
    \left\|Q_{\real\gamma}'\right\|_{\normv}
    &\leq C_{P'}\xi_{1}^{-\normv - 1}
      + (2\sigma + 1)C_{Q'}\|Q\|_{\normv}^{2\sigma}\|Q_{\real\gamma}\|_{\normv}\xi_{1}^{2\sigma\normv - 1},\\
    \left\|Q_{\imag\gamma}'\right\|_{\normv}
    &\leq C_{P'}\xi_{1}^{-\normv - 1}
      + (2\sigma + 1)C_{Q'}\|Q\|_{\normv}^{2\sigma}\|Q_{\imag\gamma}\|_{\normv}\xi_{1}^{2\sigma\normv - 1},
  \end{align*}
  with \(C_{Q'}\) as in Lemma~\ref{lemma:norm-dQ}.
\end{lemma}

\begin{proof}
  The result follows in the same way as Lemma~\ref{lemma:norm-dQ},
  using that
  \begin{equation*}
    Q_{\real\gamma}'(\xi)
    = P'(\xi)
    + P'(\xi)I_{E,\real\gamma}(\xi)
    + E'(\xi)I_{P,\real\gamma}(\xi),
  \end{equation*}
  and similarly for \(Q_{\imag\gamma}'(\xi)\).
\end{proof}

The bound for \(Q_{\kappa}'\) requires a bit more work.
\begin{lemma}\label{lemma:norm-Q-dkappa-dxi}
  Under the assumptions of Lemma~\ref{lemma:I_P-I_E} as well as
  \(\xi_1 \geq e\), we have
  \begin{align*}
    \left\|Q_{\kappa}'\right\|_{\normv}
    &\leq C_{P',\kappa}\log(\xi_{1})|\gamma|\xi_{1}^{-\normv - 1}
      + \Bigg(
      C_{Q',\kappa,1}\|Q\|_{\normv}^{2}
      + C_{Q',\kappa,2}\|Q\|_{\normv}\|Q_{\kappa}\|_{\normv}\\
    &\qquad+ C_{Q',\kappa,3}\|Q\|_{\normv}\|Q'\|_{\normv}
      + C_{Q',\kappa,4}\|Q'\|_{\normv}^{2}
      + C_{Q',\kappa,5}\|Q\|_{\normv}\|Q''\|_{\normv}
      \Bigg)\|Q\|_{\normv}^{2\sigma - 1}.
  \end{align*}
  Here
  \begin{align*}
    C_{Q',\kappa,1}
    &= \Big(
      C_{P',\kappa}C_{I_{E}}\log(\xi_{1})\xi_{1}^{-2}
      + C_{P'}C_{I_{E},\kappa}\xi_{1}^{-2}
      + C_{E',\kappa}C_{I_{P},2,1}
      + C_{E'}C_{I_{P},\kappa,1,1}
      \Big)\xi_{1}^{2\sigma\normv - 1},\\
    C_{Q',\kappa,2}
    &=
      (2\sigma + 1)(
      C_{P'}C_{I_{E}}\xi_{1}^{-2}
      + C_{E'}C_{I_{P}})\xi_{1}^{2\sigma\normv - 1},\\
    C_{Q',\kappa,3}
    &=
      (
      C_{E',\kappa}C_{I_{P},2,2}
      + C_{E'}C_{I_{P},\kappa,1,2}
      )\xi_{1}^{2\sigma\normv - 2},\\
    C_{Q',\kappa,4}
    &=
      (
      C_{E',\kappa}C_{I_{P},2,3}
      + C_{E'}C_{I_{P},\kappa,1,3}
      )\xi_{1}^{2\sigma\normv - 1},\\
    C_{Q',\kappa,5}
    &=
      (
      C_{E',\kappa}C_{I_{P},2,4}
      + C_{E'}C_{I_{P},\kappa,1,4}
      )\xi_{1}^{2\sigma\normv - 1}.
  \end{align*}
\end{lemma}

\begin{proof}
  We have
  \begin{equation*}
    Q_{\kappa}'(\xi)
    = \gamma P_{\kappa}'(\xi)
    + P_{\kappa}'(\xi)I_{E}(\xi)
    + P'(\xi)I_{E,\kappa}(\xi)
    + E_{\kappa}'(\xi)I_{P}(\xi)
    + E'(\xi)I_{P,\kappa}(\xi).
  \end{equation*}
  Similarly to Lemma~\ref{lemma:norm-Q-dkappa}, we get the termwise
  bounds
  \begin{align*}
    \xi^{\frac{1}{\sigma} - \normv}|\gamma P_{\kappa}'(\xi)|
    &\leq C_{P',\kappa}|\gamma|\log(\xi)\xi^{-\normv - 1},\\
    \xi^{\frac{1}{\sigma} - \normv}|P_{\kappa}'(\xi)I_{E}(\xi)|
    &\leq C_{P',\kappa}C_{I_{E}}
      \|Q\|_{\normv}^{2\sigma + 1}\log(\xi)\xi^{-\normv - 1}
      \xi_{1}^{(2\sigma + 1)\normv - 2},\\
    \xi^{\frac{1}{\sigma} - \normv}|P'(\xi)I_{E,\kappa}(\xi)|
    &\leq C_{P'}(
      C_{I_{E},\kappa}\|Q\|_{\normv}
      + (2\sigma + 1)C_{I_{E}}\left\|Q_{\kappa}\right\|_{\normv}
      )\|Q\|_{\normv}^{2\sigma}\xi^{-\normv - 1}\xi_{1}^{(2\sigma + 1)\normv - 2},\\
    \xi^{\frac{1}{\sigma} - \normv}|E_{\kappa}'(\xi)I_{P}(\xi)|
    &\leq C_{E',\kappa}\Big(
      C_{I_{P},2,1}\|Q\|_{\normv}^{2}
      + C_{I_{P},2,2}\|Q\|_{\normv}\|Q'\|_{\normv}\xi^{-1}\\
    &\qquad+ C_{I_{P},2,3}\|Q'\|_{\normv}^{2}
      + C_{I_{P},2,4}\|Q\|_{\normv}\|Q''\|_{\normv}
      \Big)\|Q\|_{\normv}^{2\sigma - 1}
      \xi^{2\sigma\normv - 1},\\
    \xi^{\frac{1}{\sigma} - \normv}|E'(\xi)I_{P,\kappa}(\xi)|
    &\leq C_{E'}\Big(C_{I_{P},\kappa,1,1}\|Q\|_{\normv}^{2}
      + C_{I_{P},\kappa,1,2}\|Q\|_{\normv}\|Q'\|_{\normv}\xi^{-1}
      + C_{I_{P},\kappa,1,3}\|Q'\|_{\normv}^{2}\\
    &\qquad+ C_{I_{P},\kappa,1,4}\|Q\|_{\normv}\|Q''\|_{\normv}
      + (2\sigma + 1)C_{I_{P}}\|Q\|_{\normv}\|Q_{\kappa}\|_{\normv}\\
    &\quad\Big)\|Q\|_{\normv}^{2\sigma - 1}\xi^{2\sigma\normv - 1}.
  \end{align*}
  For \(\xi \geq e\), \(\log(\xi)\xi^{-1}\) is
  decreasing; this is enough to show that all terms are decreasing in
  \(\xi\) and hence bounded above by their value at \(\xi_{1}\).
  Collecting the bounds gives us the result.
\end{proof}

Finally, for \(Q_{\epsilon}'\) we have the following lemma.
\begin{lemma}\label{lemma:norm-Q-depsilon-dxi}
  Under the assumptions of Lemma~\ref{lemma:I_P-I_E}, we have
  \begin{align*}
    \left\|Q_{\epsilon}'\right\|_{\normv}
    &\leq C_{P',\epsilon}|\gamma|\xi_{1}^{-\normv - 1}
      + \Bigg(
      C_{Q',\epsilon,1}\|Q\|_{\normv}^{2}
      + C_{Q',\epsilon,2}\|Q\|_{\normv}\|Q_{\epsilon}\|_{\normv}\\
    &\qquad+ C_{Q',\epsilon,3}\|Q\|_{\normv}\|Q'\|_{\normv}
      + C_{Q',\epsilon,4}\|Q'\|_{\normv}^{2}
      + C_{Q',\epsilon,5}\|Q\|_{\normv}\|Q''\|_{\normv}
      \Bigg)\|Q\|_{\normv}^{2\sigma - 1}.
  \end{align*}
  Here
  \begin{align*}
    C_{Q',\epsilon,1}
    &= \Big(
      C_{P',\epsilon}C_{I_{E}}\xi_{1}^{-2}
      + C_{P'}C_{I_{E},\epsilon}\xi_{1}^{-2}
      + C_{E',\epsilon}C_{I_{P},2,1}
      + C_{E'}C_{I_{P},\epsilon,1,1}
      \Big)\xi_{1}^{2\sigma\normv - 1},\\
    C_{Q',\epsilon,2}
    &=
      (2\sigma + 1)(
      C_{P'}C_{I_{E}}\xi_{1}^{-2}
      + C_{E'}C_{I_{P}})\xi_{1}^{2\sigma\normv - 1},\\
    C_{Q',\epsilon,3}
    &=
      (
      C_{E',\epsilon}C_{I_{P},2,2}
      + C_{E'}C_{I_{P},\epsilon,1,2}
      )\xi_{1}^{2\sigma\normv - 2},\\
    C_{Q',\epsilon,4}
    &=
      (
      C_{E',\epsilon}C_{I_{P},2,3}
      + C_{E'}C_{I_{P},\epsilon,1,3}
      )\xi_{1}^{2\sigma\normv - 1},\\
    C_{Q',\epsilon,5}
    &=
      (
      C_{E',\epsilon}C_{I_{P},2,4}
      + C_{E'}C_{I_{P},\epsilon,1,4}
      )\xi_{1}^{2\sigma\normv - 1}.
  \end{align*}
\end{lemma}

\begin{proof}
  We have
  \begin{equation*}
    Q_{\epsilon}'(\xi)
    = \gamma P_{\epsilon}'(\xi)
    + P_{\epsilon}'(\xi)I_{E}(\xi)
    + P'(\xi)I_{E,\epsilon}(\xi)
    + E_{\epsilon}'(\xi)I_{P}(\xi)
    + E'(\xi)I_{P,\epsilon}(\xi).
  \end{equation*}
  Similarly to Lemma~\ref{lemma:norm-Q-depsilon} we get the termwise
  bounds
  \begin{align*}
    \xi^{\frac{1}{\sigma} - \normv}|\gamma P_{\epsilon}'(\xi)|
    &\leq C_{P',\epsilon}|\gamma|\xi^{-\normv - 1},\\
    \xi^{\frac{1}{\sigma} - \normv}|P_{\epsilon}'(\xi)I_{E}(\xi)|
    &\leq C_{P',\epsilon}C_{I_{E}}
      \|Q\|_{\normv}^{2\sigma + 1}\xi^{-\normv - 1}
      \xi_{1}^{(2\sigma + 1)\normv - 2},\\
    \xi^{\frac{1}{\sigma} - \normv}|P'(\xi)I_{E,\epsilon}(\xi)|
    &\leq C_{P'}(
      C_{I_{E},\epsilon}\|Q\|_{\normv}
      + (2\sigma + 1)C_{I_{E}}\left\|Q_{\epsilon}\right\|_{\normv}
      )\|Q\|_{\normv}^{2\sigma}\xi^{-\normv - 1}\xi_{1}^{(2\sigma + 1)\normv - 2},\\
    \xi^{\frac{1}{\sigma} - \normv}|E_{\epsilon}'(\xi)I_{P}(\xi)|
    &\leq C_{E',\epsilon}\Big(
      C_{I_{P},2,1}\|Q\|_{\normv}^{2}
      + C_{I_{P},2,2}\|Q\|_{\normv}\|Q'\|_{\normv}\xi^{-1}\\
    &\qquad+ C_{I_{P},2,3}\|Q'\|_{\normv}^{2}
      + C_{I_{P},2,4}\|Q\|_{\normv}\|Q''\|_{\normv}
      \Big)\|Q\|_{\normv}^{2\sigma - 1}
      \xi^{2\sigma\normv - 1},\\
    \xi^{\frac{1}{\sigma} - \normv}|E'(\xi)I_{P,\epsilon}(\xi)|
    &\leq C_{E'}\Big(C_{I_{P},\epsilon,1,1}\|Q\|_{\normv}^{2}
    + C_{I_{P},\epsilon,1,2}\|Q\|_{\normv}\|Q'\|_{\normv}\xi^{-1}
      + C_{I_{P},\epsilon,1,3}\|Q'\|_{\normv}^{2}\\
    &\qquad+ C_{I_{P},\epsilon,1,4}\|Q\|_{\normv}\|Q''\|_{\normv}
      + (2\sigma + 1)C_{I_{P}}\|Q\|_{\normv}\|Q_{\epsilon}\|_{\normv}\\
    &\quad\Big)\|Q\|_{\normv}^{2\sigma - 1}\xi^{2\sigma\normv - 1}.
  \end{align*}
  All terms are bounded by their value at \(\xi_{1}\), and collecting
  the terms gives us the bound.
\end{proof}

\subsubsection{Enclosures at \(\xi_{1}\)}
\label{sec:enclosures-xi-1}
In this section we describe the final procedure for computing
enclosures of \(Q\) and its derivatives at \(\xi = \xi_{1}\). The idea
is to use the results from the preceding sections to get initial
bounds for the norms. These are then iteratively improved using the
fixed-point equation. To get accurate bounds for \(I_{P}\) and its
derivatives with respect to the parameters we make use of the
integration by parts formulas from Lemmas~\ref{lemma:I_P-expansion},
\ref{lemma:I_P-derivatives-expansion-1}
and~\ref{lemma:I_P-derivatives-expansion-2}.

Let us start by considering the problem of enclosing \(Q\) and \(Q'\).
From Proposition~\ref{prop:fixed-point} we get
\(\|Q\|_{\normv} \leq \rho\) and using this together with
Lemma~\ref{lemma:norm-dQ} allows us to compute bounds for
\(\|Q'\|_{\normv}\), \(\|Q''\|_{\normv}\) and \(\|Q'''\|_{\normv}\).
The bounds for the higher derivatives will be important when improving
the bound for \(Q\) in the next step. Using that
\(|Q(\xi_{1})| \leq \|Q\|_{\normv}\xi_{1}^{\normv -
  \frac{1}{\sigma}}\) gives us initial enclosures for
\(Q(\xi_{1})\) and \(Q'(\xi_{1})\). The next step is an iterative
improvement of the bounds. For \(Q\) we use the equation
\begin{equation*}
  Q(\xi_{1}) = \gamma P(\xi_{1})
  + E(\xi_{1})I_{P}(\xi_{1}),
\end{equation*}
and for \(Q'\) we use
\begin{equation*}
  Q'(\xi_{1}) = \gamma P'(\xi_{1})
  + E'(\xi_{1})I_{P}(\xi_{1}).
\end{equation*}
Note that \(I_{E}\) is not included since \(I_{E}(\xi_{1}) = 0\). For
\(P\), \(P'\), \(E\) and \(E'\) we can easily compute tight
enclosures. For \(I_{P}\), we make use of the expansion
\begin{equation*}
  I_{P}(\xi_{1}) =
  B_{W}\left(
    \sum_{k = 1}^{3} \frac{1}{(2c)^{k}} I_{P,k}(\xi_{1})
    + \frac{1}{(2c)^{3}}\hat{I}_{P,4}(\xi_{1})
  \right)
\end{equation*}
from Lemma~\ref{lemma:I_P-expansion} (taking \(n = 3\)). Here,
\(I_{P,1}(\xi_{1})\), \(I_{P,2}(\xi_{1})\) and \(I_{P,3}(\xi_{1})\)
are explicit functions depending on \(Q(\xi_{1})\), \(Q'(\xi_{1})\)
and \(Q''(\xi_{1})\). We already have initial enclosures for \(Q\) and
\(Q'\), and an enclosure for \(Q''\) can easily be computed from the
defining ODE. To enclose \(\hat{I}_{P,4}\), we make use of
Lemma~\ref{lemma:I_P-remainder-bounds}. Combining all of these
enclosures gives us new, in general tighter, enclosures of \(Q\) and
\(Q'\). These new enclosures can then be fed back into the same
process, iteratively giving us better and better enclosures. In
practice, three iterations are enough to saturate the improvement. The
remaining error is then dominated by the error in the enclosure of
\(\hat{I}_{P,4}\), which is not improved in the iterative scheme.

The enclosures for \(Q_{\real\gamma}(\xi_{1})\),
\(Q_{\real\gamma}'(\xi_{1})\), \(Q_{\imag\gamma}(\xi_{1})\),
\(Q_{\imag\gamma}'(\xi_{1})\), \(Q_{\kappa}(\xi_{1})\),
\(Q_{\kappa}'(\xi_{1})\), \(Q_{\epsilon}(\xi_{1})\), and
\(Q_{\epsilon}'(\xi_{1})\) are computed using the same approach. In
this case we use expansions from
Lemma~\ref{lemma:I_P-derivatives-expansion-1} for
\(I_{P,\real\gamma}\), \(I_{P,\imag\gamma}\), \(I_{P,\kappa,2}\) and
\(I_{P,\epsilon,2}\) and Lemma~\ref{lemma:I_P-derivatives-expansion-2}
for \(I_{P,\kappa,1}\) and \(I_{P,\epsilon,1}\). The remainder terms
are then bounded using Lemmas~\ref{lemma:I_P-remainder-bounds-1}
and~\ref{lemma:I_P-remainder-bounds-2}.

\subsection{Verifying monotonicity}\label{sec:verify-monotonicity}
As part of counting the number of critical points of \(|Q|\), one step
is verifying that \(|Q|\) is monotone for \(\xi > \xi_{1}\). In
practice, it is easier to verify that
\(|Q|^{2} = \real(Q)^{2} + \imag(Q)^{2}\) is monotone, and we thus
work with
\begin{equation*}
  \frac{d}{d\xi}|Q|^{2} = 2\real(Q\conj{Q'}).
\end{equation*}

The first step is to split \(Q\) and \(Q'\) into a leading part we can
control, and remainder terms we can bound.
\begin{lemma}\label{lemma:Q-expansion}
  For \(\xi \geq \xi_{1}\) we have
  \begin{equation*}
    Q(\xi) = p_{Q}(\xi)P(\xi)
    + R_{Q}(\xi)\xi^{(2\sigma + 1)\normv - \frac{1}{\sigma} - 4}
  \end{equation*}
  and
  \begin{equation*}
    Q'(\xi) = p_{Q}(\xi)P'(\xi)
    + R_{Q'}(\xi)\xi^{(2\sigma + 1)\normv - \frac{1}{\sigma} - 3},
  \end{equation*}
  with coefficient
  \begin{equation*}
    p_{Q}(\xi) = \gamma + I_{E}(\xi).
  \end{equation*}
  The coefficient \(p_{Q}\) satisfies
  \begin{equation*}
    |p_{Q}(\xi) - \gamma| \leq C_{p_{Q}},
  \end{equation*}
  and the remainders are bounded as
  \begin{align*}
    |R_{Q}(\xi)| &\leq C_{R_{Q}},\\
    |R_{Q'}(\xi)| &\leq C_{R_{Q'}}.
  \end{align*}
  Here the constants are given by
  \begin{align*}
    C_{p_{Q}}
    &= C_{I_{E}}\|Q\|_{\normv}^{2\sigma + 1}\xi_{1}^{(2\sigma + 1)\normv - 2},\\
    C_{R_{Q}}
    &= C_{E}\left(
      C_{I_{P},2,1}\|Q\|_{\normv}^{2}
      + C_{I_{P},2,2}\|Q\|_{\normv}\|Q'\|_{\normv}\xi_{1}^{-1}
      + C_{I_{P},2,3}\|Q'\|_{\normv}^{2}
      + C_{I_{P},2,4}\|Q\|_{\normv}\|Q''\|_{\normv}
      \right)\|Q\|_{\normv}^{2\sigma - 1},\\
    C_{R_{Q'}}
    &= C_{E'}\left(
      C_{I_{P},2,1}\|Q\|_{\normv}^{2}
      + C_{I_{P},2,2}\|Q\|_{\normv}\|Q'\|_{\normv}\xi_{1}^{-1}
      + C_{I_{P},2,3}\|Q'\|_{\normv}^{2}
      + C_{I_{P},2,4}\|Q\|_{\normv}\|Q''\|_{\normv}
      \right)\|Q\|_{\normv}^{2\sigma - 1}.
  \end{align*}
\end{lemma}

\begin{proof}
  Defining \(p_{Q}(\xi) = \gamma + I_{E}(\xi)\), we can write \(Q\)
  and \(Q'\) as
  \begin{equation*}
    Q(\xi) = p_{Q}(\xi)P(\xi) + E(\xi)I_{P}(\xi)
    \quad\text{and}\quad
    Q'(\xi) = p_{Q}(\xi)P'(\xi) + E'(\xi)I_{P}(\xi).
  \end{equation*}
  The bound on \(|p_{Q}(\xi) - \gamma|\) follows from the bound for
  \(I_{E}\) from Lemma~\ref{lemma:I_P-I_E}. We get the bound for
  \(R_{Q}\) and \(R_{Q'}\) using Lemma~\ref{lemma:bounds-list} to
  bound \(E\) and \(E'\), and Lemma~\ref{lemma:I-P-refined} to bound
  \(I_{P}\).
\end{proof}

This allows us to similarly split \(\frac{d}{d\xi}|Q(\xi)|^{2}\) into
a leading term and a remainder.
\begin{lemma}\label{lemma:abs2-Q-derivative}
  For \(\xi \geq \xi_{1}\) we have
  \begin{equation*}
    \frac{d}{d\xi}|Q|^{2}
    = 2|p_{Q}(\xi)|^{2}\real\left(P(\xi)\conj{P'(\xi)}\right)
    + 2R_{\mathrm{mon}}(\xi)\xi^{(2\sigma + 1)\normv - \frac{2}{\sigma} - 3},
  \end{equation*}
  with
  \begin{equation*}
    |R_{\mathrm{mon}}(\xi)| \leq C_{R_{\mathrm{mon}}}
  \end{equation*}
  for
  \begin{equation*}
    C_{R_{\mathrm{mon}}}
    = (|\gamma| + C_{p_{Q}})C_{P}C_{R_{Q'}}
    + (|\gamma| + C_{p_{Q}})C_{P'}C_{R_{Q}}\xi_{1}^{-2}
    + C_{R_{Q}}C_{R_{Q'}}\xi_{1}^{(2\sigma + 1)\normv - 4}.
  \end{equation*}
\end{lemma}

\begin{proof}
  The formula
  \begin{equation*}
    \frac{d}{d\xi}|Q|^{2} = 2\real(Q\conj{Q'}),
  \end{equation*}
  together with the expansions from Lemma~\ref{lemma:Q-expansion}
  gives us
  \begin{equation*}
    \begin{split}
      \frac{d}{d\xi}|Q|^{2}
      &= 2\real\left(
        \left(
        p_{Q}(\xi)P(\xi)
        + R_{Q}(\xi)\xi^{(2\sigma + 1)\normv - \frac{1}{\sigma} - 4}
        \right)
        \left(
        \conj{p_{Q}(\xi)}\conj{P'(\xi)}
        + \conj{R_{Q'}(\xi)}\xi^{(2\sigma + 1)\normv - \frac{1}{\sigma} - 3}
        \right)
        \right)\\
      &= 2|p_{Q}(\xi)|^{2}\real\left(P(\xi)\conj{P'(\xi)}\right)
        + 2\real\Big(
        p_{Q}(\xi)P(\xi)\conj{R_{Q'}(\xi)}
      + \conj{p_{Q}(\xi)}\conj{P'(\xi)}R_{Q}(\xi)\xi^{-1}\\
      &\qquad+ R_{Q}(\xi)\conj{R_{Q'}(\xi)}\xi^{(2\sigma + 1)\normv - \frac{1}{\sigma} - 4}
        \Big)\xi^{(2\sigma + 1)\normv - \frac{1}{\sigma} - 3}.
    \end{split}
  \end{equation*}
  To bound the remainder term we note that
  Lemma~\ref{lemma:bounds-list} combined with
  Lemma~\ref{lemma:Q-expansion} gives us
  \begin{multline*}
    \left|\real\left(
        p_{Q}(\xi)P(\xi)\conj{R_{Q'}(\xi)}
        + \conj{p_{Q}(\xi)}\conj{P'(\xi)}R_{Q}(\xi)\xi^{-1}
        + R_{Q}(\xi)\conj{R_{Q'}(\xi)}\xi^{(2\sigma + 1)\normv - \frac{1}{\sigma} - 4}
      \right)\right|\\
    \leq |p_{Q}(\xi)|C_{P}C_{R_{Q'}}\xi^{-\frac{1}{\sigma}}
    + |p_{Q}(\xi)|C_{P'}C_{R_{Q}}\xi^{-\frac{1}{\sigma} - 2}
    + C_{R_{Q}}C_{R_{Q'}}\xi^{(2\sigma + 1)\normv - \frac{1}{\sigma} - 4}\\
    \leq \left(
      (|\gamma| + C_{p_{Q}})C_{P}C_{R_{Q'}}
      + (|\gamma| + C_{p_{Q}})C_{P'}C_{R_{Q}}\xi^{-2}
      + C_{R_{Q}}C_{R_{Q'}}\xi^{(2\sigma + 1)\normv - 4}
    \right)\xi^{-\frac{1}{\sigma}},
  \end{multline*}
  from which the result follows.
\end{proof}

As a consequence, to prove that \(|Q|^{2}\) is monotone for
\(\xi > \xi_{1}\), it suffices to verify that
\begin{equation*}
  \left||p_{Q}(\xi)|^{2}\real\left(P(\xi)\conj{P'(\xi)}\right)\right|
  > C_{R_{\mathrm{mon}}}\xi^{(2\sigma + 1)\normv - \frac{2}{\sigma} - 3}.
\end{equation*}
Or equivalently, taking into account that \(P(\xi)\conj{P'(\xi)}\)
decays like \(\xi^{-\frac{2}{\sigma} - 1}\),
\begin{equation}\label{eq:monotone-condition}
  \left|\real\left(P(\xi)\conj{P'(\xi)}\xi^{\frac{2}{\sigma} + 1}\right)\right|
  > \frac{C_{R_{\mathrm{mon}}}}{|p_{Q}(\xi)|^{2}}\xi^{(2\sigma + 1)\normv - 2}.
\end{equation}

To verify that this condition holds for large values of \(\xi\) we
need to get an enclosure for right-hand side. For this we have the
following lemma.

\begin{lemma}\label{lemma:real-P-P-dxi-expansion}
  Let \(\xi_{1} > 1\). We have
  \begin{equation*}
    P(\xi)\conj{P'(\xi)}\xi^{\frac{2}{\sigma} + 1}
    = -2a|c^{-a}|^{2}(S_{1}(\xi)S_{2}(\xi) + R_{P\conj{P'}}(\xi)),
  \end{equation*}
  where \(S_{1}\) and \(S_{2}\) are the sums
  \begin{align*}
    S_{1} &= \sum_{k = 0}^{n - 1} \frac{(a)_{k}(a - b + 1)_{k}}{k!(-c)^{k}}\xi^{-2k},\\
    S_{2} &= \sum_{k = 0}^{n - 1} \conj{\left(\frac{(a + 1)_{k}(a - b + 1)_{k}}{k!(-c)^{k}}\right)}\xi^{-2k}
  \end{align*}
  and the remainder is bounded as
  \begin{equation*}
    |R_{P\conj{P'}}(\xi)| \leq C_{R_{P\conj{P'}}}
  \end{equation*}
  with
  \begin{align*}
    C_{R_{P\conj{P'}}}
    &= \left(\sum_{k = 0}^{n - 1} \left|\frac{(a)_{k}(a - b + 1)_{k}}{k!(-c)^{k}}\right|\xi_{1}^{-2k}\right)C_{R_{U}}(a + 1, b + 1, n, c\xi_{1}^{2})|c^{-n}|\xi_{1}^{-2n}\\
      &\qquad+ \left(\sum_{k = 0}^{n - 1} \left|\frac{(a + 1)_{k}(a - b + 1)_{k}}{k!(-c)^{k}}\right|\xi_{1}^{-2k}\right)C_{R_{U}}(a, b, n, c\xi_{1}^{2})|c^{-n}|\xi_{1}^{-2n}\\
    &\qquad+ C_{R_{U}}(a, b, n, c\xi_{1}^{2})C_{R_{U}}(a + 1, b + 1, n, c\xi_{1}^{2})|c^{-n}|^{2}\xi_{1}^{-4n}.
  \end{align*}
\end{lemma}

\begin{proof}
  Writing \(P\) and \(P'\) in terms of \(U\) we have
  \begin{align*}
    P(\xi)
    &= U(a, b, c\xi^{2}),\\
    P'(\xi)
    &= 2c\xi U'(a, b, c\xi^{2}) = -2ac\xi U(a + 1, b + 1, c\xi^{2}).
  \end{align*}
  The expansion from Lemma~\ref{lemma:U} gives us
  \begin{equation*}
    P(\xi) = c^{-a}\left(\sum_{k = 0}^{n - 1} \frac{(a)_{k}(a - b + 1)_{k}}{k!(-c)^{k}}\xi^{-2k} + R_{U}(a, b, n, c\xi^{2})c^{-n}\xi^{-2n}\right)\xi^{-2a}
  \end{equation*}
  and
  \begin{equation*}
    P'(\xi) = -2ac^{-a}\left(\sum_{k = 0}^{n - 1} \frac{(a + 1)_{k}(a - b + 1)_{k}}{k!(-c)^{k}}\xi^{-2k} + R_{U}(a + 1, b + 1, n, c\xi^{2})c^{-n}\xi^{-2n}\right)\xi^{-2a - 1}.
  \end{equation*}
  Expanding the product \(P(\xi)\conj{P'(\xi)}\), we have
  \begin{equation*}
    \begin{split}
      P(\xi)\conj{P'(\xi)}
      &= -2a|c^{-a}|^{2}\Bigg(\\
      &\qquad\left(\sum_{k = 0}^{n - 1} \frac{(a)_{k}(a - b + 1)_{k}}{k!(-c)^{k}}\xi^{-2k}\right)
        \conj{\left(\sum_{k = 0}^{n - 1} \frac{(a + 1)_{k}(a - b + 1)_{k}}{k!(-c)^{k}}\xi^{-2k}\right)}\\
      &\qquad+ \left(\sum_{k = 0}^{n - 1} \frac{(a)_{k}(a - b + 1)_{k}}{k!(-c)^{k}}\xi^{-2k}\right)\conj{R_{U}(a + 1, b + 1, n, c\xi^{2})}\conj{c}^{-n}\xi^{-2n}\\
      &\qquad+ \conj{\left(\sum_{k = 0}^{n - 1} \frac{(a + 1)_{k}(a - b + 1)_{k}}{k!(-c)^{k}}\xi^{-2k}\right)}R_{U}(a, b, n, c\xi^{2})c^{-n}\xi^{-2n}\\
      &\qquad+ R_{U}(a, b, n, c\xi^{2})\conj{R_{U}(a + 1, b + 1, n, c\xi^{2})}|c^{-n}|^{2}\xi^{-4n}
        \Bigg)\xi^{-\frac{2}{\sigma} - 1}.
    \end{split}
  \end{equation*}
  The last three terms in the inner parenthesis is exactly
  \(R_{P\conj{P'}}\), which we can bound as in Lemma~\ref{lemma:U}.
\end{proof}

From this lemma we can derive a sufficient condition for \(|Q|^{2}\)
to be monotone.

\begin{lemma}\label{lemma:monotonicity-infinity}
  If \(|\gamma| > C_{p_{Q}}\) and
  \begin{equation}\label{eq:monotonicity-infinity}
    \left|\real\left(2\conj{a}|c^{-a}|^{2}(S_{1}(\xi)S_{2}(\xi) + R_{P\conj{P'}}(\xi))\right)\right|
    > \frac{C_{R_{\mathrm{mon}}}}{|p_{Q}(\xi)|^{2}}\xi_{1}^{(2\sigma + 1)\normv - 2},
  \end{equation}
  then \(|Q|^{2}\) is monotone for \(\xi \geq \xi_{1}\).
\end{lemma}

\begin{proof}
  This is a direct consequence of the
  condition~\eqref{eq:monotone-condition} and
  Lemma~\ref{lemma:real-P-P-dxi-expansion}.
\end{proof}

To verify this condition, we need to compute an enclosure of the
left-hand side of~\eqref{eq:monotonicity-infinity}. For this, we note
that \(S_{1}(\xi)S_{2}(\xi)\) is a polynomial in \(\xi^{-1}\), we can
hence directly compute an enclosure by evaluating this polynomial on
the interval \([0, \xi_{1}^{-1}]\) using interval arithmetic.

\section{Solution at zero}
\label{sec:solution-zero}
In this section we study solutions to the initial value problem
\begin{equation*}
  \begin{split}
    (1 - i\epsilon)\left(Q'' + \frac{d - 1}{\xi}Q'\right) + i\kappa\xi Q'
    + i \frac{\kappa}{\sigma}Q - \omega Q + (1 + i\delta)|Q|^{2\sigma}Q &= 0,\\
    Q(0) &= \mu,\\
    Q'(0) &= 0.
  \end{split}
\end{equation*}
To compute \(G\)~\eqref{eq:G} and its partial derivatives we need to
compute \(Q(\xi_1)\) and \(Q'(\xi_1)\), as well as their derivatives
with respect to \(\mu\), \(\kappa\) and \(\epsilon\), at a given
\(\xi_1 > 0\). For counting the number of critical points of \(|Q|\)
we need enclosures of the curves \(Q\), \(Q'\) and \(Q''\) on an
interval \([0, \xi_{2}]\).

For the computations we split \(Q\) into its real and imaginary parts,
\(Q(\xi) = a(\xi) + ib(\xi)\). This gives us the system of equations
\begin{align*}
  a'' + \epsilon b'' + \frac{d - 1}{\xi}(a' + \epsilon b') - \kappa \xi b' - \frac{\kappa}{\sigma}b - \omega a + (a^{2} + b^{2})^{\sigma}a - \delta(a^{2} + b^{2})^{\sigma}b &= 0,\\
  b'' - \epsilon a'' + \frac{d - 1}{\xi}(b' - \epsilon a') + \kappa \xi a' + \frac{\kappa}{\sigma}a - \omega b + (a^{2} + b^{2})^{\sigma}b + \delta(a^{2} + b^{2})^{\sigma}a &= 0.
\end{align*}
with initial conditions \(a(0) = \mu\), \(b(0) = 0\) and
\(a'(0) = b'(0) = 0\).

The system is integrated using a rigorous forward numerical integrator
implemented by the CAPD library~\cite{Kapela2021}. The library
supports computations of very tight enclosures of the solution, and
its partial derivatives with respect to the parameters.

For \(d \neq 1\) the numerical integrator cannot deal with the
removable singularity at \(\xi = 0\). To handle this we must perform a
first step using a Taylor expansion at \(\xi = 0\) to compute the
solution up to \(\xi = \xi_0\) for some small \(\xi_0 > 0\). This is
described in detail below. After that the numerical integrator is used
on the interval \([\xi_0, \xi_1]\). Note that for \(d = 1\) there is
no removable singularity and the numerical integrator can be used on
the full interval \([0, \xi_1]\).

\subsection{Handling the removable singularity at zero}
\label{sec:solution-zero-removable-singularity}
For \(d \neq 1\), the equation has a removable singularity at
\(\xi = 0\). Taylor expanding \(a\) and \(b\) at \(\xi = 0\) as
\(a = \sum_{n = 0}^{\infty} a_{n}\xi^{n}\) and
\(b = \sum_{n = 0}^{\infty} b_{n}\xi^{n}\) gives us
\begin{equation*}
  a_{n + 2} = \frac{F_{1,n} - \epsilon F_{2,n}}{(n + 2)(n + d)(1 + \epsilon^{2})},\quad
  b_{n + 2} = \frac{\epsilon F_{1,n} + F_{2,n}}{(n + 2)(n + d)(1 + \epsilon^{2})}
\end{equation*}
with
\begin{align*}
  F_{1,n} &= \kappa n b_{n} + \frac{\kappa}{\sigma}b_{n} + \omega a_{n} - u_{1,n} + \delta u_{2,n},\\
  F_{2,n} &= -\kappa n a_{n} - \frac{\kappa}{\sigma}a_{n} + \omega b_{n} - u_{2,n} - \delta u_{1,n}.
\end{align*}
Here
\begin{equation*}
  u_{1,n} = \left((a^{2} + b^{2})^{\sigma}a\right)_{n} \quad\text{and}\quad
  u_{2,n} = \left((a^{2} + b^{2})^{\sigma}b\right)_{n}.
\end{equation*}

Using the above recursion with \(a_{0} = \mu\), \(b_{0} = 0\) and
\(a_{1} = b_{1} = 0\) it is straightforward to compute Taylor
expansions of arbitrarily high order. What remains is to bound the
remainder term. We will show that for \(n > N\) we have
\begin{equation*}
  |a_{n}|, |b_{n}| \leq r^{n}
\end{equation*}
for some \(r\). This allows us to bound the remainder term for \(0 <
\xi < \frac{1}{r}\) as
\begin{equation*}
  \left|\sum_{n = N + 1}^{\infty} a_{n}\xi^{n}\right|,
  \left|\sum_{n = N + 1}^{\infty} b_{n}\xi^{n}\right|
  \leq \sum_{n = N + 1}^{\infty}(r\xi)^{n}
  = \frac{(r\xi)^{N + 1}}{1 - r\xi}.
\end{equation*}
For the first and second derivatives, the bounds are
\begin{equation*}
  \left|\frac{d}{d\xi}\left(\sum_{n = N + 1}^{\infty} a_{n}\xi^{n}\right)\right|
  = \left|\sum_{n = N + 1}^{\infty} na_{n}\xi^{n - 1}\right|
  \leq \frac{1}{\xi}\sum_{n = N + 1}^{\infty}n(r\xi)^{n}
  = \frac{r(r\xi)^{N}(N + 1 - Nr\xi)}{(1 - r\xi)^{2}}
\end{equation*}
and
\begin{multline*}
  \left|\frac{d^{2}}{d\xi^{2}}\left(\sum_{n = N + 1}^{\infty} a_{n}\xi^{n}\right)\right|
  = \left|\sum_{n = N + 1}^{\infty} n(n - 1)a_{n}\xi^{n - 2}\right|\\
  \leq \frac{1}{\xi^{2}}\sum_{n = N + 1}^{\infty}n(n - 1)(r\xi)^{n}
  = \frac{r^{2}(r\xi)^{N-1}(N + N^2 + (2 - 2N^2)r\xi - (N - N^2)(r\xi)^2)}{(1 - r\xi)^{3}},
\end{multline*}
with the same bound for the \(b_n\) case.

To start with we have the following lemma, which is valid for
\(\sigma = 1\).
\begin{lemma}\label{lemma:tail-bound}
  Let \(\sigma = 1\). Let \(M\), \(N\), \(C\) and \(r\) be such that
  \(N\) is even, \(3M < N\),
  \begin{equation*}
    |a_n|, |b_n| \leq C r^n \quad\text{for}\quad n < M
  \end{equation*}
  and
  \begin{equation*}
    |a_n|, |b_n| \leq r^n \quad\text{for}\quad M \leq n \leq N.
  \end{equation*}
  If
  \begin{multline}
    \label{eq:tail-bound-inequality}
    \frac{(1 + |\epsilon|)}{(1 + \epsilon^{2})}\Bigg(
      \frac{|\kappa|}{N + d}
      + \frac{|\omega|}{(N + 2)(N + d)}\\
      + 2(1 + |\delta|)\left(
        \frac{1}{8} + \frac{1}{2N}
        + \frac{3mC(1 + 3/N)}{2(N + d)}
        + \frac{3m^{2}C^{2}}{(N + 2)(N + d)}
      \right)
    \Bigg) \leq r^2,
  \end{multline}
  where \(m = \lceil M / 2 \rceil\), then
  \begin{equation*}
    |a_{n}|, |b_{n}| \leq r^{n} \quad\text{for}\quad n > N.
  \end{equation*}
\end{lemma}

\begin{proof}
  It is enough to show that
  \begin{equation*}
    |a_{N + 1}|, |b_{N + 1}| \leq r^{N + 1},
    |a_{N + 2}|, |b_{N + 2}| \leq r^{N + 2}
  \end{equation*}
  and use that the left-hand side in \eqref{eq:tail-bound-inequality}
  is decreasing in \(N\). The result then follows for \(n > N + 2\) by
  induction.

  Since \(a_{n}\) and \(b_{n}\) are zero for odd \(n\) and \(N\) is
  even it follows that \(a_{N + 1}\) and \(b_{N + 1}\) are zero. We
  therefore only have to prove that
  \(|a_{N + 2}|, |b_{N + 2}| \leq r^{N + 2}\). We give the proof for
  \(a_{N + 2}\), the bound for \(b_{N + 2}\) follows in exactly the
  same way.

  With \(\sigma = 1\) we get \(u_{1,N} = (a^{3})_N + (ab^{2})_N\) and
  \(u_{2,N} = (a^{2}b)_N + (b^{3})_N\). We have
  \begin{equation*}
    (a^{3})_{N} = \sum_{i + j + k = N}a_{i}a_{j}a_{k}.
  \end{equation*}
  However, since \(a_{n} = 0\) for odd \(n\) and \(N\) is even we can
  limit the sum to
  \begin{equation*}
    (a^{3})_{N} = \sum_{i + j + k = N/2}a_{2i}a_{2j}a_{2k}.
  \end{equation*}
  Since the bound for \(a_{n}\) depends on if \(n\) is less than \(M\)
  or not, let us group the terms in the sum by how many of the indices
  \(i, j, k\) are less than \(m = \lceil M / 2 \rceil\):
  \begin{enumerate}
  \item The number of terms with no indices less than \(m\) is
    \begin{equation*}
      T_{0} = \sum_{i = 0}^{N/2 - 3m}\sum_{j = 0}^{N/2 - 3m - i} 1 = \frac{(N/2 - 3m + 1)(N/2 - 3m + 2)}{2}.
    \end{equation*}
    These are bounded by \(r^{N}\).
  \item For the terms with exactly one index less than \(m\) there are
    \(3\) choices for which index is the small one. If we let \(i\) be
    the small index, then there are \(N / 2 - 2m - i + 1\) choices for
    \(j\) and \(k\). Summing over all choices for \(i\) and
    multiplying by \(3\) we get
    \begin{equation*}
      T_{1} = 3\sum_{i = 0}^{m - 1} (\frac{N}{2} - 2m - i + 1) = 3\frac{m}{2}\left(N - 5m + 3\right).
    \end{equation*}
    These are bounded by \(Cr^{N}\).
  \item The number of terms with exactly two indices less than \(m\)
    is \(T_{2} = 3m^{2}\). These are bounded by \(C^{2}r^{N}\).
  \end{enumerate}
  Note that since \(3M < N\) there are no terms with three indices
  less than \(m\). In total this gives us the bound
  \begin{equation*}
    |(a^{3})_{N}| \leq (T_{0} + T_{1}C + T_{2}C^{2})r^{N}.
  \end{equation*}

  Since the bounds for \(|a_n|\) and \(|b_n|\) are the same we get the
  same bound for \(|(ab^2)_N|\), \(|(a^2b)_N|\) and \(|(b^3)_N|\). In
  particular, this gives us
  \begin{equation*}
    |u_{1,N}|, |u_{2,N}| \leq 2(T_{0} + T_{1}C + T_{2}C^{2})r^{N}.
  \end{equation*}

  Combining the above bounds for \(|u_{1,N}|\) and \(|u_{2,N}|\) with
  the bounds for \(|a_n|\) and \(|b_n|\) we get
  \begin{equation*}
    |F_{1,N}| \leq (|\kappa|(N + 1) + |\omega|)r^{N} + 2(1 + |\delta|)(T_{0} + T_{1}C + T_{2}C^{2})r^{N},
  \end{equation*}
  and the same for \(|F_{2,N}|\). This gives us
  \begin{equation*}
    \begin{split}
      |a_{N + 2}|
      &\leq \frac{(1 + |\epsilon|)(|\kappa|(N + 1) + |\omega| + 2(1 + |\delta|)(T_{0} + T_{1}C + T_{2}C^{2}))}{(N + 2)(N + d)(1 + \epsilon^{2})}r^{N}\\
      &= \frac{(1 + |\epsilon|)}{(1 + \epsilon^{2})}\left(
        \frac{|\kappa|}{N + d}
        + \frac{|\omega|}{(N + 2)(N + d)}
        + 2(1 + |\delta|)\frac{(T_{0} + T_{1}C + T_{2}C^{2})}{(N + 2)(N + d)}
        \right)r^{N}.
    \end{split}
  \end{equation*}
  We have
  \begin{equation*}
    \begin{split}
      \frac{T_{0}}{(N + 2)(N + d)}
      &= \frac{(N/2 - 3m + 1)(N/2 - 3m + 2)}{2(N + 2)(N + d)}\\
      &= \frac{1}{8}\frac{(N - 6m + 2)(N - 6m + 4)}{(N + 2)(N + d)}\\
      &\leq \frac{1}{8}\frac{N + 4}{N + d}
        \leq \frac{1}{8}\left(1 + \frac{4}{N}\right)
        = \frac{1}{8} + \frac{1}{2N},
    \end{split}
  \end{equation*}
  as well as
  \begin{equation*}
    \frac{T_{1}}{(N + 2)(N + d)}
    = 3\frac{m}{2}\frac{N - 5m + 3}{(N + 2)(N + d)}
    \leq 3\frac{m}{2}\frac{N + 3}{(N + 2)(N + d)}
    \leq \frac{3m(1 + 3/N)}{2(N + d)}.
  \end{equation*}
  This gives us
  \begin{multline*}
    |a_{N + 2}|
    \leq \frac{(1 + |\epsilon|)}{(1 + \epsilon^{2})}\Bigg(
    \frac{|\kappa|}{N + d}
      + \frac{|\omega|}{(N + 2)(N + d)}\\
      + 2(1 + |\delta|)\left(
        \frac{1}{8} + \frac{1}{2N}
        + \frac{3mC(1 + 3/N)}{2(N + d)}
        + \frac{3m^{2}C^{2}}{(N + 2)(N + d)}
      \right)
    \Bigg)r^{N}.
  \end{multline*}

  By~\eqref{eq:tail-bound-inequality}, the factor in front of
  \(r^{N}\) is bounded by \(r^{2}\), giving us our required bound
  \(|a_{N + 2}| \leq r^{N + 2}\).
\end{proof}

\subsection{Derivatives with respect to \(\mu\), \(\kappa\) and \(\epsilon\)}
Next, we handle the derivatives with respect to \(\mu\), \(\kappa\)
and \(\epsilon\). We let \(a_{\mu}\) and \(b_{\mu}\) denote the
derivatives of \(a\) and \(b\) with respect to the initial value
\(\mu\), similarly with the derivatives with respect to \(\kappa\) and
\(\epsilon\).

Since the CAPD library automatically supports computing derivatives
with respect to the parameters we do not directly need the ODEs
satisfied by these functions. However, for the removable singularity
at zero we do need Taylor expansions of these derivatives. Since the
series expansions for \(a\) and \(b\) are absolutely convergent we can
directly differentiate the expressions for \(a_{n}\) and \(b_{n}\) to
get recurrence relations for the coefficients of the derivatives.

Differentiating the recurrence relations for \(a_{n}\) and \(b_{n}\)
with respect to \(\mu\), \(\kappa\) and \(\epsilon\) respectively, we get
\begin{align*}
  a_{\mu,n + 2} &= \frac{F_{\mu,1,n} - \epsilon F_{\mu,2,n}}{(n + 2)(n + d)(1 + \epsilon^{2})},\quad
  b_{\mu,n + 2} = \frac{\epsilon F_{\mu,1,n} + F_{\mu,2,n}}{(n + 2)(n + d)(1 + \epsilon^{2})},\\
  a_{\kappa,n + 2} &= \frac{F_{\kappa,1,n} - \epsilon F_{\kappa,2,n}}{(n + 2)(n + d)(1 + \epsilon^{2})},\quad
  b_{\kappa,n + 2} = \frac{\epsilon F_{\kappa,1,n} + F_{\kappa,2,n}}{(n + 2)(n + d)(1 + \epsilon^{2})},\\
  a_{\epsilon,n + 2} &= \frac{F_{\epsilon,1,n} - \epsilon F_{\epsilon,2,n}}{(n + 2)(n + d)(1 + \epsilon^{2})},\quad
  b_{\epsilon,n + 2} = \frac{\epsilon F_{\epsilon,1,n} + F_{\epsilon,2,n}}{(n + 2)(n + d)(1 + \epsilon^{2})},
\end{align*}
with
\begin{align*}
  F_{\mu,1,n} &= \kappa n b_{\mu,n} + \frac{\kappa}{\sigma}b_{\mu,n} + \omega a_{\mu,n} - u_{\mu,1,n} + \delta u_{\mu,2,n},\\
  F_{\mu,2,n} &= -\kappa n a_{\mu,n} - \frac{\kappa}{\sigma}a_{\mu,n} + \omega b_{\mu,n} - u_{\mu,2,n} - \delta u_{\mu,1,n},\\
  F_{\kappa,1,n} &= \kappa n b_{\kappa,n} + n b_{n} + \frac{\kappa}{\sigma}b_{\kappa,n} + \frac{1}{\sigma}b_{n} + \omega a_{\kappa,n} - u_{\kappa,1,n} + \delta u_{\kappa,2,n},\\
  F_{\kappa,2,n} &= -\kappa n a_{\kappa,n} - n a_{n} - \frac{\kappa}{\sigma}a_{\kappa,n} - \frac{1}{\sigma}a_{n} + \omega b_{\kappa,n} - u_{\kappa,2,n} - \delta u_{\kappa,1,n},\\
  F_{\epsilon,1,n} &= -(n + 2)(n + d)b_{n + 2} + \kappa n b_{\epsilon,n} + \frac{\kappa}{\sigma}b_{\epsilon,n} + \omega a_{\epsilon,n} - u_{\epsilon,1,n} + \delta u_{\epsilon,2,n},\\
  F_{\epsilon,2,n} &= (n + 2)(n + d)a_{n + 2} - \kappa n a_{\epsilon,n} - \frac{\kappa}{\sigma}a_{\epsilon,n} + \omega b_{\epsilon,n} - u_{\epsilon,2,n} - \delta u_{\epsilon,1,n}.
\end{align*}
Here
\begin{align*}
  u_{\mu,1,n}
  &= \left(
    (a^{2} + b^{2})^{\sigma}a_{\mu}
    + 2\sigma(a^{2} + b^{2})^{\sigma - 1}a^{2}a_{\mu}
    + 2\sigma(a^{2} + b^{2})^{\sigma - 1}abb_{\mu}
    \right)_{n},\\
  u_{\mu,2,n}
  &= \left(
    (a^{2} + b^{2})^{\sigma}b_{\mu}
    + 2\sigma(a^{2} + b^{2})^{\sigma - 1}aba_{\mu}
    + 2\sigma(a^{2} + b^{2})^{\sigma - 1}b^{2}b_{\mu}
    \right)_{n},
\end{align*}
and similarly for \(\kappa\) and \(\epsilon\). For the initial
conditions we get \(a_{\mu,0} = 1\) and
\(a_{\mu,1} = b_{\mu,0} = b_{\mu,1} = 0\) for the derivatives with
respect to \(\mu\). For the derivatives with respect to \(\kappa\) and
\(\epsilon\) the initial conditions are all zero.

To bound the remainder term we follow the same approach as above,
showing that for \(n > N\) we have
\(|a_{\mu,n}|, |b_{\mu,n}| \leq r^{n}\) for some \(r\), and similarly
for the derivatives with respect to \(\kappa\) and \(\epsilon\).

For \(a_{\mu,n}\) and \(b_{\mu,n}\), the only difference compared to
\(a_{n}\) and \(b_{n}\) is in the form of the non-linearity. For
\(a_{\kappa,n}\), \(b_{\kappa,n}\), \(a_{\epsilon,n}\) and
\(b_{\epsilon,n}\) we have slightly modified recurrence relations. We
therefore get slightly modified versions of
Lemma~\ref{lemma:tail-bound} in all three cases.

\begin{lemma}\label{lemma:tail-bound-dmu}
  Let \(\sigma = 1\). Let \(M\), \(N\), \(C\), \(r\) satisfy the
  conditions of Lemma~\ref{lemma:tail-bound}, as well as
  \begin{equation*}
    |a_{\mu,n}|, |b_{\mu,n}| \leq C r^{n} \quad\text{for}\quad n < M
  \end{equation*}
  and
  \begin{equation*}
    |a_{\mu,n}|, |b_{\mu,n}| \leq r^{n} \quad\text{for}\quad M \leq n \leq N.
  \end{equation*}
  If the inequality~\eqref{eq:tail-bound-inequality} holds and we also
  have
  \begin{multline}
    \label{eq:tail-bound-inequality-dmu}
    \frac{(1 + |\epsilon|)}{(1 + \epsilon^{2})}\Bigg(
    \frac{|\kappa|}{N + d}
    + \frac{|\omega|}{(N + 2)(N + d)}\\
    + 6(1 + |\delta|)\left(
      \frac{1}{8} + \frac{1}{2N}
      + \frac{3mC(1 + 3/N)}{2(N + d)}
      + \frac{3m^{2}C^{2}}{(N + 2)(N + d)}
    \right)
    \Bigg) \leq r^2,
  \end{multline}
  then
  \begin{equation*}
    |a_{\mu,n}|, |b_{\mu,n}| \leq r^{n} \quad\text{for}\quad n > N.
  \end{equation*}
\end{lemma}

\begin{proof}
  The proof follows the same strategy as Lemma~\ref{lemma:tail-bound}.
  The only required adjustment is in the bounds of \(u_{\mu,1,N}\) and
  \(u_{\mu,2,N}\). In this case we have (for \(\sigma = 1\)),
  \begin{align*}
    u_{\mu,1,N} &= (a^{2}a_{\mu} + b^{2}a_{\mu} + 2a^{2}a_{\mu} + 2abb_{\mu})_{N},\\
    u_{\mu,2,N} &= (a^{2}b_{\mu} + b^{2}b_{\mu} + 2aba_{\mu} + 2b^{2}b_{\mu})_{N}.
  \end{align*}
  The triple products can be bounded in exactly the same way as
  before, giving us
  \begin{equation*}
    |u_{\mu,1,N}|, |u_{\mu,2,N}| \leq 6(T_{0} + T_{1}C + T_{2}C^{2})r^{N}.
  \end{equation*}
  The rest of the proof is the same as in
  Lemma~\ref{lemma:tail-bound}. We get the bound
  \begin{equation*}
    |a_{\mu,N + 2}|
    \leq \frac{(1 + |\epsilon|)}{(1 + \epsilon^{2})}\left(
      \frac{|\kappa|}{N + d}
      + \frac{|\omega|}{(N + 2)(N + d)}
      + 6(1 + |\delta|)\frac{(T_{0} + T_{1}C + T_{2}C^{2})}{(N + 2)(N + d)}
    \right)r^{N},
  \end{equation*}
  This can then be shown to be bounded by the left-hand side of
  inequality~\eqref{eq:tail-bound-inequality-dmu} times \(r^{N}\),
  which gives us the required bound.
\end{proof}

The versions for the derivatives with respect to \(\kappa\) and
\(\epsilon\) follow in a similar way, and we leave out the proofs.

\begin{lemma}\label{lemma:tail-bound-dkappa}
  Let \(\sigma = 1\). Let \(M\), \(N\), \(C\), \(r\) satisfy the
  conditions of Lemma~\ref{lemma:tail-bound}, as well as
  \begin{equation*}
    |a_{\kappa,n}|, |b_{\kappa,n}| \leq C r^{n} \quad\text{for}\quad n < M
  \end{equation*}
  and
  \begin{equation*}
    |a_{\kappa,n}|, |b_{\kappa,n}| \leq r^{n} \quad\text{for}\quad M \leq n \leq N.
  \end{equation*}
  If the inequality~\eqref{eq:tail-bound-inequality} holds and we also
  have
  \begin{multline}
    \label{eq:tail-bound-inequality-dkappa}
    \frac{(1 + |\epsilon|)}{(1 + \epsilon^{2})}\Bigg(
      \frac{|\kappa| + 1}{N + d}
      + \frac{|\omega|}{(N + 2)(N + d)}\\
      + 6(1 + |\delta|)\left(
        \frac{1}{8} + \frac{1}{2N}
        + \frac{3mC(1 + 3/N)}{2(N + d)}
        + \frac{3m^{2}C^{2}}{(N + 2)(N + d)}
      \right)
    \Bigg) \leq r^2,
  \end{multline}
  then
  \begin{equation*}
    |a_{\kappa,n}|, |b_{\kappa,n}| \leq r^{n} \quad\text{for}\quad n > N.
  \end{equation*}
\end{lemma}

\begin{lemma}\label{lemma:tail-bound-depsilon}
  Let \(\sigma = 1\). Let \(M\), \(N\), \(C\), \(r\) satisfy the
  conditions of Lemma~\ref{lemma:tail-bound}, as well as
  \begin{equation*}
    |a_{\epsilon,n}|, |b_{\epsilon,n}| \leq C r^{n} \quad\text{for}\quad n < M
  \end{equation*}
  and
  \begin{equation*}
    |a_{\epsilon,n}|, |b_{\epsilon,n}| \leq r^{n} \quad\text{for}\quad M \leq n \leq N.
  \end{equation*}
  If the inequality~\eqref{eq:tail-bound-inequality} holds and we also
  have
  \begin{multline}
    \label{eq:tail-bound-inequality-depsilon}
    \frac{(1 + |\epsilon|)}{(1 + \epsilon^{2})}\Bigg(
    1 +
    \frac{|\kappa|}{N + d}
    + \frac{|\omega|}{(N + 2)(N + d)}\\
    + 6(1 + |\delta|)\left(
      \frac{1}{8} + \frac{1}{2N}
      + \frac{3mC(1 + 3/N)}{2(N + d)}
      + \frac{3m^{2}C^{2}}{(N + 2)(N + d)}
    \right)
    \Bigg) \leq r^2,
  \end{multline}
  then
  \begin{equation*}
    |a_{\epsilon,n}|, |b_{\epsilon,n}| \leq r^{n} \quad\text{for}\quad n > N.
  \end{equation*}
\end{lemma}

\section{Pointwise results for CGL in Case II}
\label{sec:pointwise-results-cgl}
The results presented for the CGL equation in \textbf{Case II} in
Theorem~\ref{thm:existence-cgl-case-2} are to a large extent limited
by computational power. Instead of attempting to prove the existence
of a continuous curve of solutions, we can prove the existence of
discrete points along the numerical approximations of the curves. This
significantly reduces the computational cost, and allows us to cover a
much larger part of the curves. Of course, this comes with the drawback that we
do not actually prove that these discrete points do in fact lie on the
same curve.

The starting point for the pointwise results are the numerical
approximations of the branches given in
Figure~\ref{fig:branches-d3-numerical}. These numerical approximations
consist of a large number of discrete points that (presumably) form
the curves of solutions. The points are in general not uniformly
spread out, but cluster around the turning points, where the
continuation is more challenging. The first branch consists of 2433
discrete points, and the remaining four consist of 1179, 11094, 1187,
and 5070 points, respectively.

\begin{figure}
  \centering
  \begin{subfigure}[t]{0.45\textwidth}
    \includegraphics[width=\textwidth]{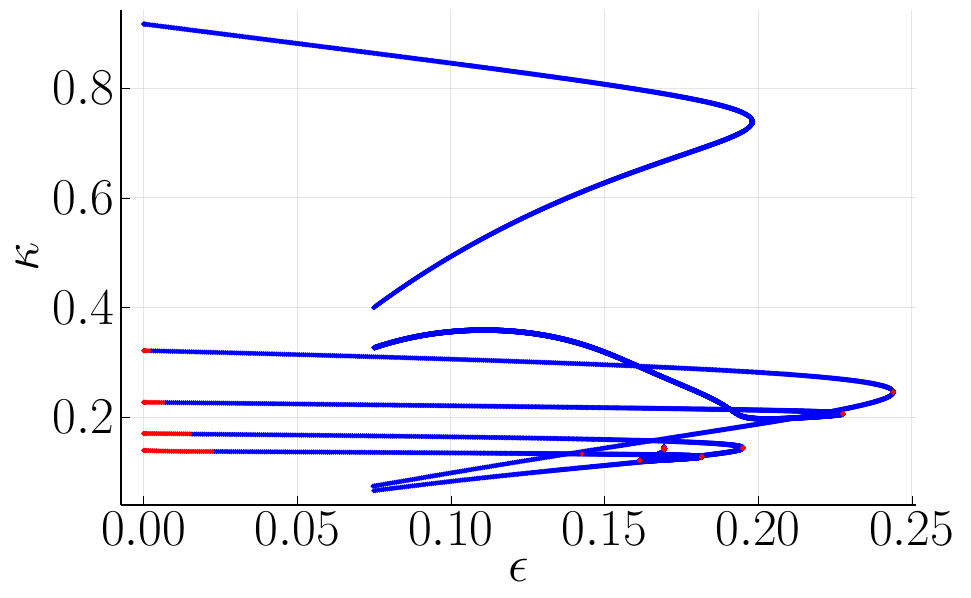}
    \caption{Fixing \(\epsilon\).}
    \label{fig:CGL-pointwise-d=3-epsilon}
  \end{subfigure}
  \hspace{0.05\textwidth}
  \begin{subfigure}[t]{0.45\textwidth}
    \includegraphics[width=\textwidth]{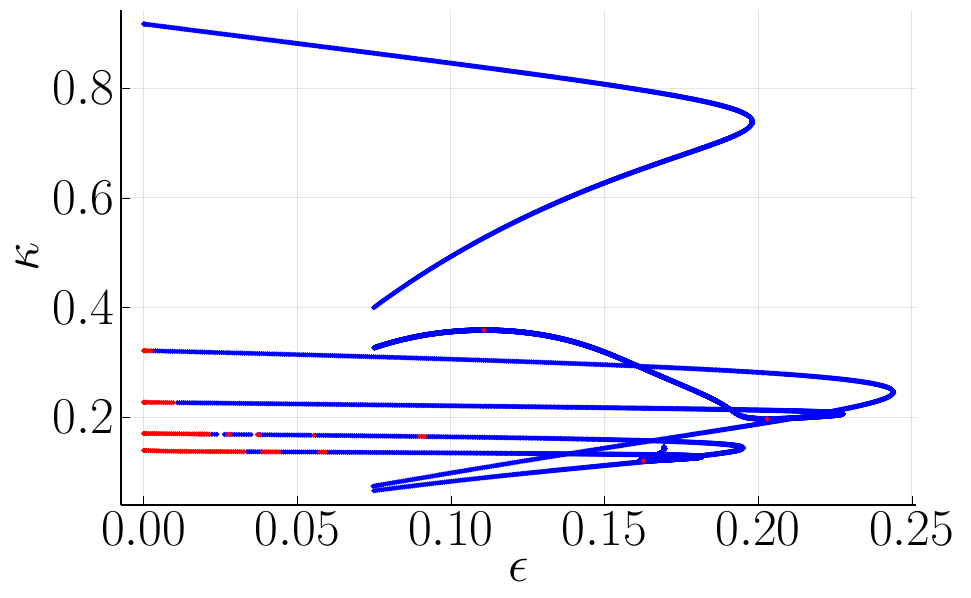}
    \caption{Fixing \(\kappa\).}
    \label{fig:CGL-pointwise-d=3-kappa}
  \end{subfigure}
  \caption{Pointwise results for \textbf{Case II}. Blue indicates
    points that have been successfully verified, and red indicates
    points that have not been verified.}
  \label{fig:CGL-pointwise-d=3}
\end{figure}

For each point on the branches we attempt to prove the existence of a
nearby solution. The results of this pointwise verification are
presented in Figure~\ref{fig:CGL-pointwise-d=3}, where blue indicates
a successfully verified point, and red a point that was not
successfully verified. In Figure~\ref{fig:CGL-pointwise-d=3-epsilon}
the verification is done by fixing \(\epsilon\), and searching for a
solution in \(\mu\), \(\gamma\) and \(\kappa\). Whereas in
Figure~\ref{fig:CGL-pointwise-d=3-kappa}, \(\kappa\) is fixed, and the
solution is searched for in \(\mu\), \(\gamma\) and \(\epsilon\).

\begin{figure}
  \centering
  \begin{subfigure}[t]{0.45\textwidth}
    \includegraphics[width=\textwidth]{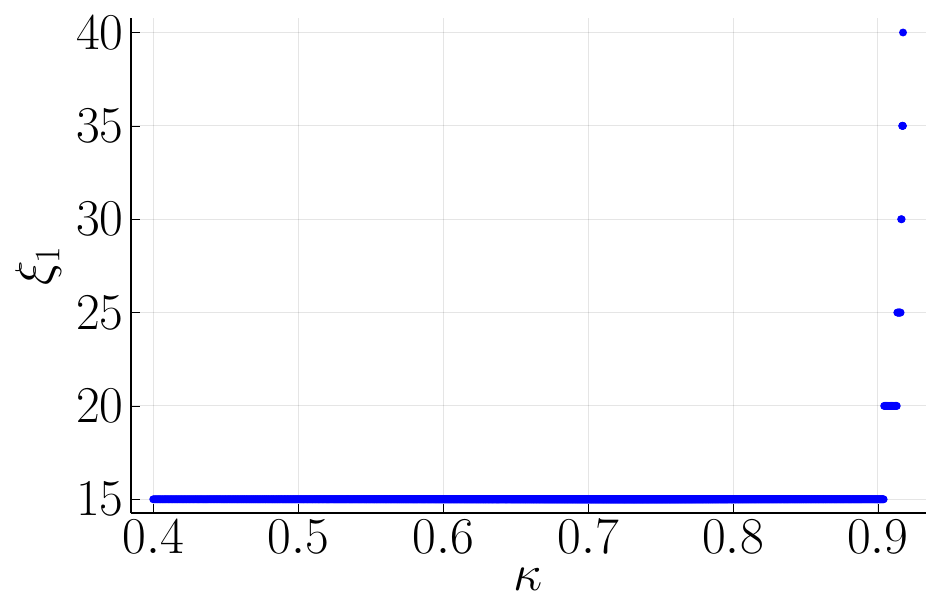}
    \caption{\(j = 1\)}
  \end{subfigure}
  \hspace{0.05\textwidth}
  \begin{subfigure}[t]{0.45\textwidth}
    \includegraphics[width=\textwidth]{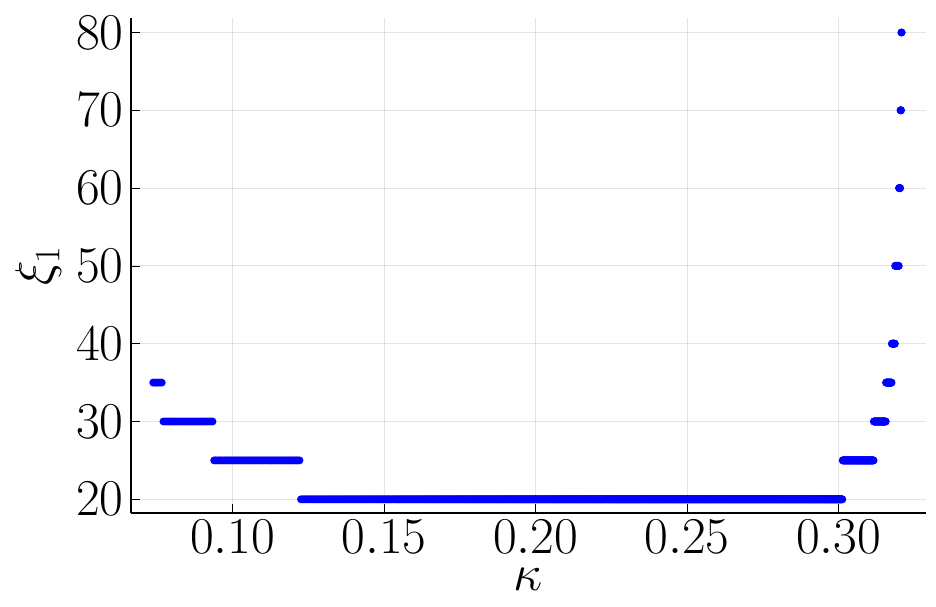}
    \caption{\(j = 2\)}
  \end{subfigure}
  \caption{Minimum \(\xi_{1}\) value for successful verification of
    points along the first two branches in \textbf{Case II}. Note that
    the \(x\)-axis is given by the \(\kappa\)-value, and hence the
    beginning of the curves are on the right-hand side. Larger
    \(\xi_{1}\) values are required in the beginning of the branches.}
  \label{fig:CGL-pointwise-xi-1-d=3}
\end{figure}

The figures show that, as expected, fixing \(\epsilon\) gives issues
near the turning points in \(\epsilon\), and fixing \(\kappa\) gives
issues near the turning points in \(\kappa\) (which only exist for
\(j = 3\) and \(j = 5\)). They also indicate that the beginning of the
branches is more difficult to verify than the later parts. We can get
an understanding of why this is the case by looking at the \(\xi_{1}\)
values used in the computations. Compared to the results in
Section~\ref{sec:nls} and~\ref{sec:branches} we do in this case not
fix the value of \(\xi_{1}\) a priori, but instead choose it
dynamically. We start with a low value and then successively increase
it until either the verification succeeds or we reach some
predetermined limit. For these results we try, in order, the
\(\xi_{1}\) values \(15, 20, 25, 30, 35, 40, 50, 60, 70\) and \(80\).
Figure~\ref{fig:CGL-pointwise-xi-1-d=3} shows the first value of
\(\xi_{1}\) at which the verification succeeded along the first two
branches. Points near the beginning (and to some extent the end) of
the branches require larger values of \(\xi_{1}\), whereas points in
the middle can be verified with much smaller values.

\begin{figure}
  \centering
  \includegraphics[width=0.45\textwidth]{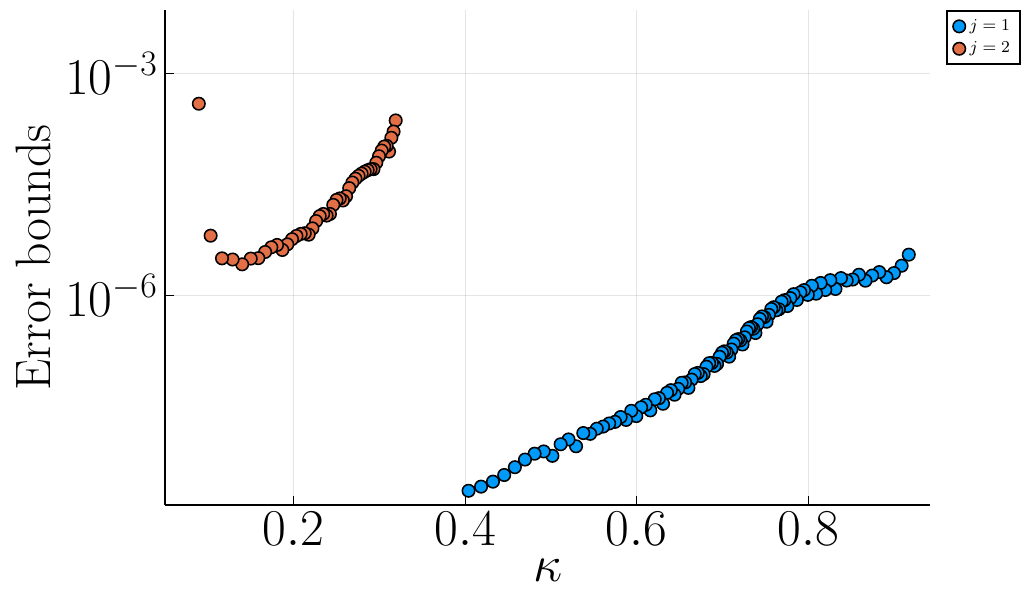}
  \caption{Error bound when enclosing \(\real(Q'_{\infty}(\xi_{1}))\)
    along the branches in \textbf{Case II}. For \(j = 1\) the value
    \(\xi_{1} = 15\) was used, and for \(j = 2\) the value
    \(\xi_{1} = 25\).}
  \label{fig:CGL-pointwise-error-bounds}
\end{figure}

The primary reason that the beginning of the branches requires larger
values of \(\xi_{1}\) seems to be that the bounds used for the
solution at infinity, \(Q_{\infty}\), are worst in the beginning of
the branches, see Figure~\ref{fig:CGL-pointwise-error-bounds}. Further
down the branches we can hence use a lower value of \(\xi_{1}\), and
still get good enough bounds to verify the solution. For the
non-verified points near the beginning of the branches one could
attempt to use even larger values for \(\xi_{1}\), this has however
only been partially successful in our case. For the second branch it
is indeed possible to verify the beginning of the branch. However, as
seen in Table~\ref{tab:results-nls-d-3}, this does require using a
much larger \(\xi_{1}\) value. For the third branch we have been
unable to verify the beginning, even with very large values for
\(\xi_{1}\). The issue being that a larger \(\xi_{1}\) puts more
strain on the integrator used for computing the solution from zero,
\(Q_{0}\), and at some point those errors start to dominate. To be
able to verify the remaining points it seems like the best approach
would be to extend the results in Section~\ref{sec:solution-infinity}
to improve on the bounds for \(Q_{\infty}\), which should allow for
using smaller values for \(\xi_{1}\).

\section{Implementation details}
\label{sec:implementation-details}
The full code on which the computer-assisted parts of the proofs are
based, as well as notebooks presenting the results and figures, is
available in the repository~\cite{CGL.jl}. Here we give an overview of
the methods used and their implementations. More details about the
structure of the code and exactly how the results are generated can be
found in that repository. The majority of the code is implemented in
Julia~\cite{Julia-2017}, except for the parts related to
CAPD~\cite{Kapela2021}, which are implemented in C++.

For the NLS equation, the computations can, at a high level, be split
into five parts:
\begin{enumerate}
\item Computing the (non-rigorous) numerical approximations;
\item Enclosing \(Q_{\infty}\) and its derivatives;
\item Enclosing \(Q_{0}\) and its derivatives;
\item Enclosing roots of \(G\);
\item Counting critical points of \(|Q|\).
\end{enumerate}
For the CGL equation, each individual segment of the branch is handled
using the same approach as for the NLS equation, but we also need to
compute numerical approximations of the branches and join all the
individual segments together.

For the non-rigorous parts, namely the numerical approximations, we
make use of classical numerical methods. For this, there are many
high-quality scientific libraries available in Julia. The most
important ones for our purposes are state-of-the-art methods for
solving ODEs with
\texttt{DifferentialEquations.jl}~\cite{Rackauckas2017}, nonlinear
solvers in \texttt{NonlinearSolve.jl}~\cite{Pal2024}, and methods for
automatic continuation of branches implemented in
\texttt{BifurcationKit.jl}~\cite{BifurcationKit} (this package
automatically handles bifurcation scenarios such as the saddle nodes
appearing in our problem).

The rigorous parts, on which the computer-assisted results are based,
make use of two different libraries, the CAPD~\cite{Kapela2021}
library and the Arb~\cite{Johansson2017arb}/FLINT~\cite{Flint}
library\footnote{In 2023 Arb was merged with the FLINT library.}. The
CAPD library is used for the rigorous numerical integration of the
ODE~\eqref{eq:problem-zero}, thus obtaining rigorous enclosures of
\(Q_{0}\) and its Jacobian, and the Arb library is used for everything
else. For the plotting of the interval enclosures, the package
\texttt{IntervalArithmetic.jl}~\cite{IntervalArithmetic.jl} is also
used. The Arb library is used through the Julia wrapper
\texttt{Arblib.jl}~\cite{Arblib.jl}, with many of the basic interval
arithmetic algorithms, such as isolating roots or enclosing maximum
values of single-variable functions, implemented in a separate
package, \texttt{ArbExtras.jl}~\cite{ArbExtras.jl}.

The enclosure of \(Q_{\infty}\) (as well as its Jacobian) is computed
following these steps:
\begin{enumerate}
\item Enclosures of the involved special functions are computed at
  \(\xi_{1}\). This includes \(P(\xi_{1})\), \(E(\xi_{1})\),
  \(J_{P}(\xi_{1})\), \(J_{E}(\xi_{1})\), \(D(\xi_{1})\) and
  \(H(\xi_{1})\), as well as derivatives with respect to \(\xi\),
  \(\kappa\) and \(\epsilon\). The computations are based on existing
  implementations in Arb for computing hypergeometric
  functions~\cite{Johansson2019}, as well as the expansion given in
  Lemma~\ref{lemma:U-a}.
\item Bounds for all the coefficients given in
  Lemma~\ref{lemma:bounds-list} are computed based on
  Lemmas~\ref{lemma:U} and~\ref{lemma:U-a} as well as the lemmas in
  Appendix~\ref{sec:function-bounds-proofs}.
\item An initial bound for \(\|Q\|_{\normv}\) is computed by bounding
  \(C_{T}\) from Lemma~\ref{lemma:fixed-point-bounds} and then finding
  \(\rho\) such that Inequalities~\eqref{eq:T-ineq-1}
  and~\eqref{eq:T-ineq-2} are satisfied. This gives the bound
  \(\|Q\|_{\normv} \leq \rho\) from the fixed-point argument in
  Proposition~\ref{prop:fixed-point}.
\item The bound for \(\|Q\|_{\normv}\) is then used to compute bounds
  for the norms of derivatives of \(Q\) with respect to \(\xi\),
  \(\real\gamma\), \(\imag\gamma\), \(\kappa\) and \(\epsilon\). This
  is done using the lemmas in Section~\ref{sec:bounds-on-norms}, which
  in turn rely on bounds for \(I_{P}\) and \(I_{E}\) from lemmas in
  Section~\ref{sec:I_P-I_E-direct-bounds}. Since the bounds for higher
  derivatives generally make use of bounds for the lower derivatives,
  this has to be done in an iterative fashion.
\item The enclosures for \(Q(\xi_{1})\) and associated derivatives are
  computed following the approach in
  Section~\ref{sec:enclosures-xi-1}.
\end{enumerate}

For the computation of \(Q_{0}\), the first step is to compute an
enclosure of \(Q_{0}(\xi_{0})\) for some small value \(\xi_{0}\) using
a Taylor expansion at zero. The code starts by trying
\(\xi_{0} = \codenumber{0.01}\), but tries successively smaller values
in case that fails. The Taylor expansion is computed using the
recurrence relation for the coefficients given in
Section~\ref{sec:solution-zero-removable-singularity}, which, combined
with the bounds for the remainder from Lemma~\ref{lemma:tail-bound},
allows us to enclose \(Q_{0}(\xi_{0})\). The solution is then
integrated along the interval \([\xi_{0}, \xi_{1}]\) using the CAPD
library: this makes use of high-degree Taylor expansions (degree
\codenumber{20} in our case) together with accurate representations of
the enclosures that are necessary to efficiently evolve the enclosures
in time. For the case \(d = 1\), there is no removable singularity at
\(\xi = 0\), and \(\xi_{0}\) is taken to be zero.

As discussed in Section~\ref{sec:detailed-example-nls}, the enclosures
of the roots of \(G\) are computed using the \emph{Krawczyk interval
  Newton method}. This involves computing, for a given set \(X\), an
enclosure of the set
\begin{equation*}
  \midint(X) - YG(\midint(X)) + (I - YJ_{G}(X))(X - \midint(X)),
\end{equation*}
where enclosures of \(G\) and \(J_{G}\) are computed through
\(Q_{\infty}\) and \(Q_{0}\) as described above. The matrix \(Y\) is a
preconditioner and is computed by taking the midpoint of the
enclosure for \(J_{G}(X)\) and inverting the resulting matrix using
regular non-rigorous floating-point arithmetic. To find a set \(X\)
for which the Newton condition~\eqref{eq:newton-condition} is
satisfied, we start with an initial approximation, \(X_{0}\), and then
compute
\begin{equation*}
  \tilde{X}_{1} = \midint(X_{0}) - YG(\midint(X_{0})) + (I - YJ_{G}(X_{0}))(X_{0} - \midint(X_{0})).
\end{equation*}
The set \(\tilde{X}_{1}\) is then expanded slightly by increasing its
diameter by a small factor, and we call the expanded set \(X_{1}\).
This process is then iterated, and we stop if at any step we have
\(\tilde{X}_{i + 1} \subseteq \operatorname{int}(X_{i})\), in which
case we set \(X = X_{i}\). Initially, \(X_{0}\) represents a thin
interval (a point), but due to the nonzero error bounds introduced in
the computations, the later iterations will always be non-thin
intervals. If the verification has not succeeded after
\codenumber{10} steps, the process stops and reports failure. The
process also stops early if at any point the enclosure of either
\(G(\midint(X_{i}))\) or \(J_{G}(X_{i})\) is non-finite, or if the
computation of \(Y\) fails. On success, the set \(X\) returned by this
approach generally just barely satisfies the Newton condition. This
means that the existence and uniqueness sets are about the same size.
To expand the set of uniqueness, we can successively widen \(X\) until
\(J_{G}\) stops being invertible; this generally gives us a much
larger set for which we have uniqueness.

The computation of critical points closely follows the approach
detailed in Section~\ref{sec:detailed-example-nls}, and there is not
much more to add. The bounds on \(C_{p_{\mathrm{mon}}}\) and
\(C_{R_{\mathrm{mon}}}\) are computed in the same way as those used
for \(Q_{\infty}\). The enclosures on the interval \([0, \xi_{2}]\)
are computed using a combination of a Taylor expansion at zero and
CAPD, as described for \(Q_{0}\) above.

\section*{Acknowledgments}
The authors are grateful for helpful comments from Denis Gaidashev and
Javier Gómez-Serrano. We are also thankful to Fredrik Johansson for
discussions about the evaluation of hypergeometric functions, and to
Daniel Wilczak for valuable input on the CAPD library. We thank Erik
Wahlén for helpful discussions, Svetlana Roudenko for providing
valuable insight into the NLS literature, as well as Vladimír Šverák
for fruitful discussions, in particular regarding the behavior of the
branches in \textbf{Case II}.

\noindent \textbf{Author contributions} All authors contributed
equally to the work.

\medskip

\noindent \textbf{Funding} Initial computations were enabled by
resources provided by the National Academic Infrastructure for
Supercomputing in Sweden (NAISS) and the Swedish National
Infrastructure for Computing (SNIC) at the PDC Center for High
Performance Computing, KTH Royal Institute of Technology, partially
funded by the Swedish Research Council through grant agreements no.
2022-06725 and no. 2018-05973. The final calculations made use of
resources provided by the Minnesota Supercomputing Institute (MSI) at
the University of Minnesota. J-Ll Figueras has been partially
supported by the VR-grant 2019-04591.

\medskip

\noindent \textbf{Data availability} The code on which the
computer-assisted parts of the proofs are based and notebooks
generating the results and figures are available at the following URL:
\url{https://github.com/Joel-Dahne/CGL.jl}.

\medskip

\noindent \textbf{Declaration of AI usage} AI, primarily Claude Opus 5
and earlier versions, has been used as part of proofreading the
manuscript and reviewing the code. The ideas and the writing of both
the manuscript and the code are entirely the work of the authors, who
take full responsibility for the content of both.

\appendix

\section{Auxiliary results on explicit asymptotic bounds}
\label{sec:function-bounds-proofs}
The primary purpose of this appendix is to give the explicit bounds on
which Lemma~\ref{lemma:bounds-list} is based. The functions under
consideration are given in terms of the confluent hypergeometric
function \(U\) and its derivatives, for which asymptotic bounds are
given in Lemmas~\ref{lemma:U} and~\ref{lemma:U-a}. The main procedure
in the below proofs is to explicitly write the functions in terms of
\(U\), and then apply the known bounds for \(U\).

In addition to the notation in Section~\ref{sec:solution-infinity} we
let
\begin{equation*}
  C_{U^{(m)}}(a, b, n, z_{1}) = C_{U}(a + m, b + m, n, z_{1})(a)_{m},
\end{equation*}
so that
\begin{equation}
  \label{eq:U-dz-bound}
  U^{(m)}(a, b, z) \leq C_{U^{(m)}}(a, b, n, z_{1})|z^{-a - m}|.
\end{equation}
We also let
\begin{equation*}
  p_{U,k}(a, b) = \frac{(a)_{k}(a - b + 1)_{k}}{k!}
\end{equation*}
denote the coefficients in the asymptotic expansion of \(U\).

Lemmas~\ref{lemma:U} and~\ref{lemma:U-a} have conditions on the
parameters \(a\), \(b\) and \(\xi\) for the bounds to be valid. In the
presentation below we do not make these conditions explicit. Instead,
whenever a constant such as \(C_{U}\) appears it is implicitly
understood that the resulting bound only holds if the parameters
satisfies the required conditions. That these conditions are actually
satisfied is checked whenever one of these constants are computed in
the code.

\begin{lemma}\label{lemma:C-P-E}
  Let \(\xi_{1} > 1\). For \(0 \leq m \leq 3\),
  \(n \in \mathbb{Z}_{\geq 0}\) and \(\xi \geq \xi_{1}\) we have
  \begin{equation*}
    |P^{(m)}(\xi)| \leq C_{P^{(m)}}\xi^{-\frac{1}{\sigma} - m},\quad
    |E^{(m)}(\xi)| \leq C_{E^{(m)}}e^{\real(c)\xi^{2}}\xi^{\frac{1}{\sigma} - d + m}.
  \end{equation*}
  Here the constants are given by
  \begin{align*}
    C_{P} &= |c^{-a}|C_{U}(a, b, n, c\xi_{1}^{2}),\\
    C_{P'}
    &= |2c^{-a}|C_{U'}(a, b, n, c\xi_{1}^{2}),\displaybreak[3] \\
    C_{P''}
    &= |2ac^{-a}|\Bigg(
      \sum_{k = 0}^{n - 1}
      \left|(2(a + 1)p_{U,k}(a + 2, b + 2) - p_{U,k}(a + 1, b + 1))(c\xi_{1}^{2})^{-k}\right|\\
    &\qquad+ (2|a + 1|C_{R_{U}}(a + 2, b + 2, n, c\xi_{1}^{2}) + C_{R_{U}}(a + 1, b + 1, n, c\xi_{1}^{2}))|c\xi_{1}^{2}|^{-n}
      \Bigg),\displaybreak[3] \\
    C_{P'''}
    &= |4a(a + 1)c^{-a}|\Bigg(
      \sum_{k = 0}^{n - 1}
      \left|(-2(a + 2)p_{U,k}(a + 3, b + 3) + 3p_{U,k}(a + 2, b + 2))(c\xi_{1}^{2})^{-k}\right|\\
    &\qquad+ (2|a + 2|C_{R_{U}}(a + 3, b + 3, n, c\xi_{1}^{2}) + 3C_{R_{U}}(a + 2, b + 2, n, c\xi_{1}^{2}))|c\xi_{1}^{2}|^{-n}
      \Bigg)
  \end{align*}
  and
  \begin{align*}
    C_{E} &= |(-c)^{a - b}|C_{U}(b - a, b, n, -c\xi_{1}^{2}),\\
    C_{E'}
    &= |2c(-c)^{a - b}|C_{U}(b - a, b, n, -c\xi_{1}^{2})
      + |2c(-c)^{a - b - 1}|C_{U'}(b - a, b, n, -c\xi_{1}^{2})\xi_{1}^{-2},\\
    C_{E''}
    &= |(4c^{2} + c\xi_{1}^{-2})(-c)^{a - b}|C_{U}(b - a, b, n, -c\xi_{1}^{2})\\
    &\qquad+ |(8c^{2} + 2c\xi_{1}^{-2})(-c)^{a - b - 1}|C_{U'}(b - a, b, n, -c\xi_{1}^{2})\xi_{1}^{-2}\\
    &\qquad+ |4c^{2}(-c)^{a - b - 2}|C_{U''}(b - a, b, n, -c\xi_{1}^{2})\xi_{1}^{-4},\\
    C_{E'''}
    &= |(8c^{3} + 12c^{2}\xi_{1}^{-2})(-c)^{a - b}|C_{U}(b - a, b, n, -c\xi_{1}^{2})\\
    &\qquad+ |(24c^{3} + 24c^{2}\xi_{1}^{-2})(-c)^{a - b - 1}|C_{U'}(b - a, b, n, -c\xi_{1}^{2})\xi_{1}^{-2}\\
    &\qquad+ |(24c^{3} + 12c^{2}\xi_{1}^{-2})(-c)^{a - b - 2}|C_{U''}(b - a, b, n, -c\xi_{1}^{2})\xi_{1}^{-4}\\
    &\qquad+ |8c^{3}(-c)^{a - b - 3}|C_{U'''}(b - a, b, n, -c\xi_{1}^{2})\xi_{1}^{-6}.
  \end{align*}
\end{lemma}

\begin{proof}
  Writing \(P\) and \(E\) in terms of \(U\) we get
  \begin{align*}
    P(\xi)
    &= U(a, b, c\xi^{2}),\\
    P'(\xi)
    &= 2c\xi U'(a, b, c\xi^{2}),\\
    P''(\xi)
    &= 4c^{2}\xi^{2}U''(a, b, c\xi^{2}) + 2cU'(a, b, c\xi^{2}),\\
    P'''(\xi)
    &= 8c^{3}\xi^{3}U'''(a, b, c\xi^{2}) + 12c^{2}\xi U''(a, b, c\xi^{2})
  \end{align*}
  and
  \begin{align*}
    E(\xi)
    &= e^{c\xi^{2}}U(b - a, b, -c\xi^{2}),\\
    E'(\xi)
    &= 2ce^{c\xi^{2}}\xi(U(b - a, b, -c\xi^{2}) - U'(b - a, b, -c\xi^{2})),\\
    E''(\xi)
    &= e^{c\xi^{2}}\Big(
      (4c^{2}\xi^{2} + 2c)U(b - a, b, -c\xi^{2})
      - (8c^{2}\xi^{2} + 2c)U'(b - a, b, -c\xi^{2})\\
    &\qquad\qquad+ 4c^{2}\xi^{2}U''(b - a, b, -c\xi^{2})
      \Big),\\
    E'''(\xi)
    &= e^{c\xi^{2}}\Big(
      (8c^{3}\xi^{3} + 12c^{2}\xi)U(b - a, b, -c\xi^{2})
      - (24c^{3}\xi^{3} + 24c^{2}\xi)U'(b - a, b, -c\xi^{2})\\
    &\qquad\qquad+ (24c^{3}\xi^{3} + 12c^{2}\xi)U''(b - a, b, -c\xi^{2})
      - 8c^{3}\xi^{3}U'''(b - a, b, -c\xi^{2})
      \Big).
  \end{align*}

  For \(P\), \(P'\) and all derivatives of \(E\), it is enough to
  bound everything termwise, using the bounds
  \begin{equation*}
    |U(a, b, c\xi^{2})| \leq |c^{-a}|C_{U}(a, b, n, c\xi_{1}^{2})\xi^{-2\real(a)}
  \end{equation*}
  and
  \begin{equation*}
    |U(b - a, b, -c\xi^{2})| \leq |(-c)^{a - b}|C_{U}(b - a, b, n, -c\xi^{2})\xi^{2\real(a - b)},
  \end{equation*}
  and the bound~\eqref{eq:U-dz-bound} for the higher derivatives.

  For \(P''\) and \(P'''\), a termwise bound gives large
  overestimations. All terms are of similar magnitude, and it is
  therefore important to take into account their cancellations. We
  instead expand them using the asymptotic expansions and collect
  terms with the same exponents for \(\xi\). The resulting expansions
  are then bounded termwise.
\end{proof}

\begin{lemma}\label{lemma:C-P-E-dkappa}
  Let \(\xi_{1} > 1\). For \(n \in \mathbb{Z}_{\geq 0}\) and
  \(\xi \geq \xi_{1}\) we have
  \begin{equation*}
    |P_{\kappa}^{(m)}(\xi)| \leq C_{P^{(m)},\kappa}\log(\xi)\xi^{-\frac{1}{\sigma} - m},\quad
    |E_{\kappa}^{(m)}(\xi)| \leq C_{E^{(m)},\kappa}e^{\real(c)\xi^{2}}\xi^{\frac{1}{\sigma} - d + 2 + m},
  \end{equation*}
  with \(0 \leq m \leq 2\) for \(P_{\kappa}\) and \(0 \leq m \leq 1\)
  for \(E_{\kappa}\). The constants are given by
  \begin{align*}
    C_{P,\kappa}
    &= |c^{-a - 1}c_{\kappa}|C_{U'}(a, b, n, c\xi_{1}^{2})\log(\xi_{1})^{-1}
      + |c^{-a}a_{\kappa}|C_{U,a}(a, b, n, c\xi_{1}^{2})(2 + |\log(c)|\log(\xi_{1})^{-1}),\\
    C_{P',\kappa}
    &= |2c^{-a - 1}c_{\kappa}|C_{U'}(a, b, n, c\xi_{1}^{2})\log(\xi_{1})^{-1}
      + |2c^{-a - 1}c_{\kappa}|C_{U''}(a, b, n, c\xi_{1}^{2})\log(\xi_{1})^{-1}\\
    &\qquad+ |2c^{-a}a_{\kappa}a|C_{U,a}(a + 1, b + 1, n, c\xi_{1}^{2})(2 + |\log(c)|\log(\xi_{1})^{-1})\\
    &\qquad+ |2c^{-a}a_{\kappa}|C_{U}(a + 1, b + 1, n, c\xi_{1}^{2})\log(\xi_{1})^{-1},\\
    C_{P'',\kappa}
    &= |2c^{-a - 1}c_{\kappa}|C_{U'}(a, b, n, c\xi_{1}^{2})\log(\xi_{1})^{-1}
      + |10c^{-a - 1}c_{\kappa}|C_{U''}(a, b, n, c\xi_{1}^{2})\log(\xi_{1})^{-1}\\
    &\qquad+ |4c^{-a - 1}c_{\kappa}|C_{U'''}(a, b, n, c\xi_{1}^{2})\log(\xi_{1})^{-1}
      + |2c^{-a}a_{\kappa}|C_{U}(a + 1, b + 1, n, c\xi_{1}^{2})\log(\xi_{1})^{-1}\\
    &\qquad+ |2c^{-a}a_{\kappa}a|C_{U,a}(a + 1, b + 1, n, c\xi_{1}^{2})(2 + |\log(c)|\log(\xi_{1})^{-1})\\
    &\qquad+ |4c^{-a}a_{\kappa}(2a + 1)|C_{U}(a + 2, b + 2, n, c\xi_{1}^{2})\log(\xi_{1})^{-1}\\
    &\qquad+ |4c^{-a}a_{\kappa}a(a + 1)|C_{U,a}(a + 2, b + 2, n, c\xi_{1}^{2})(2 + |\log(c)|\log(\xi_{1})^{-1}),\displaybreak[3] \\
    C_{E,\kappa}
    &= |(-c)^{a - b}c_{\kappa}|C_{U}(b - a, b, n, -c\xi_{1}^{2})\\
    &\qquad+ |(-c)^{a - b}a_{\kappa}|C_{U,a}(b - a, b, n, -c\xi_{1}^{2})(2 + |\log(-c)|\log(\xi_{1})^{-1})\log(\xi_{1})\xi_{1}^{-2}\\
    &\qquad+ |(-c)^{a - b - 1}c_{\kappa}|C_{U'}(b - a, b, n, -c\xi_{1}^{2})\xi_{1}^{-2},\\
    C_{E',\kappa}
    &= 2|(-c)^{a - b}c_{\kappa}|C_{U}(b - a, b, n, -c\xi_{1}^{2})(|c| + \xi_{1}^{-2})\\
    &\qquad+ 2|(-c)^{a - b - 1}c_{\kappa}|C_{U'}(b - a, b, n, -c\xi_{1}^{2})(|2c| + \xi_{1}^{-2})\xi_{1}^{-2}\\
    &\qquad+ 2|(-c)^{a - b}a_{\kappa}c|C_{U,a}(b - a, b, n, -c\xi_{1}^{2})(2 + |\log(-c)|\log(\xi_{1})^{-1})\log(\xi_{1})\xi_{1}^{-2}\\
    &\qquad+ 2|(-c)^{a - b - 2}c_{\kappa}c|C_{U''}(b - a, b, n, -c\xi_{1}^{2})\xi_{1}^{-2}\\
    &\qquad+ 2|(-c)^{a - b - 1}a_{\kappa}c|C_{U}(b - a + 1, b + 1, n, c\xi_{1}^{2})\xi_{1}^{-4}\\
    &\qquad+ 2|(b - a)(-c)^{a - b - 1}a_{\kappa}c|C_{U,a}(b - a + 1, b + 1, n, c\xi_{1}^{2})(2 + |\log(-c)|\log(\xi_{1})^{-1})\log(\xi_{1})\xi_{1}^{-4}.
  \end{align*}
\end{lemma}

\begin{proof}
  Writing the functions in terms of \(U\) and using that
  \begin{align*}
    U_{a}'(a, b, c\xi^{2}) &= -U(a + 1, b + 1, c\xi^{2}) - aU_{a}(a + 1, b + 1, c\xi^{2}),\\
    U_{a}''(a, b, c\xi^{2}) &= (2a + 1)U(a + 2, b + 2, c\xi^{2}) + a(a + 1)U_{a}(a + 2, b + 2, c\xi^{2}),
  \end{align*}
  we get
  \begin{align*}
    P_{\kappa}(\xi)
    &= c_{\kappa}U'(a, b, c\xi^{2})\xi^{2} + a_{\kappa}U_{a}(a, b, c\xi^{2}),\\
    P_{\kappa}'(\xi)
    &= 2c_{\kappa}U'(a, b, c\xi^{2})\xi
      + 2c_{\kappa}cU''(a, b, c\xi^{2})\xi^{3}
      - 2a_{\kappa}cU(a + 1, b + 1, c\xi^{2})\xi\\
    &\qquad- 2a_{\kappa}acU_{a}(a + 1, b + 1, c\xi^{2})\xi,\\
    P_{\kappa}''(\xi)
    &= 2c_{\kappa}U'(a, b, c\xi^{2})\xi^{2}
      + 10c_{\kappa}cU''(a, b, c\xi^{2})\xi^{2}
      + 4c_{\kappa}c^{2}U'''(a, b, c\xi^{2})\xi^{4}\\
    &\qquad- 2a_{\kappa}cU(a + 1, b + 1, c\xi^{2})
      - 2aa_{\kappa}cU_{a}(a + 1, b + 1, c\xi^{2})\\
    &\qquad+ 4a_{\kappa}(2a + 1)c^{2}U(a + 2, b + 2, c\xi^{2})
      + 4a_{\kappa}a(a + 1)c^{2}U_{a}(a + 2, b + 2, c\xi^{2}),\displaybreak[3]\\
    E_{\kappa}(\xi)
    &= e^{c\xi^{2}}(
      c_{\kappa}U(b - a, b, -c\xi^{2})\xi^{2}
      - a_{\kappa}U_{a}(b - a, b, -c\xi^{2})
      - c_{\kappa}U'(b - a, b, -c\xi^{2})\xi^{2}
      ),\\
    E_{\kappa}'(\xi)
    &= 2e^{c\xi^{2}}\xi\Big(
      c_{\kappa}(1 + c\xi^{2})U(b - a, b, -c\xi^{2})
      - c_{\kappa}(1 + 2c\xi^{2})U'(b - a, b, -c\xi^{2})\\
    &\qquad- a_{\kappa}cU_{a}(b - a, b, -c\xi^{2})
      + cc_{\kappa}U''(b - a, b, -c\xi^{2})\xi^{2}\\
    &\qquad- a_{\kappa}cU'(b - a + 1, b + 1, -c\xi^{2})
      - (b - a)a_{\kappa}cU_{a}(b - a + 1, b + 1, -c\xi^{2})
      \Big).
  \end{align*}
  We compute the bounds termwise. For \(E_\kappa\) and \(E_\kappa'\)
  this gives a very tight upper bound. For \(P_\kappa\) the termwise
  bound gives an overestimation by about a factor \(1.5\) in our
  cases, and for \(P_\kappa'\) the overestimation is by about a factor
  \(2\) and for \(P_\kappa''\) about \(5\). This is, however, good
  enough for our purposes.
\end{proof}

\begin{lemma}\label{lemma:C-P-E-depsilon}
  Let \(\xi_{1} > 1\). For \(n \in \mathbb{Z}_{\geq 0}\) and
  \(\xi \geq \xi_{1}\) we have
  \begin{equation*}
    |P_{\epsilon}^{(m)}(\xi)| \leq C_{P^{(m)},\epsilon}\xi^{-\frac{1}{\sigma} - m},\quad
    |E_{\epsilon}^{(m)}(\xi)| \leq C_{E^{(m)},\epsilon}e^{\real(c)\xi^{2}}\xi^{\frac{1}{\sigma} - d + 2 + m},
  \end{equation*}
  with \(0 \leq m \leq 2\) for \(P_{\epsilon}\) and
  \(0 \leq m \leq 1\) for \(E_{\epsilon}\). The constants are given by
  \begin{align*}
    C_{P,\epsilon}
    &= |c^{-a - 1}c_{\epsilon}|C_{U'}(a, b, n, c\xi_{1}^{2}),\\
    C_{P',\epsilon}
    &= |2c^{-a - 1}c_{\epsilon}|C_{U'}(a, b, n, c\xi_{1}^{2})
      + |2c^{-a - 1}c_{\epsilon}|C_{U''}(a, b, n, c\xi_{1}^{2}),\\
    C_{P'',\epsilon}
    &= |2c^{-a - 1}c_{\epsilon}|C_{U'}(a, b, n, c\xi_{1}^{2})
      + |10c^{-a - 1}c_{\epsilon}|C_{U''}(a, b, n, c\xi_{1}^{2})
      + |4c^{-a - 1}c_{\epsilon}|C_{U'''}(a, b, n, c\xi_{1}^{2}),\displaybreak[3]\\
    C_{E,\epsilon}
    &= |(-c)^{a - b}c_{\epsilon}|C_{U}(b - a, b, n, -c\xi_{1}^{2})
      + |(-c)^{a - b - 1}c_{\epsilon}|C_{U'}(b - a, b, n, -c\xi_{1}^{2})\xi_{1}^{-2},\\
    C_{E',\epsilon}
    &= 2|(-c)^{a - b}c_{\epsilon}|C_{U}(b - a, b, n, -c\xi_{1}^{2})(|c| + \xi_{1}^{-2})\\
    &\qquad+ 2|(-c)^{a - b - 1}c_{\epsilon}|C_{U'}(b - a, b, n, -c\xi_{1}^{2})(|2c| + \xi_{1}^{-2})\xi_{1}^{-2}\\
    &\qquad+ 2|(-c)^{a - b - 2}c_{\epsilon}c|C_{U''}(b - a, b, n, -c\xi_{1}^{2})\xi_{1}^{-2},\\
  \end{align*}
\end{lemma}

\begin{proof}
  Follows in the same way as Lemma~\ref{lemma:C-P-E-dkappa}, using that
  \(a_{\epsilon} = 0\).
\end{proof}

\begin{lemma}\label{lemma:C-J-P-J-E}
  Let \(\xi_{1} > 1\). For \(n \in \mathbb{Z}_{\geq 0}\) and
  \(\xi \geq \xi_{1}\) we have
  \begin{equation*}
    |J_{P}(\xi)| \leq C_{J_{P}}e^{-\real(c)\xi^{2}}\xi^{-\frac{1}{\sigma} + d - 1},\quad
    |J_{E}(\xi)| \leq C_{J_{E}}\xi^{\frac{1}{\sigma} - 1},
  \end{equation*}
  with the constants given by
  \begin{align*}
    C_{J_{P}} &= |B_{W}|C_{P},\\
    C_{J_{E}} &= |B_{W}|C_{E}.
  \end{align*}
\end{lemma}

\begin{proof}
  The bounds follow directly from Lemma~\ref{lemma:C-P-E}.
\end{proof}

\begin{lemma}\label{lemma:C-J-P-kappa-J-E-kappa}
  Let \(\xi_{1} > 1\). For \(n \in \mathbb{Z}_{\geq 0}\) and
  \(\xi \geq \xi_{1}\) we have
  \begin{align*}
    |J_{P,\kappa}(\xi)| &\leq C_{J_{P},\kappa}e^{-\real(c)\xi^{2}}\xi^{-\frac{1}{\sigma} + d + 1},\\
    |J_{E,\kappa}(\xi)| &\leq C_{J_{E},\kappa}\log(\xi)\xi^{\frac{1}{\sigma} - 1},\\
    |J_{P,\epsilon}(\xi)| &\leq C_{J_{P},\epsilon}e^{-\real(c)\xi^{2}}\xi^{-\frac{1}{\sigma} + d + 1},\\
    |J_{E,\epsilon}(\xi)| &\leq C_{J_{E},\epsilon}\xi^{\frac{1}{\sigma} - 1},
  \end{align*}
  with the constants given by
  \begin{align*}
    C_{J_{P},\kappa}
    &= C_{P}(
      |c_{\kappa}B_{W}| + |B_{W,\kappa}|\xi_{1}^{-2}
      )
      + C_{P,\kappa}|B_{W}|\log(\xi_{1})\xi_{1}^{-2}
      ,\\
    C_{J_{E},\kappa}
    &= C_{U,a}(b - a, b, n, -c\xi_{1}^{2})|B_{W}(-c)^{a - b}a_{\kappa}|(2 + |\log(-c)|\log(\xi_{1})^{-1})\\
    &\qquad+ C_{U}(b - a, b, n, -c\xi_{1}^{2})|B_{W,\kappa}(-c)^{a - b}|\log(\xi_{1})^{-1}\\
    &\qquad+ C_{U'}(b - a, b, n, -c\xi_{1}^{2})|B_{W,\kappa}(-c)^{a - b - 1}c_{\kappa}|\log(\xi_{1})^{-1},\\
    C_{J_{P},\epsilon}
    &= C_{P}(
      |c_{\epsilon}B_{W}| + |B_{W,\epsilon}|\xi_{1}^{-2}
      )
      + C_{P,\epsilon}|B_{W}|\xi_{1}^{-2}
      ,\\
    C_{J_{E},\epsilon}
    &= C_{U}(b - a, b, n, -c\xi_{1}^{2})|B_{W,\epsilon}(-c)^{a - b}|
      + C_{U'}(b - a, b, n, -c\xi_{1}^{2})|B_{W,\epsilon}(-c)^{a - b - 1}c_{\epsilon}|,\\
  \end{align*}
\end{lemma}

\begin{proof}
  We have
  \begin{align*}
    J_{P,\kappa}(\xi)
    &= e^{-c\xi^{2}}\xi^{d + 1}(
      P(\xi)(-c_{\kappa}B_{W} + B_{W,\kappa}\xi^{-2})
      + P_{\kappa}(\xi)B_{W}\xi^{-2}
      ),\\
    J_{E,\kappa}(\xi)
    &= U(b - a, b, -c\xi^{2})B_{W,\kappa}\xi^{d - 1}
      - U_{a}(b - a, b, -c\xi^{2})a_{\kappa}B_{W}\xi^{d - 1}
      - U'(b - a, b, -c\xi^{2})c_{\kappa}B_{W}\xi^{d + 1},\\
    J_{P,\epsilon}(\xi)
    &= e^{-c\xi^{2}}\xi^{d + 1}(
      P(\xi)(-c_{\epsilon}B_{W} + B_{W,\epsilon}\xi^{-2})
      + P_{\epsilon}(\xi)B_{W}\xi^{-2}
      ),\\
    J_{E,\epsilon}(\xi)
    &= U(b - a, b, -c\xi^{2})B_{W,\epsilon}\xi^{d - 1}
      - U'(b - a, b, -c\xi^{2})c_{\epsilon}B_{W}\xi^{d + 1}.
  \end{align*}
  They can then be bounded termwise.
\end{proof}

\begin{lemma}\label{lemma:C-D-H}
  Let \(\xi_{1} > 1\). For \(0 \leq m \leq 2\),
  \(n \in \mathbb{Z}_{\geq 0}\) and \(\xi \geq \xi_{1}\) we have
  \begin{equation*}
    |D^{(m)}(\xi)| \leq C_{D^{(m)}}\xi^{-\frac{1}{\sigma} - m},\quad
    |H^{(m)}(\xi)| \leq C_{H^{(m)}}\xi^{-\frac{1}{\sigma} - m},
  \end{equation*}
  with the constants given by
  \begin{align*}
    C_{D}
    &= |c_{\kappa}B_{W}|C_{P}
      + (
      |B_{W,\kappa}|C_{P}
      + |B_{W}|C_{P,\kappa}\log(\xi_{1})
      )\xi_{1}^{-2},\\
    C_{D'}
    &= |c_{\kappa}B_{W}|C_{P'}
      + \Big(
      |B_{W,\kappa}|C_{P'}
      + 2|B_{W,\kappa}|C_{P}\\
    &\qquad+ \big(
      |B_{W}|C_{P',\kappa}
      + 2|B_{W}|C_{P,\kappa}
      \big)\log(\xi_{1})
      \Big)\xi_{1}^{-2},\\
    C_{D''}
    &= |c_{\kappa}B_{W}|C_{P''}\\
    &\qquad+ \Big(
      |B_{W,\kappa}|C_{P''}
      + 4|B_{W,\kappa}|C_{P'}
      + 6|B_{W,\kappa}|C_{P}\\
    &\qquad\qquad+ \big(
      |B_{W}|C_{P'',\kappa}
      + 4|B_{W}|C_{P',\kappa}
      + 6|B_{W}|C_{P,\kappa}
      \big)\log(\xi_{1})
      \Big)\xi_{1}^{-2},\displaybreak[3]\\
    C_{H}
    &= |c_{\epsilon}B_{W}|C_{P}
      + (
      |B_{W,\epsilon}|C_{P}
      + |B_{W}|C_{P,\epsilon}
      )\xi_{1}^{-2},\\
    C_{H'}
    &= |c_{\epsilon}B_{W}|C_{P'}
      + \left(
      |B_{W,\epsilon}|C_{P'}
      + 2|B_{W,\epsilon}|C_{P}
      + |B_{W}|C_{P',\epsilon}
      + 2|B_{W}|C_{P,\epsilon}
      \right)\xi_{1}^{-2},\\
    C_{H''}
    &= |c_{\epsilon}B_{W}|C_{P''}\\
    &\qquad+ \Big(
      |B_{W,\epsilon}|C_{P''}
      + 4|B_{W,\epsilon}|C_{P'}
      + 6|B_{W,\epsilon}|C_{P}\\
    &\qquad\qquad+ |B_{W}|C_{P'',\epsilon}
      + 4|B_{W}|C_{P',\epsilon}
      + 6|B_{W}|C_{P,\epsilon}
      \Big)\xi_{1}^{-2}.
  \end{align*}
\end{lemma}

\begin{proof}
  For \(D\) we have
  \begin{align*}
    D'(\xi)
    &= -c_{\kappa}B_{W}P'(\xi)
      + B_{W,\kappa}P'(\xi)\xi^{-2}
      - 2B_{W,\kappa}P(\xi)\xi^{-3}
      + B_{W}P_{\kappa}'(\xi)\xi^{-2}
      - 2B_{W}P_{\kappa}(\xi)\xi^{-3},\\
    D''(\xi)
    &= -c_{\kappa}B_{W}P''(\xi)
      + B_{W,\kappa}P''(\xi)\xi^{-2}
      - 4B_{W,\kappa}P'(\xi)\xi^{-3}
      + 6B_{W,\kappa}P(\xi)\xi^{-4}\\
    &\qquad+ B_{W}P_{\kappa}''(\xi)\xi^{-2}
      - 4B_{W}P_{\kappa}'(\xi)\xi^{-3}
      + 6B_{W}P_{\kappa}(\xi)\xi^{-4}.
  \end{align*}
  The expressions for \(H\) are similar, with differentiation with
  respect to \(\kappa\) replaced with differentiation with respect to
  \(\epsilon\). All the bounds are then computed termwise.
\end{proof}

Finally, we have a bound that is used in the proof of
Lemma~\ref{lemma:fixed-point-bounds}.

\begin{lemma}
  For all \(z_1, z_2 \in \mathbb{C}\) and \(\sigma \geq 1\) we have
  \begin{equation*}
    \left||z_{1}|^{2\sigma}z_{1} - |z_{2}|^{2\sigma}z_{2}\right| \leq
    M_\sigma |z_{1} - z_{2}|(|z_{1}|^{2\sigma} + |z_{2}|^{2\sigma}).
  \end{equation*}
  For \(\sigma = 1\) we can take
  \begin{equation*}
    M_{1} = \frac{\sqrt{2}}{4 - 2\sqrt{2}} + 1.
  \end{equation*}
  In general we can take
  \begin{equation*}
    M_{\sigma} = \max\left(
      \max_{0 \leq t \leq t_{0}} g(t),
      g_{t_{0}}(t_{1}),
      \max_{t_{1} \leq t \leq 2} g(t),
      2
    \right),
  \end{equation*}
  for any \(0 \leq t_0 \leq 1 \leq t_1 \leq 2\). Here
  \begin{equation*}
    g(t) = \frac{1 - t^{2\sigma}}{(1 - t)(1 + t^{2\sigma})} + 1
    \quad\text{and}\quad
    g_{t_0}(t) = \frac{1 - t^{2\sigma}}{(1 - t)(1+t_0^{2\sigma})} + 1.
  \end{equation*}
  \label{lemma:M}
\end{lemma}

\begin{proof}
  Let
  \begin{equation*}
    f(z_{1}, z_{2}) = \frac{
      \left||z_{1}|^{2\sigma}z_{1} - |z_{2}|^{2\sigma}z_{2}\right|
    }{
      |z_{1} - z_{2}|(|z_{1}|^{2\sigma} + |z_{2}|^{2\sigma})
    }.
  \end{equation*}
  We want \(M_\sigma\) to be an upper bound for \(|f(z_1, z_2)|\). Writing
  \(z_2 = sz_1\), with \(s \in \mathbb{C}\) gives us
  \begin{equation*}
    f(z_1, z_2) = \frac{|1 - |s|^{2\sigma}s|}{|1 - s|(1 + |s|^{2\sigma})}.
  \end{equation*}
  An upper bound is given by
  \begin{equation*}
    \frac{|1 - |s|^{2\sigma}|}{|1 - s|(1 + |s|^{2\sigma})} + \frac{|s|^{2\sigma}}{1 + |s|^{2\sigma}}
    \leq \frac{|1 - |s|^{2\sigma}|}{|1 - s|(1 + |s|^{2\sigma})} + 1.
  \end{equation*}
  If we let \(t = |s|\) and use that \(|1 - s| \geq |1 - |s||\) we get
  that a uniform upper bound for \(|f(z_1, z_2)|\) is given by
  \(\sup_{t \geq 0} g(t)\), with \(g(t)\) as in the statement of the
  lemma. For \(\sigma = 1\) this simplifies to
  \begin{equation*}
    g(t) = \frac{1 + t}{1 + t^{2\sigma}} + 1,
  \end{equation*}
  for which the upper bound is attained at \(t = \sqrt{2} - 1\) and
  given by \(\sqrt{2} / (4 - 2\sqrt{2}) + 1\). This gives us \(M_1\).
  For \(\sigma \neq 1\) we start by noting that
  \begin{equation*}
    g(t)
    \leq \frac{1 + t^{2\sigma}}{|1 - t|(1 + t^{2\sigma})} + 1
    = \frac{1}{|1 - t|} + 1.
  \end{equation*}
  In particular, we have \(g(2) \leq 2\), and hence
  \begin{equation*}
    \sup_{t \geq 0} g(t) \leq \max\left(\max_{0 \leq t \leq 2} g(t), 2\right).
  \end{equation*}
  To handle the removable singularity at \(t = 1\) we let
  \(g_{t_{0}}(t)\) be as in the statement of the lemma and note that
  \(g(t) \leq g_{t_0}(t)\) for any \(t \geq t_0\). Furthermore, we
  have that \(g_{t_0}\) is increasing in \(t\). In particular, this
  gives us
  \begin{equation*}
    \max_{t_{0} \leq t \leq t_{1}} g(t) \leq g_{t_{0}}(t_{1})
  \end{equation*}
  for any \(0 \leq t_0 \leq 1 \leq t_1 \leq 2\). Combining the bounds
  for \(0 \leq t \leq t_0\), \(t_0 \leq t \leq t_1\), \(t_1 \leq t \leq
  2\) and \(t \geq 2\) gives us the result.

  What is left is to prove that $g_{t_0}$ is increasing in $t$. For
  doing so, observe that its derivative is
  \[
  g'_{t_0}(t)=\dfrac{(2\sigma-1)t^{2\sigma}-2\sigma t^{2\sigma-1}+1}{(1-t)^2(1+t_0^{2\sigma})}.
  \]
  Hence, $g_{t_0}$ is monotone increasing if its derivative's
  numerator is non-negative (it is equal to $1$ for $t=0$, equal to
  zero at $t=1$ and positive for large $t$). Finally, this numerator's
  derivative is equal to \(2\sigma(2\sigma-1)t^{2\sigma-2}(t-1)\),
  which is negative for $t\in[0,1)$, zero at $t=1$ and positive for
  $t\in(1, \infty)$.
\end{proof}

\printbibliography

\begin{tabular}{c@{\hspace{1cm}}c}
  \begin{tabular}[t]{l}
    \textbf{Joel Dahne} \\
    School of Mathematics \\
    University of Minnesota \\
    127 Vincent Hall, 206 Church St. SE \\
    Minneapolis, MN 55455, USA \\
    Email: jdahne@umn.edu
  \end{tabular}
  &
  \begin{tabular}[t]{l}
    \textbf{Jordi-Llu\'is Figueras} \\
    Department of Mathematics \\
    Uppsala University \\
    Box 480 \\
    751 06 Uppsala, Sweden \\
    Email: figueras@math.uu.se
  \end{tabular}
\end{tabular}

\end{document}